\documentclass{article}

\usepackage{vmargin}
\usepackage{epsfig}
\usepackage[francais]{babel}
\usepackage[latin1]{inputenc}
\usepackage[T1]{fontenc}
\usepackage[tbtags]{amsmath}
\usepackage{amsthm}
\usepackage{amssymb,bbm,mathrsfs}
\usepackage{array,multirow}
\usepackage{alltt}
\usepackage[all]{xy}
\usepackage{url}
\usepackage{times}
\usepackage{stmaryrd}
\usepackage{calc}
\usepackage{microtype}

\newtheorem{introtheo}{Théorème}
\newtheorem{theo}{Théorème}[section]
\newtheorem{lemme}[theo]{Lemme}
\newtheorem{prop}[theo]{Proposition}
\newtheorem{cor}[theo]{Corollaire}

\newtheorem{quest}[theo]{Question}
\theoremstyle{definition}
\newtheorem{deftn}[theo]{Définition}

\theoremstyle{remark}
\newtheorem{rem}[theo]{Remarque}

\def\leq{\leqslant}
\def\geq{\geqslant}

\def\N{\mathbb{N}}
\def\Z{\mathbb{Z}}
\def\Q{\mathbb{Q}}
\def\R{\mathbb{R}}
\def\C{\mathbb{C}}
\def\F{\mathbb{F}}

\def\a{\mathfrak{a}}
\def\d{\mathfrak{d}}
\def\m{\mathfrak{m}}
\def\t{\mathfrak{t}}
\def\Sk{\mathfrak{S}}
\def\O{\mathcal{O}}
\def\E{\mathcal{E}}

\def\calM{\mathcal{M}}

\def\calT{\mathcal{T}}
\def\calF{\mathcal{F}}

\def\calK{\mathcal{K}}

\def\calR{\mathcal{R}}
\def\calS{\mathcal{S}}

\def\frakM{\mathfrak{M}}
\def\frakN{\mathfrak{N}}

\def\eps{\varepsilon}
\def\ueps{\underline \eps}
\def\ueta{\eta}
\def\upi{\underline \pi}

\def\Tau{\mathcal{A}}

\def\coker{\text{\rm coker}}
\def\id{\text{\rm id}}
\def\GL{\text{\rm GL}}
\def\tr{\text{\rm tr}}

\def\Frac{\text{\rm Frac}\,}
\def\Card{\text{\rm Card}\,}
\def\Hom{\text{\rm Hom}}
\def\Gal{\text{\rm Gal}}
\def\sep{\text{\rm sep}}
\def\ur{\text{\rm ur}}

\def\mod{\,\text{\rm mod}\,}
\def\Mod{\text{\rm Mod}}
\def\Rep{\text{\rm Rep}}
\def\st{\text{\rm st}}
\def\res{\text{\rm res}}
\def\np{\text{\rm np}}
\def\tnp{\text{\rm -np}}

\def\ent{\text{\rm int}}

\def\alg{\text{\rm alg}}
\def\triv{\text{\rm triv}}
\def\perf{\text{\rm perf}}
\def\cris{\text{\rm cris}}
\def\DP{\text{\rm dp}}

\def\Libre{\text{\rm Libre}}

\title{Représentations galoisiennes $p$-adiques et $(\varphi,\tau)$-modules}
\author{Xavier Caruso}
\date{Septembre 2012}

\begin{document}

\maketitle

\begin{abstract}
Étant donné un nombre premier impair $p$ et un corps $p$-adique $K$, on 
développe dans cet article, un analogue de la théorie des 
$(\varphi,\Gamma)$-modules de Fontaine en remplaçant la $p$-extension 
cyclotomique par l'extension $K_\infty$ de $K$ obtenue en ajoutant un 
système compatible de racines $p^n$-ièmes d'une uniformisante $\pi$ 
fixée. Ceci nous conduit à une nouvelle classification des 
représentations $p$-adiques de $G_K = \Gal (\bar K/K)$ \emph{via} des 
$(\varphi,\tau)$-modules. Nous établissons ensuite un lien entre la 
théorie des $(\varphi,\tau)$-modules à celle des 
$(\varphi,N_\nabla)$-modules de Kisin. Comme corollaire, nous répondons 
à une question de Tong Liu en démontrant que, lorsque $K$ est une 
extension finie de $\Q_p$, toute représentation de $E(u)$-hauteur finie 
de $G_K$ est potentiellement semi-stable.
\end{abstract}

\renewcommand{\abstractname}{Abstract}

\begin{abstract}
Let $p$ be an odd prime number and $K$ be a $p$-adic field. In this
paper, we develop an analogue of Fontaine's theory of
$(\varphi,\Gamma)$-modules replacing the $p$-cyclotomic extension by the
extension $K_\infty$ obtained by adding to $K$ a compatible system of
$p^n$-th roots of a fixed uniformizer $\pi$ of $K$. As a result, we
obtain a new classification of $p$-adic representations of $G_K =
\Gal(\bar K/K)$ by some $(\varphi, \tau)$-modules. We then make a
link between the theory of $(\varphi,\tau)$-modules discussed above
and the so-called theory of $(\varphi,N_\nabla)$-modules developped by 
Kisin. As a corollary, we answer a question of Tong Liu: we prove that, 
if $K$ is a finite extension of $\Q_p$, every representation of $G_K$ of 
$E(u)$-finite height is potentially semi-stable.
\end{abstract}

\setcounter{tocdepth}{2}
\tableofcontents

\bigskip

\noindent   
\hrulefill

\bigskip

Soit $p$ un nombre premier impair. Soit $K$ un corps de caractéristique
nulle, complet pour une valuation discrète, dont le corps résiduel $k$
est parfait de caractéristique $p$. On fixe $\bar K$ une clôture
algébrique de $K$ et on s'intéresse aux représentations $p$-adiques du
groupe de Galois $G_K = \Gal(\bar K/K)$. Plusieurs théories ont déjà été
développées pour étudier ces représentations, et notamment celle des
$(\varphi, \Gamma)$-modules dûe à Fontaine \cite{fontaine}. Dans cette
théorie, la $p$-extension cyclotomique de $K$, notée
$K(\zeta_{p^\infty})$ et obtenue en ajoutant à $K$ les racines
$p^n$-ièmes de l'unité, joue un rôle essentiel. En effet, une première
étape cruciale consiste à démontrer que les $\Q_p$-représentations
(resp. les $\Z_p$-représentations libres ou de torsion) de type fini du
groupe de Galois $H_K = \Gal(\bar K / K(\zeta_{p^\infty}))$ sont
classifiées par certains $\varphi$-modules\footnote{Les définitions
précises seront données par la suite. On pourra se contenter pour
l'instant de savoir que les $\varphi$-modules sont des modules munis
d'un opérateur semi-linéaire, généralement noté $\varphi$.} sur le corps
$\E$ (resp. sur l'anneau $\E^\ent$) défini comme suit
$$\E = \Big\{ \, \sum_{i \in \Z} a_i u^i \quad \Big| \quad a_i \in
W[1/p], \, (a_i) \text{ bornée},\, \lim_{i \to -\infty} a_i = 0 \, 
\Big\}$$
$$\text{(resp. }\E^\ent = \Big\{ \, \sum_{i \in \Z} a_i u^i \quad 
\Big| \quad 
a_i \in W, \,\, \lim_{i \to -\infty} a_i = 0 \, \Big\} \text{)}$$
où $W = W(k)$ désigne l'anneau des vecteurs de Witt à coefficients dans $k$.
À partir de là, on retrouve l'action complète de $G_K$ en ajoutant à ce
$\varphi$-module une action du quotient $\Gamma = G_K/H_K =
\Gal(K(\zeta_{p^\infty}) / K)$, obtenant au final un objet appelé
$(\varphi,\Gamma)$-module.

Certains travaux récents de Breuil et Kisin (voir \cite{breuil} et
\cite{kisin}) ont montré que, plus encore que la $p$-extension
cyclotomique considérée précédemment, l'extension $K_\infty$ de $K$,
obtenue en ajoutant un système compatible de racines $p^n$-ièmes d'une
uniformisante $\pi$ fixée, joue un rôle particulier pour l'étude des
représentations semi-stables de $G_K$.
Le but de cet article est de développer un analogue de la théorie des
$(\varphi,\Gamma)$-modules à partir de l'extension $K_\infty$. La
première étape, qui consiste à classifier les représentations de
$G_\infty = \Gal(\bar K/ K_\infty)$, fonctionne comme dans le cas
classique : à une telle représentation est associée un $\varphi$-module
sur $\E$ ou $\E^\ent$ (selon que les coefficients soient pris dans
$\Q_p$ ou $\Z_p$), et celui-ci suffit à la décrire complètement. Par
contre, une fois arrivé à ce niveau, il n'est plus possible de copier
à la lettre la méthode de Fontaine, tout simplement parce que l'extension
$K_\infty/K$ n'est pas galoisienne. Cela n'a donc aucun sens d'ajouter
au $\varphi$-module précédent l'action résiduelle de $\Gal(K_\infty/K)$
puisque ce dernier groupe n'est pas défini ! 

Suivant certains travaux de Liu (voir \cite{liu1} et \cite{liu2}), on
peut toutefois procéder comme suit. On considère un élément $\tau \in
G_K$ agissant trivialement sur $K(\zeta_{p^\infty})$ --- ou même
seulement sur $K(\zeta_p)$ --- tel que $\tau$ et $G_\infty$ engendrent
ensemble $G_K$. Dans ces conditions, connaître le $\varphi$-module $M$
et l'action de $\tau$ suffit à reconstruire l'action complète de $G_K$
sur $T$. Cette action supplémetaire de $\tau$ a un pendant au niveau des
$\varphi$-modules. Elle ne correspond certes pas à un simple
endomorphisme de $M$ (comme cela aurait été le cas si l'extension
$K_\infty/K$ avait été galoisienne), mais à un endomorphisme
semi-linéaire de ${\E^\ent_\tau} \otimes_{\E^\ent} M$ où
${\E^\ent_\tau}$ est une certaine $\E^\ent$-algèbre munie d'une action
de $G_K$. Ceci nous conduit à définir un $(\varphi, \tau)$-module sur
$(\E^\ent, {\E^\ent_\tau})$ (resp. sur $(\E,\E_\tau)$ où $\E_\tau =
\E^\ent_\tau[1/p]$) comme la donnée d'un $\varphi$-module de type fini
$M$ sur $\E^\ent$ (resp. $\E$) muni d'une application supplémentaire
$\tau : {\E^\ent_\tau} \otimes_{\E^\ent} M \to {\E^\ent_\tau}
\otimes_{\E^\ent} M$ vérifiant un certain nombre de conditions (voir
définition \ref{def:phitau} pour plus de précisions). Nous démontrons
alors le théorème suivant (se reporter au \S\S \ref{subsec:equivmodp} et
\ref{subsec:defphitau} pour la définition des foncteurs).

\begin{introtheo}
\label{introtheo:equivphitau}
Il existe des équivalences de catégories :
\begin{eqnarray*}
\left\{ \begin{array}{c}
\Q_p\text{-représentations de }\\
\text{dimension finie de }G_K 
\end{array} \right\}
& \stackrel{\sim}{\longrightarrow} & 
\left\{ \begin{array}{c}
(\varphi, \tau)\text{-modules sur }(\E, \E_\tau)
\end{array} \right\}
\\
\left\{ \begin{array}{c}
\Z_p\text{-représentations}\\
\text{libres de type fini de }G_K 
\end{array} \right\}
& \stackrel{\sim}{\longrightarrow} & 
\left\{ \begin{array}{c}
(\varphi, \tau)\text{-modules} \\
\text{libres sur }(\E^\ent, {\E^\ent_\tau})
\end{array} \right\}
\end{eqnarray*}
et, pour tout entier $n$ :
\begin{eqnarray*}
\left\{ \begin{array}{c}
\Z_p\text{-représentations de type}\\
\text{fini de }G_K \text{ annulées par } p^n 
\end{array} \right\}
& \stackrel{\sim}{\longrightarrow} & 
\left\{ \begin{array}{c}
(\varphi, \tau)\text{-modules}
\text{ sur} \\
(\E^\ent, {\E^\ent_\tau})
\text{ annulés par }p^n
\end{array} \right\}.
\end{eqnarray*}
\end{introtheo}


Dans le \S \ref{sec:reseaux}, nous introduisons et étudions la notion de 
$(\varphi,\tau)$-réseau. Si $M$ est un $(\varphi,\tau)$-module, un 
$(\varphi,\tau)$-réseau de $M$ est un sous-$W[[u]]$-module de type fini 
de $M$, qui engendre $M$ comme $\E^\ent$-module, qui est stable par 
$\varphi$ et dont l'extension des scalaires à $\Sk_\tau$ (un certain 
sous-anneau de $\E^\ent_\tau$ qui sera défini dans le corps du texte) 
est stable par $\tau$. Une notion importante liée aux réseaux est celle 
de hauteur : étant donné un élément $U \in \Sk_\tau$, on dit qu'un 
réseau $\frakM$ est de hauteur divisant $U$ si le conoyau de
$$\id \otimes \varphi : \Sk_\tau \otimes_{\varphi,\Sk} \frakM \to
\Sk_\tau \otimes_\Sk \frakM$$
est annulé par $U$. Il s'avère que l'existence d'un 
$(\varphi,\tau)$-module de hauteur divisant $U$ à l'intérieur d'un
$(\varphi,\tau)$-module implique des bornes explicites sur la 
ramification de la représentation galoisienne correspondante. Plus
exactement, nous démontrerons le théorème \ref{introtheo:bornes}
ci-après.

\begin{introtheo}
\label{introtheo:bornes}
On note $G_\infty^{(\mu)}$ et $G_K^{(\mu)}$ les filtrations de
ramification en numérotation supérieure\footnote{Pour le groupe $G_K$,
il s'agit de la ramification usuelle telle que définie, par exemple,
dans \cite{fontaine}. Sur le groupe $G_\infty$, la ramification est
définie de même \emph{via} l'isomorphisme de corps des normes $G_\infty
\simeq \Gal(F_0^\sep/F_0)$ où $F_0 = \E^\ent/p\E^\ent \simeq k((u))$.
Pour plus de précisions à ce sujet, on renvoie au \S
\ref{subsec:ramif}.}
des groupes de Galois respectifs $G_\infty$ et
$G_K$. On fixe un nombre entier $n \geq 1$ et un élément $U \in \Sk$ qui
n'est pas multiple de $p$. On pose $h = v_R(U \mod p)$.

\smallskip

\noindent {\bf 1.}
Soit $T$ une $\Z_p$-représentation de $G_\infty$ qui est de type fini
comme $\Z_p$-module et annulée par $p^n$. On suppose que le
$\varphi$-module étale sur $\E^\ent$ associé à $T$ admet un
$\varphi$-réseau de hauteur divisant $U$. Alors pour tout $\mu > 
\max(1, \frac{h p^n}{p-1})$, le sous-groupe $G_\infty^{(\mu)}$ agit 
trivialement sur $T$.

\smallskip

\noindent {\bf 2.}
Il existe une constante $c(K)$ ne dépendant que du corps $K$ telle que
l'assertion suivante soit vraie : pour toute $\Z_p$-représentation $T$
de $G_K$ de type fini comme $\Z_p$-module, annulée par $p^n$, dont la
restriction à $G_\infty$ correspond à un $\varphi$-module étale sur
$\E^\ent$ de hauteur divisant $U$, et pour tout 
$\mu > c(K) + e \cdot \max (\frac 1{p-1}, n + \log_p (\frac h e))$, 
le groupe $G_K^{(\mu)}$ agit trivialement sur $T$.
\end{introtheo}

\noindent
{\it Remarque.}
La constante $c(K)$ dépend de façon assez explicite de la ramification
absolue de $K$ ; par exemple, lorsque $e$ est premier avec $p$
(\emph{i.e.} lorsque $K$ est absolument modérément ramifié), elle peut
être choisie égale à $1 + e + \frac e{p-1}$.

\medskip

Nous nous intéressons enfin plus particulièrement aux 
$(\varphi,\tau)$-modules qui sont de hauteur divisant $E(u)^r$ pour un 
certain entier $r$ (on dit aussi : \emph{de $E(u)$-hauteur $\leq r$}). 
L'intérêt de cette notion a été récemment mis en lumière dans un premier 
temps par Breuil dans \cite{breuil} puis par Kisin dans \cite{kisin} qui 
a démontré dans \emph{loc. cit.} qu'une représentation semi-stable était 
nécessairement de $E(u)$-hauteur finie. Réciproquement, dans le \S 
\ref{sec:Eu}, nous démontrons le théorème suivant :

\begin{introtheo}
\label{introtheo:Eu}
On suppose que $K$ est une extension finie de $\Q_p$.
On note $s$ le plus grand entier tel que $K$ contienne une racine 
primitive $p^s$-ième de l'unité. Alors, toute représentation de 
$E(u)$-hauteur finie de $G_K$ devient semi-stable en restriction au 
sous-groupe (distingué) $G_s = \Gal(\bar K/K(\!\sqrt[p^s] \pi))$.
\end{introtheo}

Il est à noter que les bornes de ramification données par le théorème 
\ref{introtheo:bornes} jouent un rôle absolument essentiel dans la 
démonstration du théorème \ref{introtheo:Eu} (il en est en fait le point 
de départ). En effet, elles fournissent des bornes sur l'action de 
$\tau$ sur le $(\varphi,\tau)$-module associé à $V$ qui sont exactement 
celles dont on a besoin pour construire un opérateur $\log \tau$ sur le 
$(\varphi,\tau)$-module associée à $V$ (en se basant la formule 
$\sum_{i=1}^\infty \frac{(\id-\tau)^i} i$). À partir de cet opérateur, 
on fabrique ensuite un $(\varphi,N)$-module filtré à la Fontaine duquel 
émergera finalement une représentation semi-stable. Il ne restera alors 
plus qu'à démontrer que cette dernière représentation coïncide avec 
celle dont on est partie, au moins en restriction au sous-groupe $G_s$.

\medskip

Nous démontrons enfin au \S \ref{subsec:bornesst} un raffinement du
théorème \ref{introtheo:bornes} pour les représentations semi-stable
qui n'est autre que la conjecture 1.2.(1) de \cite{carliu}.

\medskip

Ce travail puise son origine et son inspiration dans de nombreuses
discussions avec Tong Liu ; à travers ces quelques lignes, l'auteur
souhaite lui témoigner tous ses remerciements.
L'auteur est également reconnaissant à l'Agence Nationale de la 
Recherche (ANR) pour son soutien financier par l'intermédiaire du projet 
CETHop (Calculs Effectifs en Théorie de Hodge $p$-adique) référencé 
ANR-09-JCJC-0048-01.

\numberwithin{equation}{section}

\section{La théorie générale des $(\varphi,\tau)$-modules}
\label{sec:phitau}

\begin{figure}[t]
\begin{center}
\includegraphics{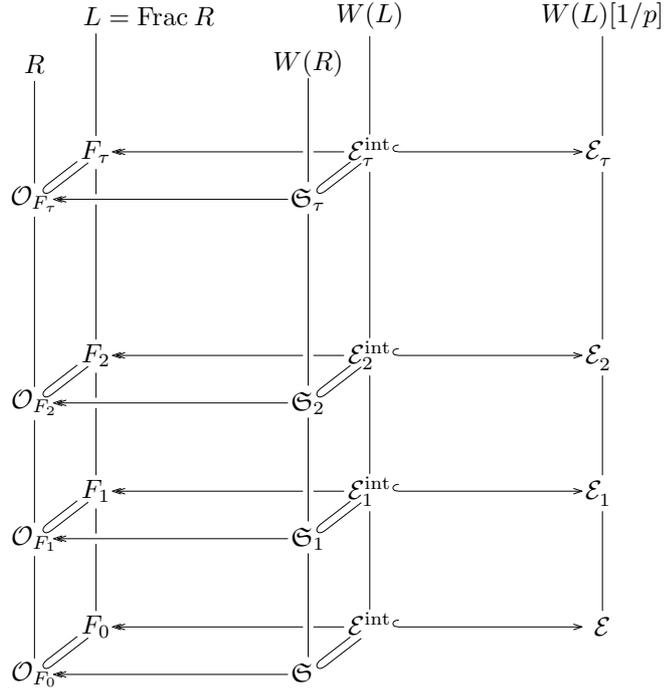}

\caption{Principaux anneaux intervenant dans la théorie des
$(\varphi,\tau)$-modules}
\label{fig:anneaux}

\bigskip

\begin{minipage}{13cm}
La partie gauche (resp. centrale, resp. droite) du diagramme concerne la
théorie sur $\F_p$ (resp. $\Z_p$, resp. $\Q_p$). Les anneaux qui
apparaissent au second plan concernent la théorie générale des
$(\varphi,\tau)$-modules ; ils sont présentés et étudiés dans le \S
\ref{sec:phitau}. Au premier plan, au contraire, nous trouvons les anneaux
nécessaires à l'étude des $(\varphi,\tau)$-réseaux ; ceux-ci sont 
introduits et utilisés dans le \S \ref{sec:reseaux}.
\end{minipage}
\end{center}
\end{figure}

L'objectif principal de cette première partie est de définir les
foncteurs qui réalisent les équivalences de catégories énoncées dans le
théorème \ref{introtheo:equivphitau} de l'introduction (et donc de
définir, aussi, en particulier, l'anneau ${\E^\ent_\tau}$), puis de
démontrer ce théorème.
Tout au long de cette section et de la suivante, nous allons petit à
petit définir un certain nombre d'anneaux. La figure \ref{fig:anneaux}
présente un diagramme qui fait apparaître les plus importants d'entre
eux, ainsi que les morphismes essentiels les reliant. Nous espérons que
celui-ci pourra faciliter la lecture de cet article.

Comme dans l'introduction, on fixe un nombre premier impair
$p$ et un corps $K$ muni d'une valuation discrète pour laquelle il est
complet. On suppose que $K$ est de caractéristique nulle et que son
corps résiduel $k$ est parfait de caractéristique $p$. On rappelle
également que $W$ désigne l'anneau des vecteurs de Witt à
coefficients dans $k$ ; son corps des fractions $W[1/p]$ s'identifie
canoniquement au plus gros sous-corps de $K$ absolument non ramifiée.
L'extension $K/W[1/p]$ est totalement ramifiée et on note $e$ son degré.
On considère $(\zeta_{p^s})$ un système compatibles de racines
primitives $p^s$-ièmes de l'unité et on note $K(\zeta_{p^\infty}) =
\bigcup_s K(\zeta_{p^s})$. L'extension $K(\zeta_{p^\infty})$ est
galoisienne et son groupe de Galois $\Gamma$ s'identifie \emph{via} le
caractère cyclotomique $\chi$ à un sous-groupe ouvert de $\Z_p^\times$.

\subsection{L'extension $K_\infty$ : définition et propriétés
galoisiennes}
\label{subsec:Kinfty}

Nous commençons par quelques préliminaires, en grande partie déja
connus, concernant l'extension $K_\infty$ de Breuil-Kisin. Soit $\pi$
une uniformisante de $\O_K$. On choisit une fois pour toutes une système
compatible $(\pi_s)$ de racines $p^s$-ièmes de $\pi$. On pose $K_s =
K(\pi_s)$ pour tout entier $s$ et $K_\infty = \bigcup_s K_s$. Les
extensions $K_s/K$ ne sont pas galoisiennes, sauf éventuellement pour
les premiers entiers $s$. De façon plus précise, la clôture galoisienne
de $K_s/K$ est l'extension $K_s(\zeta_{p^s})$. Ainsi $K_s/K$ est
galoisienne si, et seulement si $\zeta_{p^s} \in K_s$. L'extension
$K_\infty/K$, quant à elle, n'est \emph{jamais} galoisienne et sa
clôture galoisienne est l'extension composée $K_\infty(\zeta_{p^\infty})
= K_\infty {\cdot} K(\zeta_{p^\infty})$. On pose $G_\infty = \Gal(\bar
K/K_\infty) \subset G_K$ et $H_\infty = \Gal(\bar
K/K_\infty(\zeta_{p^\infty})) \subset G_\infty$.

\subsubsection{Le cocycle $c$}

Bien que $G_\infty$ ne soit pas distingué dans $G_K$, on peut définir
une application continue (qui bien sûr, n'est pas un morphisme de
groupes) $c : G_K \to \Z_p$ dont le \og noyau \fg\ (\emph{i.e.} l'image
réciproque de $0 \in \Z_p$) s'identifie à $G_\infty$. Pour cela, on
remarque qu'étant donné un élément $g \in G_K$ et un entier $s$, il
existe un unique élément $c_s(g) \in \Z/p^s\Z$ tel que $g(\pi_s) =
\zeta_{p^s}^{c_s(g)} \pi_s$. La famille des $(c_s(g))_{s \geq 1}$
vérifie la congruence $c_{s+1}(g) \equiv c_s(g) \pmod {p^s}$, et
définit donc un élément de $\Z_p,$ que l'on appelle $c(g)$. On vérifie
sans peine que $c$ est continue et que le fait que $c(g)$ s'annule
signifie exactement que $g$ agit trivialement sur chacun des $\pi_s$,
c'est-à-dire sur $K_\infty$. Ainsi, comme nous le disions précédemment,
on a $c^{-1}(0) = G_\infty$.
L'application $c$ n'est certes pas un morphisme de groupes mais vérifie
malgré tout une relation de $1$-cocycle à savoir :
\begin{equation}
\label{eq:cgh}
\forall g,h \in G_\infty, \quad c(gh) = c(g) + \chi(g) c(h).
\end{equation}
On en déduit que, si $\Z_p \rtimes \Z_p^\times$ désigne le produit
semi-direct de $\Z_p$ par $\Z_p^\times$ où $\Z_p^\times$ agit sur $\Z_p$
par multiplication, l'application $\chi_\infty : G_K \to \Z_p \rtimes
\Z_p^\times$, $g \mapsto (c(g), \chi(g))$ est un morphisme de groupes dont
le noyau est exactement $H_\infty$. En particulier, le $1$-cocycle $c :
G_K \to \Z_p$ se factorise par $G_K/H_\infty$, et ce dernier groupe se
plonge dans $\Z_p \rtimes \Z_p^\times$ \emph{via} $\chi_\infty$. Voici un
diagramme qui résume les liens entre les différents groupes que l'on
vient d'introduire :
$$\xymatrix{
0 \ar[r] & \Gal(K_\infty(\zeta_{p^\infty}) / K_\infty) \ar[r] \ar[d]
& G_K/H_\infty \ar[r]^-{\chi} \ar@{^(->}[d]^-{\chi_\infty} & \Z_p^\times 
\ar@{=}[d] \\
0 \ar[r] & \Z_p \ar[r] &
\Z_p \rtimes \Z_p^\times \ar[r] & \Z_p^\times \ar[r] & 0
}$$
Soulignons que les groupes qui apparaissent sur le diagramme sont
tous profinis, tandis que les morphismes sont tous continus pour la
topologie profinie. Par ailleurs, il résulte du diagramme que le
morphisme $\Gal(K_\infty(\zeta_{p^\infty}) / K_\infty) \to \Z_p$ est
injectif. Le lemme 5.1.2 de \cite{liu-compo} (qui stipule que les extensions
$K_\infty$ et $K(\zeta_{p^\infty})$ sont linéairement disjointes sur
$K$) montre même que c'est un isomorphisme. On déduit également du même
lemme que $\chi(G_\infty) = \chi(G_K)$. On remarque encore que si $g \in
G_K$ est tel que $\chi(g) \equiv 1 \pmod p$, alors $\chi_\infty(g)$
appartient à l'unique pro-$p$-Sylow de $\Z_p \rtimes \Z_p^\times$ ; il
en résulte que $\chi(g)^{p^s}$ converge vers l'élément neutre de ce groupe
lorsque $s$ tend vers l'infini, et donc que la suite des $g^{p^s}
\mod H_\infty$ converge elle aussi (vers l'élément neutre de
$G_K/H_\infty$). En particulier, si $a$ est un entier $p$-adique,
l'élément $g^a$ est bien défini dans $G_K/H_\infty$, et il en est donc
de même de $c(g^a)$. Une récurrence à partir de \eqref{eq:cgh} suivie
d'un passage à la limite conduit enfin à la formule
\begin{equation}
\begin{array}{ccrcll}
\label{eq:cga}
c(g^a) = c(g) \cdot [a]_{\chi(g)}
& \text{où} &
[a]_{\chi(g)} & = & a & \text{si } \chi(g) = 1.
\medskip \\
& & & = & \displaystyle \frac{\chi(g)^a-1}{\chi(g)-1} 
& \text{si } \chi(g) \neq 1.
\end{array}
\end{equation}
Dans la suite, le nombre $[a]_{\chi(g)}$ sera appelé le
\emph{$\chi(g)$-analogue} de $a$.

\subsubsection{L'élément $\tau$}

Soit $\tau : K_\infty \to \bar K$ le $K$-plongement défini par
$\tau(\pi_s) = \zeta_{p^s} \pi_s$. Les extensions $K_\infty$ et
$K(\zeta_p)$ étant linéairement disjointes sur $K$, on peut prolonger
$\tau$ à $\bar K$ de façon à ce que $\chi(\tau) \equiv 1 \pmod p$. En
réalité, le lemme 5.1.2 de \cite{liu-compo} nous dit même que l'on peut
imposer $\chi(\tau) = 1$. Toutefois, bien que cela complique légèrement
les calculs, nous préférons continuer à travailler avec un $\tau$ plus
général vérifiant seulement $\chi(\tau) \equiv 1 \pmod p$. Il suit de la
définition de $\tau$ que $c(\tau) = 1$. En outre, en vertu de ce qui a
été dit juste en dessous de cet alinéa, on peut définir $c(\tau^a)$ et
$\chi(\tau)^a$ pour tout $a \in \Z_p$. La formule \eqref{eq:cga} s'écrit
alors simplement $c(\tau^a) = [a]_{\chi(\tau)}$.

\begin{lemme}
\label{lem:bijanalogue}
L'application $\Z_p \to \Z_p$, $a \mapsto [a]_{\chi(\tau)}$ est une
bijection. De plus, pour tout $a \in \Z_p$, on a $a \equiv
[a]_{\chi(\tau)} \pmod p$ et les valuations $p$-adiques de $a-1$ et de
$[a]_{\chi(\tau)}-1$ sont égales.
\end{lemme}

\begin{proof}
Si $\chi(\tau) = 1$, l'application est l'identité et il n'y a rien
à démontrer. Dans le cas contraire, pour démontrer que $a \mapsto
[a]_{\chi(\tau)}$ est bijective, il suffit de remarquer que
l'inverse est donnée par $b \mapsto \frac{\log (1 + (\chi(\tau) - 1)
b)}{\log \chi(\tau)}$, application qui est bien définie car $\chi(\tau)
\equiv 1 \pmod p$ et $p > 2$.
Pour la deuxième assertion, on écrit $\chi(\tau) = 1 + p x$. On a alors
$$[a]_{\chi(\tau)} = a + \frac{a(a-1)} 2 \cdot p x +
\frac{a(a-1)(a-2)} 6 \cdot p^2 x^2 + \cdots + \frac{a(a-1) \cdots
(a-n+1)} {n!} \cdot p^n x^n + \cdots$$
et il est déjà clair que $a \equiv [a]_{\chi(\tau)} \pmod p$. On
remarque ensuite que de $a \equiv 1 \pmod {p^s}$, il suit
$[a]_{\chi(\tau)} \equiv a \pmod{p^{s+1}}$ car $\frac{p^n} {n!}$ est
toujours multiple de $p$ dans $\Z_p$ (on rappelle que l'on suppose $p >
2$) et $a-1$ est lui, par définition, multiple de $p^s$. La conclusion
en résulte.
\end{proof}


D'après le lemme, pour tout élément $g \in G_K$, il existe un unique $a
\in \Z_p$ tel que $[a]_{\chi(\tau)} = \chi(g)$. Notons $\chi_\tau (g)$
cet élément ; on définit de cette manière une fonction $\chi_\tau : G_K
\to \Z_p$ qui, d'après la congruence $a \equiv [a]_{\chi(\tau)} \pmod
p$, prend ses valeurs dans $\Z_p^\times$. On prendra garde au fait que
$\chi_\tau$ n'est pas un morphisme de groupes ; cependant, son \og noyau
\fg\ (\emph{i.e.} l'image réciproque de $1 \in \Z_p^\times$) est le
sous-groupe $H_K$. L'égalité entre valuations établie dans le lemme
\ref{lem:bijanalogue} assure en outre que, pour tout entier $m$, l'image
réciproque du sous-groupe $1 + p^m\Z_p$ de $\Z_p^\times$ est le
sous-groupe $\Gal(\bar K / K(\zeta_{p^m}))$ de $G_K$.

\begin{lemme}
\label{lem:strucGK}
\begin{enumerate}
\item Tout élément $g \in G_K/H_\infty$ s'écrit de façon unique
sous la forme $\tau^a g'$ avec $a \in \Z_p$ et $g' \in 
G_\infty/H_\infty$

\item Pour tout $g \in G_\infty/H_\infty$, le produit
$\tau^{-\chi_\tau(g)} g \tau$ calculé dans $G_K/H_\infty$ appartient à
$G_\infty/H_\infty$.
\end{enumerate}
\end{lemme}

\begin{proof}
Pour le premier alinéa, il s'agit de montrer que l'équation $c(\tau^{-a}
g) = 0$ a une unique solution dans $\Z_p$. Or on a
$c(\tau^{-a} g) = [-a]_{\chi(\tau)} + \chi(\tau)^{-a} c(g)
= \chi(\tau)^{-a} \big( c(g) - [a]_{\chi(\tau)} \big)$, à partir de quoi 
le lemme \ref{lem:bijanalogue} permet de conclure. Pour le deuxième 
alinéa, il suffit de vérifier que $c(\tau^{-\chi_\tau(g)} g \tau) = 0$, ce 
qui se fait comme précédemment.
\end{proof}

\begin{prop}
\label{prop:prolongmorp}
Soit $\Gamma$ un groupe topologique et $\rho : G_\infty/H_\infty \to
\Gamma$ un morphisme de groupes continu. Se donner un prolongement
continu $\tilde \rho : G_K/H_\infty \to \Gamma$ de $\rho$ revient à se
donner un élément $\tau_\Gamma \in \Gamma$ (l'image de $\tau$) tel que
pour tout $g \in G_\infty/H_\infty$ vérifiant $\chi_\tau(g) \in \N$, on
ait :
\begin{equation}
\label{eq:commutation}
\rho(g) \cdot \tau_\Gamma = \tau_\Gamma^{\chi_\tau(g)} \cdot 
\rho(\tau^{-\chi_\tau(g)} g \tau).
\end{equation}
\end{prop}

\begin{proof}
D'après le premier point du lemme \ref{lem:strucGK}, il est clair que
$\tau_\Gamma = \tilde \rho(\tau)$ détermine entièrement $\tilde \rho$. 
L'unicité en résulte. Pour l'existence, on montre dans un premier temps que $\lim_{s \to \infty}
\tau_\Gamma^{p^s} = 1$. Étant donné que $\chi(G_\infty) = \chi (G_K)$ est
ouvert dans $\Z_p^\times$, il existe un entier $s_0$ tel que 
$\chi(G_\infty) \supset 1 + p^{s_0} \Z_p$. Comme d'après le lemme
\ref{lem:bijanalogue}, l'application $a \mapsto [a]_{\chi(\tau)}$ 
réalise une bijection de $1 + p^s \Z_p$ dans lui-même pour tout $s$, 
il existe une suite $(g_s)_{s \geq s_0}$ d'éléments de $G_\infty/
H_\infty$ convergeant vers l'identité et telle que $\chi_\tau(g_s) = 
p^s + 1$. En appliquant la relation \eqref{eq:commutation} avec $g_s$,
il vient 
$\tau_\Gamma^{p^s} = \rho(g_s) \cdot \tau_\Gamma \cdot \rho(\tau^{-\chi_\tau(g_s)} 
g_s \tau) \cdot \tau_\Gamma^{-1}$.
Un passage à la limite pour $s$ tendant vers l'infini donne alors la
conclusion annoncée. Ceci nous permet de définir $\tau_\Gamma^a$ pour tout
élément $a \in \Z_p$ et en passant à la limite dans 
\eqref{eq:commutation}, on montre que cette dernière relation est
valable pour tout $g \in G_\infty/H_\infty$. Il suffit maintenant de 
montrer que l'application $\tilde \rho$ définie par
$$\forall a \in \Z_p, \, \forall g \in G_\infty/H_\infty, \quad
\tilde \rho(\tau^a g) = \tau_\Gamma^a \rho(g)$$
est un morphisme de groupes. Il s'agit donc de vérifier que si $g \tau^a
= \tau^b h$ avec $a, b \in \Z_p$ et $g, h \in G_\infty/H_\infty$, alors 
$\rho(g) \tau_\Gamma^a = \tau_\Gamma^b \rho(h)$. Tout d'abord, en appliquant $c$ à 
l'égalité $g \tau^a = \tau^b h$, on obtient $\chi(g) \cdot [a]_{\chi(\tau)} = 
[b]_{\chi(\tau)}$. Si $a = 1$, on voit à présent que l'égalité que l'on
a à démontrer n'est autre que l'hypothèse. Pour $a \in \N$, l'égalité se
démontre par récurrence, tandis qu'enfin, pour $a \in \Z_p$, on utilise 
un argument de passage à la limite.
\end{proof}

Bien entendu, le cas qui nous intéressera particulièrement dans cet
article est celui où le groupe $\Gamma$ est le groupe des automorphismes
linéaires ou semi-linéaires d'un certain espace. Dans cette situation,
la proposition dit exactement ce qu'il faut ajouter à une représentation
de $G_\infty/H_\infty$ pour définir une représentation de
$G_K/H_\infty$. 

\subsubsection{Comment se dispenser de quotienter par $H_\infty$}
\label{subsec:quotHinf}

Dans le lemme \ref{lem:strucGK} et la proposition
\ref{prop:prolongmorp}, nous avons systématiquement quotienter par
$H_\infty$. Toutefois, on peut s'affrachir de cela en prenant soin de
choisir au préalable un élément $\tau \in G_K$ vérifiant $\lim_{s \to
\infty} \tau^{p^s} = \id$. Un tel choix est toujours possible car, étant
donné que l'extension $K_\infty/K$ n'admet pas de sous-extension
modérément ramifiée, on peut choisir $\tau$ dans le sous-groupe
d'inertie sauvage ; comme celui-ci est un pro-$p$-groupe, la convergence
requise en résulte. Une fois ce choix fait, les démonstrations du lemme
\ref{lem:strucGK} et de la proposition \ref{prop:prolongmorp} s'étendent
mot pour mot en remplaçant partout $G_K/H_\infty$ par $G_K$ et
$G_\infty/H_\infty$ par $G_\infty$. En particulier, on voit que se
donner une action de $G_K$ revient, dans ce cas, à se donner une action
de $G_\infty$ et un automorphisme $\tau_\Gamma$ (correspondant à
l'action de $\tau$) vérifiant la relation de commutation
\eqref{eq:commutation}. Toutefois, dans la suite, nous nous contenterons
d'appliquer la proposition \ref{prop:prolongmorp} au groupe quotient
$G_K/H_\infty$, et il ne nous sera donc pas nécessaire de particulariser
ainsi le choix de $\tau$.

\medskip

Un fait notable, par contre, est qu'un corollaire de la généralisation
que l'on vient d'évoquer permet de donner une description de $G_K$ en
termes de $G_\infty$ à l'aide d'une construction de type \og produit
semi-direct \fg. On choisit pour cela un élément $\tau$ dans le groupe
d'inertie sauvage vérifiant à la fois $c(\tau) = 1$ et $\chi(\tau) = 0$
(un tel choix est possible). On a alors les propriétés suivantes :
\begin{itemize}
\item tout élément de $G_K$ s'écrit de façon unique sous la forme
$\tau^a g$ avec $a \in \Z_p$ et $g \in G_\infty$ ;
\item si $g \in G_\infty$ et $a \in \Z_p$, alors $h = \tau^{-a\chi(g)} 
g \tau^a \in G_\infty$ et on a bien sûr la relation $g \tau^a = \tau^{a 
\chi(g)} h$.
\end{itemize}
Ainsi si l'on se donne $G_\infty$, la restriction du caractère
cyclotomique $\chi : G_\infty \to \Z_p^\times$ et l'application $\psi :
G_\infty \to G_\infty$, $g \mapsto \tau^{-\chi(g)} g \tau$, on peut
reconstruire le groupe $G_K$ tout entier en considérant l'ensemble $\Z_p
\times G_\infty$ muni de la loi de groupe suivante :
$$(a,g) \cdot (b,h) = (a + b\chi(g), \psi^b(g) h).$$
On prendra garde néanmoins au fait que l'application $\psi$ n'est pas 
un morphisme de groupes ; elle vérifie cependant une relation, à savoir
$\psi(gh) = \psi^{\chi(h)}(g) \psi(h)$.

\subsection{Les $(\varphi,\tau)$-modules en caractéristique $p$}

Dans ce paragraphe, nous mettons au point la théorie des
$(\varphi,\tau)$-modules en caractéristique $p$, c'est-à-dire que nous
établissons la troisième équivalence de catégories du théorème
\ref{introtheo:equivphitau} de l'introduction lorsque $n=1$ (les anneaux
$\E^\ent$ et ${\E^\ent_\tau}$ étant alors remplacés par des variantes
plus simples). La méthode est simple et naturelle : elle consiste à
mettre ensemble ce qui vient d'être fait et le théorème de Fontaine et
Fontaine-Wintenberger de classification des représentations de
$G_\infty$ \emph{via} les $\varphi$-modules étales. Nous commençons par
quelques rappels à ce sujet.

\subsubsection{Rappels sur la classification des représentations
de $G_\infty$}

On considère l'anneau $R = \varprojlim \O_{\bar K}/p$ où les morphismes
de transition sont donnés par l'élévation à la puissance $p$ : un
élément $x \in R$ est donc une suite $(x_s)_{s \geq 0}$ telle que
$x_{s+1}^p = x_s$ pour tout $s$. On note $\Frac R$ le corps des
fractions de $R$ ; c'est un corps algébriquement clos. À côté de cela,
il est clair que $R$ est muni d'une action canonique de $G_K$ qui
s'étend à $\Frac R$. Par ailleurs, si l'on note $v_K$ la valuation sur
$\bar K$ normalisée par $v_K(K^\star) = \Z$, on démontre que la formule
$v_R(x) = \lim_{s \to \infty} p^s v_K(x_s)$ définit une valuation (non
discrète) sur $R$ pour laquelle il est complet. La valuation $v_R$
s'étend naturellement à $\Frac R$ et on note encore $v_R$ ce
prolongement.  Le corps résiduel $k$ s'injecte canoniquement dans $R$ :
à un élément $\lambda \in k$, on associe la suite des $\lambda^{1/p^s}$
vus comme éléments de $\O_{\bar K}/p$. La famille $(\zeta_{p^s})$ 
(resp. $(\pi_s)$) de racines $p^s$-ièmes de l'unité (resp. 
de $\pi$) définit, elle aussi, un élément de $R$ que l'on note $\ueps$ 
(resp. $\upi$). On a $v_R(\ueps) = 0$ et
$v_R(\upi) = 1$. On pose $\ueta = 1 - \ueps$ de sorte que $v_R(\ueta) =
\frac {ep}{p-1}$.  Par construction le groupe $H_K$ agit trivialement
sur $\ueps$ et $\ueta$ tandis que $G_\infty$, lui, agit trivialement sur
$\upi$. Dans la suite, on plongera systématiquement l'anneau $k[[u]]$
dans $R$ en envoyant $u$ sur $\upi$. Ce morphisme s'étend aux corps des
fractions et définit ainsi un plongement du corps $F_0 = k((u))$ dans
$\Frac R$. On appelle $F_0^\sep$ la clotûre séparable de $F_0$ dans
$\Frac R$. Étant donné que $G_\infty$ agit trivialement sur $F_0$ vu
dans $\Frac R$, tout élément de $G_\infty$ stabilise $F_0^\sep$ et on
obtient comme ceci un morphisme $G_\infty \to \Gal(F_0^\sep /F_0)$. La
théorie du corps des normes de Fontaine et Wintenberger (voir
\cite{wintenberger}) affirme que c'est en fait un isomorphisme. Comme
corollaire du théorème de Hilbert 90, Fontaine démontre alors la
proposition suivante.

\begin{prop}
\label{prop:fontaine}
Soit $T$ une $\F_p$-représentation du groupe $G_\infty$. On suppose
que $T$ est de dimension finie $d$ sur $\F_p$. Alors l'espace 
$\Hom_{\F_p[G_\infty]}(T, F_0^\sep)$ est de dimension $d$ sur $F_0$ et, 
plus précisément, le morphisme naturel
$$F_0^\sep \otimes_{F_0} \Hom_{\F_p[G_\infty]}(T, F_0^\sep) \to
\Hom_{\F_p}(T, F_0^\sep)$$
est un isomorphisme.
\end{prop}

\begin{proof}
C'est un cas particulier\footnote{Dans cette référence, on travaille
déjà avec des coefficients dans $\Z_p$, ce que, de notre côté, nous ne
ferons que dans le \S \ref{subsec:relevement} (voir théorème
\ref{theo:fontaine}).} de la proposition A.1.2.6 et de la remarque 
A.1.2.7 de \cite{fontaine}.
\end{proof}

\begin{rem}
\label{rem:injFontaine}
Le vrai contenu de la proposition réside dans la surjectivité de
l'application. L'injectivité, quant à elle, est un fait beaucoup plus
général, qui reste valable si l'on remplace $G_\infty$ par n'importe
quel groupe topologique $H$, $F_0^\sep$ par n'importe quel corps $L$ de
caractéristique $p$ sur lequel $H$ agit continument, et $F_0$ par $L^H$
(la démonstration étant en tout point identique).
\end{rem}

À partir de la proposition précédente, Fontaine déduit un théorème de 
classification des $\F_p$-représentations de dimension finie de 
$G_\infty$. Pour l'énoncer, on doit d'abord définir la notion de 
\emph{$\varphi$-module étale sur $F_0$} : il s'agit de la donnée d'un 
espace vectoriel $M$ de dimension finie sur $F_0$ muni d'une application 
$\varphi_M : M \to M$ semi-linéaire par rapport au Frobenius sur $F_0$ 
(défini comme l'élévation à la puissance $p$), et dont l'image contient 
une base de $M$. Étant donné un $\varphi$-module étale $M$ sur $F_0$, il 
est souvent commode de considérer le linéarisé de $\varphi_M$ défini 
comme l'application $\id \otimes \varphi_M : F_0 \otimes_{\varphi, F_0} 
M \to M$. On a alors affaire à une application $F_0$-linéaire, et la 
condition selon laquelle l'image de $\varphi_M$ contient une base de $M$ 
se traduit en disant que $\id \otimes \varphi_M$ est un isomorphisme. Si 
$T$ est une représentation de $G_\infty$, alors 
$\Hom_{\F_p[G_\infty]}(T, F_0^\sep)$ est naturellement muni d'un 
endomorphisme $\varphi$ déduit du Frobenius usuel sur $F_0^\sep$, et il 
est facile de vérifier que c'est en fait un $\varphi$-module étale sur 
$F_0$.

\begin{theo}
\label{theo:fontaine}
Les foncteurs suivants sont des équivalences de catégories
quasi-inverses l'une de l'autre :
\begin{eqnarray*}
\left\{ \begin{array}{c}
\F_p\text{-représentations de }\\
\text{dimension finie de }G_\infty 
\end{array} \right\}
& \stackrel{\sim}{\longrightarrow} & 
\left\{ \begin{array}{c}
\varphi\text{-modules étales sur } F_0
\end{array} \right\}
\\
T & \mapsto & \Hom_{\F_p[G_\infty]}(T, F_0^\sep) \\
\Hom_{F_0, \varphi} (M, F_0^\sep) & \mapsfrom & M
\end{eqnarray*}
où $\Hom_{F_0, \varphi}$ signifie que l'on considère les morphismes
$F_0$-linéaires commutant à l'action de $\varphi$ (celui-ci agissant
par l'élévation à la puissance $p$ sur $F_0^\sep$).
\end{theo}

\begin{proof}
Voir proposition A.1.2.6 et remarque A.1.2.7 de \cite{fontaine}.
\end{proof}

\subsubsection{Une équivalence de catégories}
\label{subsec:equivmodp}

On pose à partir de maintenant pour simplifier les écritures $L = \Frac 
R$. On note $F_\tau$ le sous-corps de $L$ formé des éléments fixes par 
$H_\infty$ ; il contient manifestement $(F_0^\sep)^{H_\infty}$ et donc 
en particulier $F_0$.

\begin{deftn}
\label{def:phitaumodp}
Un \emph{$(\varphi,\tau)$-module sur $(F_0,F_\tau)$} est la donnée de
\begin{itemize}
\item un $\varphi$-module étale sur $F_0$, noté $M$ ;
\item un endomorphisme $\tau$-semi-linéaire $\tau_M : F_\tau \otimes_{F_0}
M \to F_\tau \otimes_{F_0} M$ qui commute à $\varphi_{F_\tau} \otimes
\varphi_M$ (où $\varphi_{F_\tau}$ est le Frobenius usuel sur $F_\tau$) et qui
vérifie, pour tout $g \in G_\infty/H_\infty$ tel que $\chi_\tau(g) \in 
\N$, la relation suivante :
\begin{equation}
\label{eq:commuttau}
\forall x \in M, \quad 
(g \otimes \id) \circ \tau_M (x) = \tau_M^{\chi_\tau(g)} (x).
\end{equation}
\end{itemize}
\end{deftn}

\noindent
On est en droit de se demander d'où vient la différence entre la
relation précédente et celle de la proposition \ref{prop:prolongmorp}
dans laquelle il apparaissait un terme supplémentaire dans le membre de
droite. La raison en est --- et nous attirons l'attention du lecteur sur
ce point --- que l'on ne demande à l'égalité \eqref{eq:commuttau} de
n'être satisfaite que pour $x \in M$ et non pas pour tout $x \in 
F_\tau \otimes_{F_0} M$. La semi-linéarité de $\tau$ montre en
réalité que l'identité \eqref{eq:commuttau} est équivalente à
l'égalité suivante entre applications :
$$(g \otimes \id) \circ \tau_M = \tau_M^{\chi_\tau(g)} \circ 
((\tau^{-\chi_\tau(g)} g \tau) \otimes \id)$$
ce qui correspond bien à la formule de la proposition
\ref{prop:prolongmorp}.

On montre comme dans la preuve de cette même proposition que pour tout
$(\varphi,\tau)$-module $M$, les applications $\tau_M^{p^s}$ forment une
suite qui converge vers l'identité. On peut donc définit $\tau_M^a$ pour
$a \in \Z_p$ et la relation \eqref{eq:commuttau} est alors satisfaite
pour tout $g \in G_\infty$ sans la restriction $\chi_\tau(g) \in \N$.
Ceci montre en particulier que $\tau_M^{-1}$ est défini, c'est-à-dire
que $\tau_M : F_\tau \otimes_{F_0} M \to F_\tau \otimes_{F_0} M$ est une
bijection.
Dans la suite, lorsque cela ne prêtera pas à confusion, on notera
simplement $\tau$ à la place de $\tau_M$ et de même, on écrira souvent
$\varphi$ à la place de $\varphi_M$.

Le but de ce paragraphe est de démontrer que la catégorie des
$\F_p$-représentations de dimension finie de $G_K$ est équivalente à la
catégorie des $(\varphi,\tau)$-modules sur $(F_0, F_\tau)$. À ce sujet,
nous aimerions signaler au lecteur que Tavares Ribeiro s'est déjà
intéressé, dans le premier chapitre de sa thèse \cite{ribeiro}, à des
questions semblables. Toutefois, le point de vue que nous adoptons est
un peu différent et en fait plus proche des travaux de Liu ; notamment,
nous avons constamment le souci de garder explicitement la trace du
$\varphi$-module $M$ décrivant l'action du sous-groupe $G_\infty$ alors
que celui-ci n'apparaît pas du tout dans les travaux de Tavares Ribeiro.

\paragraph{Construction d'un foncteur}

On considère, dans un premier temps, une $\F_p$-représentation $T$ de
$G_K$ de dimension finie $d$ et on note $\calM(T) =
\Hom_{\F_p[G_\infty]}(T, F_0^\sep)$ le $\varphi$-module associé par le
théorème \ref{theo:fontaine}. On sait que $\calM(T)$ est de dimension
$d$ sur $F_0$. Pour définir une structure de $(\varphi,\tau)$-module sur
$\calM(T)$, il reste à construire un automorphisme $\tau$ de
$F_\tau \otimes_{F_0} \calM(T)$. On commence par un lemme qui
donne une description alternative de cet espace.

\begin{lemme}
\label{lem:isomextscal}
L'application naturelle $F_\tau \otimes_{F_0} \calM(T) \to
\Hom_{\F_p[H_\infty]} (T, L)$ est un isomorphisme.
\end{lemme}

\begin{proof}
Par la remarque \ref{rem:injFontaine}, on sait que l'application $L
\otimes_{F_\tau} \Hom_{\F_p[H_\infty]} (T, L) \to \Hom_{\F_p} (T,L)$
est injective, et donc que l'espace $\Hom_{\F_p[H_\infty]} (T, L)$ est
de dimension au plus $d$ sur $F_\tau$. Il suffit donc de montrer
que le morphisme du lemme est injectif. Or, celui-ci s'écrit 
$\beta \circ (F_\tau \otimes_{(F_0^\sep)^{H_\infty}} \alpha)$ où 
$\alpha$ et $\beta$ sont les morphismes canoniques suivants :
$$\begin{array}{rl}
\alpha \, : &
(F_0^\sep)^{H_\infty} \otimes_{F_0} \Hom_{\F_p[G_\infty]} (T, F_0^\sep)
\longrightarrow \Hom_{\F_p[H_\infty]} (T, F_0^\sep) \medskip \\
\beta \, : &
F_\tau \otimes_{(F_0^\sep)^{H_\infty}} 
\Hom_{\F_p[H_\infty]} (T, F_0^\sep) \longrightarrow
\Hom_{\F_p[H_\infty]} (T, L).
\end{array}$$
Il suffit donc de montrer que $\alpha$ et $\beta$ sont injectifs. Pour 
$\alpha$, cela résulte du diagramme commutatif 
$$\xymatrix {
(F_0^\sep)^{H_\infty} \otimes_{F_0} \Hom_{\F_p[G_\infty]} (T, F_0^\sep)
\ar[r]^-{\alpha} \ar@{^(->}[d] & \Hom_{\F_p[H_\infty]} (T, F_0^\sep)
\ar@{^(->}[d] \\
F_0^\sep \otimes_{F_0} \Hom_{\F_p[G_\infty]} (T, F_0^\sep)
\ar[r] & \Hom_{\F_p} (T, F_0^\sep) }$$
et de la proposition \ref{prop:fontaine} qui affirme que la flèche du
bas est un isomorphisme. On en vient maintenant à $\beta$. Il suffit de
montrer que si $f_1, \ldots, f_n$ est une famille d'éléments de
$\Hom_{\F_p[H_\infty]} (T, F_0^\sep)$ qui est liée sur $F_\tau$,
alors elle l'est déjà sur $(F_0^\sep)^{H_\infty}$. Soit
$\lambda_1 f_1 + \cdots + \lambda_n f_n = 0$ une relation de liaison
avec tous les $\lambda_i$ dans $F_\tau$. 
Puisque $T$ est un ensemble fini, il existe un sous-groupe distingué
d'indice fini $H \subset H_\infty$ qui agit trivialement sur $T$, ce 
qui signifie que les fonctions $f_i$ prennent leurs valeurs dans
$(F_0^\sep)^H$. Les extensions $L^H/F_\tau$ et $(F_0^\sep)^H /
(F_0^\sep)^{H_\infty}$ sont alors galoisiennes de groupes de Galois 
isomorphes à $H_\infty/H$. 
On considère une forme linéaire $\ell : L^H \to (F_0^\sep)^H$ qui envoie
$\lambda_1$ sur un élément dont la trace sur $(F_0^\sep)^{H_\infty}$ ne
s'annule pas, et on définit une application $\ell_{H_\infty} : L^H \to
(F_0^\sep)^H$ en moyennant $\ell$ comme suit :
$$\ell_{H_\infty}(x) = \sum_{h \in H/H_\infty} h \, \ell(h^{-1} x).$$
On vérifie sans difficulté que $\ell_{H_\infty}$ est encore une forme
linéaire. Elle est en outre $H_\infty$-équivariante, et donc applique 
$F_\tau$ sur $(F_0^\sep)^{H_\infty}$. Par ailleurs, puisque 
$\lambda_1$ est fixé par $H_\infty$, son image par $\ell_{H_\infty}$
est égale à $\tr_{(F_0^\sep)^H/(F_0^\sep)^{H_\infty}} 
\ell(\lambda_1)$, et n'est donc pas nulle.
Ainsi, en appliquant $\ell_{H_\infty}$ à l'égalité $\lambda_1 f_1 +
\cdots + \lambda_n f_n = 0$, on obtient une relation de liaison non
triviale entre les $f_i$ qui est à coefficients dans
$(F_0^\sep)^{H_\infty}$, ce qui est bien ce que l'on cherchait.
\end{proof}

\begin{rem}
La preuve que l'on vient de donner montre plus généralement que le 
lemme précédent vaut pour un sous-groupe fermé quelconque $H$ de
$G_\infty$ et un sous-corps $L$ de $\Frac R$ contenant $F_0^\sep$ 
et stable seulement par $H$.
\end{rem}

Il n'est maintenant plus difficile de définir l'automorphisme $\tau$. À
cette fin, on note que pour $\sigma \in G_K$ et $f \in
\Hom_{\F_p[H_\infty]}(T,L)$, l'application $g : x \mapsto \sigma
f(\sigma^{-1} x)$ est encore $H_\infty$-équivariante ; la formule
précédente définit donc une action de $G_K$ sur
$\Hom_{\F_p[H_\infty]}(T,L) = F_\tau \otimes_{F_0} \calM(T)$ qui
est, comme on le vérifie directement, semi-linéaire par rapport à la
structure de $F_\tau$-espace vectoriel. De plus, pour cette
action, le sous-groupe $H_\infty$ agit trivialement (ce qui signifie que
l'action se factorise par $G_K/H_\infty$), tandis que le groupe
$G_\infty$, de son côté, agit trivialement sur le sous-ensemble
$\calM(T)$. On définit enfin l'application $\tau$ comme l'automorphisme
de $F_\tau \otimes_{F_0} \calM(T)$ donné par l'action de $\tau$.
Les remarques que l'on vient de faire, combinée au second point du lemme
\ref{lem:strucGK} assurent que l'on obtient bien comme ceci un
$(\varphi,\tau)$-module dans le sens de la définition
\ref{def:phitaumodp}.

\paragraph{Construction d'un quasi-inverse}

On part à présent d'un $(\varphi,\tau)$-module $M$ sur $(F_0, F_\tau)$
et on cherche à lui associer une représentation $\calT(M)$ de $G_K$ qui
soit de dimension finie sur $\F_p$. Bien entendu, en tant que
représentation de $G_\infty$, $\calT(M)$ est la représentation associé
au $\varphi$-module sous-jacent, \emph{i.e.}  $\calT(M) =
\Hom_{F_0,\varphi} (M, F_0^\sep)$. Il reste donc à expliquer comment
étendre cette action à $G_K$ en utilisant l'automorphisme $\tau$. La clé
est le lemme suivant.

\begin{lemme}
\label{lem:isomextscal2}
L'application 
$\calT(M) = \Hom_{F_0,\varphi} (M, F_0^\sep) \to \Hom_{F_\tau, \varphi} 
(F_\tau \otimes_{F_0} M, L)$ déduite de l'extension des
scalaires de $F_0$ à $F_\tau$ est un isomorphisme.
\end{lemme}

\begin{proof}
On fixe une base $(e_1, \ldots, e_d)$ de $M$ et on appelle $G$ l'unique
matrice pour laquelle l'égalité $(\varphi(e_1), \ldots, \varphi(e_d)) =
(e_1, \ldots, e_d) G$ est vérifiée. Se donner un élément de
$\Hom_{F_0,\varphi} (M, F_0^\sep)$ (resp. de $\Hom_{F_\tau,
\varphi}(F_\tau \otimes_{F_0} M, L)$) revient à se donner les images de
$e_i$ qui sont des éléments $x_i \in F_0^\sep$ (resp. $x_i \in L$)
vérifiant le système d'équations $(x_1^p, \ldots, x_d^p) = (x_1, \ldots,
x_d) G$. Or, un tel système a au plus $p^d$ solutions dans n'importe
quel corps, et on sait qu'il en a déjà ce nombre de $F_0^\sep$. Toutes
les solutions dans $L$ sont donc dans $F_0^\sep$ et le lemme est
démontré.
\end{proof}

Étant donné les conditions satisfaites par $\tau$, la proposition
\ref{prop:prolongmorp} s'applique et montre qu'il existe une unique
action de $G_K/H_\infty$ sur le produit tensoriel $F_\tau \otimes_{F_0}
M$ pour laquelle les éléments $g \in G_\infty/H_\infty$ agissent par $(g
\otimes \id)$ et l'élémént $\tau$ agit \emph{via} l'automorphisme
$\tau$. En composant par la projection canonique $G_K \to G_K/H_\infty$,
on obtient une action de $G_K$ sur $F_\tau \otimes_{F_0} M$. Par
ailleurs, $G_K$ agit également sur $L$ et donc sur l'espace $\calT(M) =
\Hom_{F_\tau, \varphi} (F_\tau \otimes_{F_0} M, L)$ \emph{via} la
formule usuelle $\sigma \cdot f : x \mapsto \sigma \, f(\sigma^{-1} x)$.
La forme particulière de l'action de $G_\infty$ sur $F_\tau
\otimes_{F_0} M$ montre immédiatement que l'action qu'on vient de
définir sur $\calT(M)$ prolonge celle de $G_\infty$.

\begin{theo}
\label{theo:equivphitaumodp}
Les deux foncteurs $\calM$ et $\calT$ précédents
induisent des équivalences de catégories inverses l'une de l'autre entre
la catégorie des $\F_p$-représentations de dimension finie de $G_K$ et
la catégorie des $(\varphi,\tau)$-modules sur $(F_0, F_\tau)$.
\end{theo}

\begin{proof}
On sait déjà, par le théorème \ref{theo:fontaine}, que les morphismes
canoniques $M \to \calM(\calT(M))$ et $T \to \calT(\calM(T))$ sont des
isomorphismes pour tout $(\varphi,\tau)$-module $M$ sur $(F_0,
F_\tau)$ et toute $\F_p$-représen\-tation $T$ de dimension finie de
$G_K$. Il ne reste donc qu'à vérifier que le premier commute à 
l'action de $\tau$ tandis que le deuxième est $G_K$-équivariant, ce
qui ne pose aucune difficulté.
\end{proof}

\begin{rem}
\label{rem:Ks}
En considérant non pas l'action de $\tau$ mais celle de $\tau^{p^s}$,
on obtient de la même façon une équivalence de catégories entre, d'une
part, la catégorie des $\F_p$-représentations du groupe $G_s$ et,
d'autre part, la catégorie des $(\varphi,\tau^{p^s})$-modules sur
$(F_0, F_\tau)$ dont les objets sont la donnée de :
\begin{itemize}
\item un $\varphi$-module étale sur $F_0$, noté $M$ ;
\item un automorphisme $\tau^{p^s}$-semi-linéaire $\tau^{(p^s)}_M :
F_\tau \otimes_{F_0} M \to F_\tau \otimes_{F_0} M$ qui
commute à $\varphi$ et tel que, pour tout $g \in G_\infty/H_\infty$
$$\forall x \in M, \quad 
(g \otimes \id) \circ \tau^{(p^s)}_M (x) = (\tau^{(p^s)}_M)^a (x)$$
où $a$ est l'unique élément de $\Z_p$ tel que $\chi(g) \cdot [p^s]
_{\chi(\tau)} = [a]_{\chi(\tau)}$.
\end{itemize}
\end{rem}

\subsubsection{Quelques mots sur le corps $F_\tau$}
\label{subsec:LHinf}

Le corps $F_\tau$ est important car c'est celui qui sert de base à
l'action de $\tau$ sur un $(\varphi,\tau)$-module. Il semble donc
crucial de bien le comprendre. Or, malheureusement, comme nous allons 
le voir dans ce paragraphe, sa structure est loin d'être simple.

\begin{prop}
\label{prop:LH}
Soit $H$ un sous-groupe fermé de $G_\infty$. Alors $L^H$ est l'adhérence
(dans $L$) du perfectisé de $(F_0^\sep)^H$.

De plus, si on note $\m_R$ l'idéal maximal de $R$, la projection
canonique induit, pour tout $x \in R$, un morphisme surjectif $R^H \to
(R/x \m_R)^H$.
\end{prop}

\begin{proof}
La première partie de la proposition s'obtient en mettant ensemble les
deux ingrédients suivants : le corps $F_0^\sep$ est dense dans $L$ (ce
qui est contenu dans la théorie du corps des normes) et le théorème
principal de \cite{ax} qui décrit les points fixes, sous l'action du
groupe de Galois, du complété de la clôture algébrique d'un corps local.
La seconde assertion se démontre de manière analogue en utilisant, à la
place du théorème principal de \cite{ax}, la proposition 2 de ce même
article (p.~424), qui est un peu plus précise.
\end{proof}

La proposition précédente s'applique en particulier avec le groupe
$H_\infty$ et montre donc que $F_\tau$ s'identifie à l'adhérence du
perfectisé du corps $F_\infty = (F_0^\sep)^{H_\infty}$. L'élément
(important) $\eta = 1 - \ueps \in R$, qui est clairement stable par
l'action de $H_\infty$, doit donc s'écrire comme une série faisant
intervenir certains éléments du perfectisé de $F_\infty$. Cependant,
obtenir une telle écriture de façon explicite ne semble pas du tout
facile. 
Voici une autre façon d'appréhender $F_\tau$ (qui montre encore que
décrire les éléments de ce corps est une question délicate). Pour tout
entier $m$, on note $H_m = \Gal(\bar K / K_\infty(\zeta^{p^m}))$ et $F_m
= (F_0^\sep)^{H_m}$. Comme $H_m$ est d'indice fini dans $G_\infty$,
l'extension $F_m/F_0$ est finie. La réunion de ces
extensions, que l'on note $F_\alg$, définit donc un sous-corps de
$F_\tau$ qui est algébrique sur $F_0$. Par ailleurs, on peut considérer
le morphisme $\iota : k[[X,Y]] \to F_\tau$ obtenu en envoyant $X$ sur
$u$ et $Y$ sur $\eta$. La proposition \ref{prop:iotainj} ci-après montre
que le corps des fractions de l'image de $\iota$, noté $k((u,\eta))$,
définit un sous-corps de $F_\tau$ qui est isomorphe à un corps de séries
formelles en deux variables. C'est donc en particulier une extension
purement transcendante de $F_0$. Les sous-corps $F_\alg$ et $k((u,
\eta))$ apparaissent donc, d'un point de vue algébrique, comme deux
constituants \og orthogonaux \fg\ de $F_\tau$. D'un point de vue
topologique, par contre, ces corps semblent s'entremêler de façon 
subtile.

\begin{prop}
\label{prop:iotainj}
Le morphisme $\iota : k[[X,Y]] \to R$, $X \mapsto u$, $Y \mapsto \ueta$ 
est injectif.
\end{prop}

\begin{proof}
Soit une série formelle $F \in k[[X,Y]]$ telle que $F(u, \ueta) =
0$ dans $R$. En faisant agir le groupe de Galois $G_K$, on obtient
l'annulation de $F(u (1 + \ueta)^{c(g)}, (1+\ueta)^{\chi(g)} - 1)$ pour
tout $g \in G_K$. Il en résulte les égalités :
$$F(u + u \ueta^{p^n}, \ueta) =  0 
\quad \text{et} \quad
F(u, \ueta + \ueta^{p^n} + \ueta^{p^n+1}) = 0$$
pour tout entier $n$ suffisamment grand. Pour tout entier $i$, soit
$\partial_X^{[i]}$ l'application $\frac 1 {i!} \frac{\partial^i}
{\partial X^i}$ agissant sur $k[[X,Y]]$ (qui est bien définie). On 
dispose de la formule de Taylor suivante :
\begin{equation}
\label{eq:taylor}
F(u + u \ueta^{p^n}, \ueta) = F(u,\eta) + u \ueta^{p^n} \partial_X
^{[1]} F(u, \eta) + \cdots + u^i \ueta^{i p^n} \partial_X ^{[i]} 
F(u, \eta) + \cdots
\end{equation}
grâce à laquelle on déduit, à partir des annulations précédemment citées, 
que $v_R(\partial_X^{[1]} F(u,\eta)) \geq \frac{ep^{n+1}}{p-1}$. Comme
ceci est vrai pour tout $n$, il vient $\partial_X^{[1]} F(u,\eta) = 0$.
Sachant cela, l'égalité \eqref{eq:taylor} implique maintenant que
$v_R(\partial_X^{[2]} F(u,\eta)) \geq \frac{ep^{n+1}}{p-1}$, et donc
que $\partial_X^{[2]} F(u,\eta)$ s'annule lui aussi. Par récurrence,
on démontre de la même manière que $\partial_X^{[i]} F(u,\eta) = 0$ pour
tout $i$.
Un raisonnement analogue à partir de $F(u, \ueta + \ueta^{p^n} +
\ueta^{p^n+1}) = 0$ montre que $\partial_Y^{[j]} F(u,\eta) = 0$ (avec
des notations évidentes) pour tout $j$. En appliquant cela non pas à la
série $F$ mais à $\partial_X^{[i]} F$, on obtient même $\partial_Y^{[j]}
\partial_X^{[i]} F(u,\eta) = 0$ pour tous $i$ et $j$. En réduisant cette
dernière égalité dans $\bar k$, on s'aperçoit alors que le coefficient
de $X^i Y^j$ dans la série formelle $F$ s'annule lui aussi. Il en 
résulte que la série $F$ est, elle-même, nulle.
\end{proof}

Terminons ce paragraphe en revenant un instant sur les corps $F_m$. 
L'extension $F_1/F_0$ se comprend assez bien, et on sait même la
décrire complètement lorsque le corps $K$ est absolument non ramifié
(\emph{i.e.} si $e = 1$). En effet, dans ce cas, une uniformisante de
$K(\zeta_p)$ est donnée par une racine $(p-1)$-ième de $-p$, que l'on
note $\varpi$. Si $\lambda$ désigne un élément de $k$ tel que $p
\equiv -\pi \lambda \pmod {p^2}$, la suite d'éléments de $\O_{\bar K}/p$
suivante :
$$\Big( \frac 1 {\lambda^{1/p^n}} \Big)^{\! 1 + p + \cdots + p^{n-1}} 
\cdot
\Bigg( \frac{\varpi}{\pi_n^{1 + p + \cdots + p^{n-1}}} \mod p \Bigg)$$
définit un élément $v$ de $R$, qui appartient manifestement à $F_1$. Par
ailleurs, un calcul immédiat montre que $v^{p-1} = \lambda u$. Comme 
$F_1$ est de degré $p-1$ sur $F_0$, il s'ensuit que $F_1 = F_0[v] = F_0
[\sqrt[p-1] {\lambda u}]$. Lorsque $e > 1$, l'extension $F_1/F_0$ est
encore totalement et modérément ramifiée, et son degré est égal à celui
de l'extension $K(\zeta_p)/K$. Par contre, dès que $m \geq 1$, il 
semble bien plus difficile de comprendre l'extension $F_{m+1}/F_m$. On 
sait néanmoins qu'elle est soit triviale, soit de degré $p$. D'après
la théorie d'Artin-Schreier, dans le deuxième cas, elle s'écrit comme
le corps de rupture d'un polynôme de la forme $X^p - X - a_m$ avec
$a_m \in F_m$. L'étude de la ramification sauvage de $F_{m+1}/F_m$
permet d'accéder à la valuation de $a_m$, mais ne permet pas de
répondre à la question plus générale suivante qui nous paraît
intéressante.

\begin{quest}
\label{quest:Fm}
Est-il possible, peut-être seulement sous l'hypothèse $e = 1$, de
décrire explicitement un élément $a_m \in F_m$ tel que $F_{m+1}$ soit le
corps de rupture sur $F_m$ du polynôme $X^p - X - a_m$ ?

Plus généralement étant donné deux entiers $m$ et $m'$ avec $1 \leq m 
< m'$, peut-on décrire l'extension $F_{m'}/F_m$ en termes d'extensions
d'Artin-Schreier-Witt ?
\end{quest}

\noindent
Comme me l'a signalé Berger, la compatibilité entre corps des normes et
corps de classe (démontrée dans \cite{laubie}, \S 3) semble \emph{a
priori} une bonne piste pour étudier ce problème. Les calculs restent,
malgré tout, encore à faire.

\subsubsection{Une variante déperfectisée}
\label{subsec:variante}

Le théorème \ref{theo:equivphitaumodp} que l'on a démontré précédemment
reste valable --- et la démonstration est identique -- si le corps $L$
est remplacé par $L_\np = k((u,\eta))^\sep$ (la clôture séparable de
$k((u,\eta))$ dans $\Frac R$). (Dans la notation, \og $\np$ \fg\
signifie \og non parfait \fg.) En particulier, si l'on pose
$F_{\tau,\np} = L_\np^{H_\infty}$, on a une équivalence de catégories
entre la catégorie des $\F_p$-représentations galoisiennes de $G_K$ et
la catégorie des $(\varphi,\tau)$-modules sur $(F_0, F_{\tau,\np})$.

De surcroît, dans certains cas, le calcul des points fixes de $L_\np$
sous l'action de sous-groupes de $G_\infty$ est plus simple, comme le
montre la proposition suivante.

\begin{prop}
\label{prop:Lnpfixe}
Si $H$ est un sous-groupe d'indice \emph{fini} de $G_\infty$, alors
$L_\np^H = (F_0^\sep)^H$.
\end{prop}

\begin{proof}
Étant donné que $F_0 = k((u))$ est inclus dans $k((u,\eta))$, on a
clairement $(F_0^\sep)^H \subset L_\np^H$. Pour démontrer l'inclusion
réciproque, on considère un élément $x \in L_\np^H$. Comme $x$ est
stable par l'action de $H$, son polynôme minimal sur $k((u,\eta))$ est à
coefficients dans $k((u,\eta))^H$. Soit $v$ une uniformisante de
$(F_0^\sep)^H$ de façon à ce que $(F_0^\sep)^H = k((v))$. On a alors
$k((u,\eta)) \subset k((v,\eta)) \subset k((v))((\eta))$ et $H$ agit
encore sur ce dernier espace en fixant $v$ et en envoyant $\eta$ sur
$(1+\eta)^{\chi(\cdot)} - 1$. Étant donné que $H$ est d'incide fini dans
$G_\infty$, son image par le caractère cyclotomique $\chi$ n'est pas
triviale et on a donc $k((v))((\eta))^H = k((v))$. On en déduit que
$k((u,\eta))^H = k((v))^H = F_0$ et, donc, que le polynôme minimal de
$x$ est à coefficients dans $F_0$, d'où on déduit que $x$ est élément de
$F_0^\sep$. Comme $x$ est en outre fixe par $H$, on a bien $x \in
(F_0^\sep)^H$.
\end{proof}

On prendra garde au fait que, dans la proposition ci-dessus, l'hypothèse
\og d'indice fini \fg\ est essentielle. En particulier, la conclusion de
la proposition ne vaut malheureusement pas si $H = H_\infty$ (qui est
pourtant le cas qui nous intéresserait le plus). Il est d'ailleurs
facile de s'en convaincre puisque l'élément $\eta$ appartient
manifestement à $L_\np^{H_\infty}$, sans pour autant être dans
$(F_0^\sep)^{H_\infty}$ car, comme cela a été vu, il n'est pas
algébrique sur $F_0$.

\begin{rem}
De la même façon, on peut, à la place de $L$ ou de $L_\np$, utiliser des
versions partiellement déperfectisée comme le corps $L_{u\tnp} =
k((u,\eta^{1/p^\infty}))$ ou $L_{\eta\tnp} = k((u^{1/p^\infty},\eta))$.
Le cas de $L_{u\tnp}$ sera utile dans la suite.
\end{rem}

\subsection{Relèvement modulo $p^n$ et en caractéristique nulle}
\label{subsec:relevement}

Le théorème \ref{introtheo:equivphitau} de l'introduction se déduit enfin 
du théorème \ref{theo:equivphitaumodp} pour $L = F_\infty$ à l'aide
d'arguments classiques de dévissage. Ce sont ces arguments que nous nous
proposons de présenter dans cette partie. Afin surtout de mettre en
place les notations, nous commençons par rappeler brièvement comment
ceux-ci fonctionnent dans le cas classique des $\varphi$-modules et des
représentations de $G_\infty$.

\subsubsection{La théorie de Fontaine modulo $p^n$ et en caractéristique
nulle}

L'idée de base consiste à remplacer l'anneau $R$ par l'anneau des
vecteurs de Witt $W(\Frac R)$ ; c'est naturellement une $W$-algèbre
munie d'une action de $G_K$.
L'anneau $\E^\ent$ défini dans l'introduction par, rappelons-le :
$$\E^\ent = \Big\{ \, \sum_{i \in \Z} a_i u^i \quad \big| \quad 
a_i \in W, \,\, \lim_{i \to -\infty} a_i = 0 \, \Big\}$$
se plonge dans $W(\Frac R)$ en envoyant $u$ sur le représentant de
Teichmüller $[\upi]$. Muni de la valuation $p$-adique, $v_\E (\sum_{i
\in \Z} a_i u^i) = \min_{i \in \Z} v_p(a_i)$, $\E^\ent$ est un anneau de
valuation discrète qui admet $F_0$ pour corps résiduel. On note
${\E^{\ent,\ur}}$ l'unique sous-algèbre étale (infinie) de $W(\Frac R)$ ayant
pour corps résiduel $F_0^\sep \subset \Frac R$. Si l'on pose $\E = \Frac
\E^\ent$ et $\E^\ur = \Frac {\E^{\ent,\ur}}$, le groupe de Galois de
l'extension $\E^\ur/\E$ s'identifie à celui de l'extension
résiduelle $F_0^\sep/F_0$ et donc finalement à $G_\infty$. 
L'anneau $W(\Frac R)[1/p]$ est naturellement muni d'un opérateur de
Frobenius et celui-ci définit par restriction des endomorphismes de 
$\E^\ent$, ${\E^{\ent,\ur}}$, $\E$ et $\E^\ur$. Sur $\E^\ent$, par exemple, on
voit aisément qu'il agit en appliquant le Frobenius traditionnel aux
coefficients et en envoyant $u$ sur $u^p$.  On définit un 
\emph{$\varphi$-module étale} sur $\E^\ent$ (resp. sur $\E$) comme la donnée
d'un $\E^\ent$-module de type fini (resp. d'un $\E$-espace vectoriel de
dimension finie) $M$ munie d'une application $\varphi : M \to M$
semi-linéaire par rapport au Frobenius et dont l'image engendre $M$.

\begin{theo}
Les foncteurs suivants sont des équivalences de catégories
quasi-inverses l'une de l'autre :
\begin{eqnarray*}
\left\{ \begin{array}{c}
\Q_p\text{-représentations de }\\
\text{dimension finie de }G_\infty 
\end{array} \right\}
& \stackrel{\sim}{\longrightarrow} & 
\left\{ \begin{array}{c}
\varphi\text{-modules étales sur } \E
\end{array} \right\}
\\
T & \mapsto & \Hom_{\Q_p[G_\infty]}(T, \E^\ur) \\
\Hom_{\E, \varphi} (M, \E^\ur) & \mapsfrom & M
\end{eqnarray*}
\end{theo}

On dispose d'énoncés analogues pour les $\Z_p$-représentations libres
d'une part, et annulées par $p^n$ d'autre part, que nous laissons au
lecteur le soin d'écrire. On fait remarquer quand même qu'afin d'obtenir
ces énoncés, il convient de remplacer l'anneau de périodes $\E^\ur$ par
${\E^{\ent,\ur}}$ et ${\E^{\ent,\ur}}/p^n{\E^{\ent,\ur}}$
respectivement. On peut également écrire une équivalence de catégories
mettant en jeu à gauche la catégorie de toutes les
$\Z_p$-représentations annulées par une puissance de $p$ (non précisée) ; 
dans le cas, l'anneau de périodes $\E^\ur$ doit être remplacé par le
produit tensoriel $\Q_p/\Z_p \otimes_{\Z_p} {\E^{\ent,\ur}}$ (qui n'est
pas un anneau, mais seulement un ${\E^{\ent,\ur}}$-module) ou, ce qui
revient au même, par le quotient $\E^\ur / {\E^{\ent,\ur}}$.

\subsubsection{Définition générale des $(\varphi,\tau)$-modules}
\label{subsec:defphitau}

On pose $\E^\ent_\tau = W(F_\tau)$ et $\E_\tau = \Frac \E^\ent_\tau$. Le
corps $F_\tau$ étant parfait, la valuation $p$-adique fait de
$\E^\ent_\tau$ un anneau de valuation discrète, complet, dont le corps
résiduel s'identifie à $F_\tau$. En outre, $\E_\tau$ s'obtient à partir
de $\E^\ent_\tau$ simplement en inversant $p$. Tous les anneaux que l'on
vient de définir sont munis d'un endomorphisme de Frobenius que l'on
note $\varphi$ ou $\varphi_A$ ($A$ étant l'anneau sur lequel le
Frobenius agit) dans les cas où il sera important de le préciser. La
définition des $(\varphi,\tau)$-modules est désormais la même qu'en
caractéristique $p$ (voir définition \ref{def:phitaumodp}) à part que
les anneaux de base sont ceux que l'on vient de définir. On la redonne
ci-dessous pour plus de clarté.

\begin{deftn}
\label{def:phitau}
Un \emph{$(\varphi,\tau)$-module sur $(\E^\ent, {\E^\ent_\tau})$} (resp.
$(\E, \E_\tau)$) est la donnée de
\begin{itemize}
\item un $\varphi$-module étale sur $\E^\ent$ (resp. sur $\E$), noté $M$ ;
\item un endomorphisme $\tau$-semi-linéaire $\tau_M : {\E^\ent_\tau} 
\otimes_{\E^\ent} M \to {\E^\ent_\tau} \otimes_{\E^\ent} M$ qui commute à 
$\varphi_{{\E^\ent_\tau}} \otimes \varphi_M$ et
vérifie, pour tout $g \in G_\infty/H_\infty$ tel que $\chi_\tau(g) \in
\N$, la relation suivante :
\begin{equation}
\label{eq:commuttau2}
\forall x \in M, \quad
(g \otimes \id) \circ \tau_M (x) = \tau_M^{\chi_\tau(g)} (x).
\end{equation}
\end{itemize}
\end{deftn}

On souhaite à présent démontrer le théorème \ref{introtheo:equivphitau}
de l'introduction, c'est-à-dire que les catégories de
$(\varphi,\tau)$-modules sont équivalentes aux catégories
correspondantes de représentations galoisiennes. L'étape essentielle
pour cela consiste à étendre les lemmes \ref{lem:isomextscal} et
\ref{lem:isomextscal2} (qui constituaient la clé de la démonstration
dans le cas de caractéristique $p$) à la nouvelle situation relevée.
À partir de maintenant, on supposera toujours implicitement que les
$\Q_p$-représentations (resp. $\Z_p$-représentations) considérées sont
de dimension finie sur $\Q_p$ (resp. de type fini comme $\Z_p$-module).
Si $T$ est une telle représentation du groupe $G_\infty$, on note
$\calM(T)$ le $\varphi$-module étale sur $\E^\ent$ ou sur $\E$ qui lui
est associé par la théorie de Fontaine. De même, si $M$ est un
$\varphi$-module étale défini sur $\E^\ent$ ou sur $\E$, on note
$\calT(M)$ la représentation $p$-adique qui lui correspond.

\begin{lemme}
\label{lem:relpn}
Pour toute $\Z_p$-représentation de torsion (resp. $\Z_p$-représentation 
libre, resp. $\Q_p$-représentation) $T$ de $G_\infty$, l'application 
naturelle 
$$\begin{array}{c}
{\E^\ent_\tau} \otimes_{\E^\ent} \calM(T) \to \Hom_{\Z_p[H_\infty]} 
(T, \Q_p/\Z_p \otimes_{\Z_p} W(L))
\medskip \\
\text{(resp. }
{\E^\ent_\tau} \otimes_{\E^\ent} \calM(T) \to \Hom_{\Z_p[H_\infty]} 
(T, W(L)),
\medskip \\
\text{resp. } \E_\tau \otimes_{\E} \calM(T) \to
\Hom_{\Q_p[H_\infty]} (T, W(L)[1/p]) \text{)}
\end{array}$$
est un isomorphisme.

Pour tout $\varphi$-module étale $M$ défini sur $\O_\E$ et de 
torsion (resp. défini sur $\O_\E$ et libre, resp. défini sur
$\E$), l'application naturelle
$$\begin{array}{c}
\calT(M) \to \Hom_{{\E^\ent_\tau}, \varphi} 
({\E^\ent_\tau} \otimes_{\E^\ent} M, \Q_p/\Z_p \otimes_{\Z_p} W(L)) 
\medskip \\ 
\text{(resp. }
\calT(M) \to \Hom_{{\E^\ent_\tau}, \varphi} 
({\E^\ent_\tau} \otimes_{\E^\ent} M, W(L)),
\medskip \\
\text{resp. } \calT(M) \to \Hom_{\E_\tau, \varphi} 
(\E_\tau \otimes_{\E} M, W(L)[1/p]) \text{)}
\end{array}$$
est un isomorphisme.
\end{lemme}

\begin{proof}
On ne démontre que la première partie du lemme, la seconde étant
totalement analogue. On prouve tout d'abord le résultat lorsque $T$
est une $\Z_p$-représentation annulée par $p^n$. On raisonne par
récurrence sur $n$. Pour $n = 1$, le résultat à démontrer est exactement 
l'assertion du lemme \ref{lem:isomextscal} ; il n'y a donc plus rien à
faire. Pour passer de $n$ à $n+1$, on considère $T$ une représentation
annulée par $p^{n+1}$. Elle s'insère dans la suite exacte $0 \to pT
\to T \to T/pT \to 0$ qui donne naissance au diagramme suivant :
$$\xymatrix  @C=20pt {
0 \ar[r] &
{\E^\ent_\tau} \otimes_{\E^\ent} \calM(T/pT) \ar[d] \ar[r] &
{\E^\ent_\tau} \otimes_{\E^\ent} \calM(T) \ar[d] \ar[r] &
{\E^\ent_\tau} \otimes_{\E^\ent} \calM(pT) \ar[d] \ar[r] & 0 \\
0 \ar[r] &
\Hom_{\Z_p[H_\infty]} (T/pT, CW(L)) \ar[r] &
\Hom_{\Z_p[H_\infty]} (T, CW(L)) \ar[r] &
\Hom_{\Z_p[H_\infty]} (pT, CW(L)) }$$
où $CW(L) = \Q_p/\Z_p \otimes_{\Z_p} W(L)$.
La ligne du haut est exacte car, d'une part, le foncteur $\calM$ l'est
et, d'autre part, ${\E^\ent_\tau}$ est plat sur $\E^\ent$. Les flèches
horizontales de gauche et de droite sont des isomorphismes par hypothèse
de récurrence, et finalement la ligne du bas est exacte par exactitude à
gauche du foncteur $\Hom$. Une chasse au diagramme montre alors que la
flèche verticale centrale est aussi un isomorphisme.

Le cas où $T$ est une $\Z_p$-représentation quelconque (toujours
supposée de type fini comme $\Z_p$-module) s'obtient alors par passage à
la limite, tandis que celui où $T$ est une $\Q_p$-représentation s'en
déduit en inversant $p$.
\end{proof}

Il nous reste à définir des foncteurs dans les deux sens entre la
catégorie des $\Z_p$-représentations (resp. $\Q_p$-représentations) de
$G_K$ et celle des $(\varphi,\tau)$-modules sur $(\E^\ent,
{\E^\ent_\tau})$ (resp. sur $(\E, \E_\tau)$) puis à montrer que ceux-ci
réalisent des équivalences de catégories inverses l'une de l'autre. Pour
cela, à la lumière du lemme précédent, il suffit de reprendre presque
\emph{verbatim} les constructions du \S \ref{subsec:equivmodp}, ce que
nous laissons en exercice au lecteur. On réécrit toutefois explicitement
les deux propriétés essentielles qui façonnent la construction de ces
foncteurs : \emph{primo}, le $\varphi$-module sous-jacent au
$(\varphi,\tau)$-module associé à une représentation $T$ est $\calM(T)$
et \emph{secundo}, dans le cas où $T$ est définie sur $\Z_p$ par
exemple, l'action de $\tau$ sur ${\E^\ent_\tau} \otimes_{\E^\ent}
\calM(T)$ provient \emph{via} le lemme \ref{lem:relpn} de son action 
naturelle sur l'espace $\Hom_{\Z_p[H_\infty]} (T, W(L))$.

\subsubsection{Variante déperfectisée}
\label{subsec:variante2}

De même que, comme nous l'avons expliqué au \S \ref{subsec:variante}, il
était possible en caractéristique $p$ de remplacer le corps $L$ par le
corps plus petit $L_\np = k((u,\eta))^\sep$, on peut, dans la situation
relevée considérée ici, remplacer $W(L)$ par un anneau plus petit qui
n'est pas parfait. Pour construire ce remplaçant, on définit en premier
lieu l'anneau $\calF_\np^\ent$ comme le complété $p$-adique du localisé
de $W[[u,\eta]]$ en l'idéal premier $p$ (cet idéal est bien premier car
le quotient $W[[u,\eta]] / p W[[u,\eta]]$ s'identifie à $k[[u,\eta]]$
qui est manifestement intègre). Clairement, $\calF_\np^\ent$ est un
anneau de Cohen de $k((u,\eta))$. On définit alors
$\calF_\np^{\ent,\sep}$ comme le complété $p$-adique de l'unique
extension étale infinie de $\calF_\np^\ent$ incluse dans $W(L)$ dont le
corps résiduel s'identifie à $L_\np$ ; c'est notre subtitut à $W(L)$. On
pose $\E^\ent_{\np,\tau} = (\calF_\np^{\ent,\sep}) ^{H_\infty}$ et
$\E_{\np,\tau} = \E^\ent_{\np,\tau}[1/p] = \Frac \E^\ent_{\np,\tau}$.

Les arguments du \S \ref{subsec:defphitau} se généralisent directement à
cette nouvelle situation ; on en déduit que le théorème
\ref{introtheo:equivphitau} de l'introduction est encore valable si on
remplace $\E^\ent_\tau$ et $\E_\tau$ par $\E^\ent_{\np,\tau}$ et
$\E_{\np,\tau}$ respectivement. En d'autres termes, l'opérateur $\tau$
n'est pas uniquement défini sur $\E^\ent_\tau$ (resp. sur $\E_\tau$)
mais sur l'anneau plus petit $\E^\ent_{\np,\tau}$ (resp.
$\E_{\np,\tau}$).

\medskip

De la même façon qu'au \S \ref{subsec:variante}, on peut aussi
introduire des variantes partiellement déperfectisées de $W(L)$,
$\E_\tau^\ent$ et $\E_\tau$. Par exemple, en autorisant par exemple les
racines $p^n$-ièmes de $\eta$ (mais pas celles de $u$), on obtient par
comme ceci les anneaux $\calF_{u\tnp}^\ent$,
$\calF_{u\tnp}^{\ent,\sep}$, $\E^\ent_{u\tnp,\tau}$ et $\E_{u\tnp,\tau}$
qui joueront une rôle dans la suite de cet article (voir la
démonstration de la proposition \ref{prop:existres}).

\subsubsection{Deux exemples}
\label{subsec:exemples}

Le premier exemple que nous aimerions présenter est celui d'une
représentation $T$ (définie au choix sur $\F_p$, $\Z_p$ ou $\Q_p$) dont
la restriction au sous-groupe $G_\infty$ est triviale, c'est-à-dire dont
le $\varphi$-module correspondant est trivial. Cette hypothèse est en
réalité très restrictive car, si l'action du sous-groupe $G_\infty$ est
triviale, il en est nécessairement de même de tous ses conjugués. Or, si
$s$ désigne le plus grand entier pour lequel le corps $K$ admet une
racine primitive $p^s$-ième de l'unité, les conjugués de $G_\infty$
engendrent ensemble le sous-groupe (distingué) d'indice fini $G_s$. En
particulier, si $K$ ne contient pas de racine primitive $p$-ième de
l'unité (par exemple si son indice de ramification absolu $e$ est
strictement plus petit que $p-1$), une représentation de $G_K$, dont la
restriction à $G_\infty$ est triviale, est, elle même, triviale. Dans le
cas général, l'action se factorise par le quotient $G_K/G_s$ qui est un
groupe cyclique de cardinal $p^s$ engendré par l'image de $\tau$. 

Autrement dit, se donner une représentation $T$ dont la restriction à
$G_\infty$ est triviale, revient à se donner un automorphisme $\tau$ de
$T$ d'ordre $p^s$. Le $(\varphi,\tau)$-module $M$ associé à $T$ se
décrit alors comme suit (la vérification est immédiate et laissée au
lecteur) : on a $M = \E^\ent \otimes_{\Z_p} T$ et l'automorphisme
$\tau_M$ sur $\E_\tau^\ent \otimes_{\E^\ent} M = \E_\tau^\ent
\otimes_{\Z_p} T$ est $\tau \otimes \tau$.

\bigskip

Venons-en maintenant à notre second exemple, qui est celui du caractère
cyclotomique, ou plus généralement d'une de ses puissances. En réalité,
cet exemple a été traité par Liu dans son article \cite{liu2} (voir
exemple 3.2.3) au cours de son étude des réseaux dans les
représentations semi-stables. Comme le calcul n'est pas évident et
demande de connaître un peu de théorie de Hodge $p$-adique, nous
préférons nous contenter ici de donner le résultat en renvoyant à
\emph{loc. cit.} pour la preuve.

Soit $\t \in W(R)$ un élément non divisible par $p$ vérifiant
$\varphi(\t) = c^{-1} E(u) \t$ où $c = \frac{E(0)} p$ est un élément
inversible dans $\Z_p$ (on rappelle, à toutes fins utiles, que $E(u)$
désigne le polynôme minimal sur $W[1/p]$ de l'uniformisante $\pi$
choisie, et que c'est donc un polynôme d'Eisenstein). L'existence d'un
tel élément $\t$ découle du calcul de \emph{loc. cit.} mais peut aussi se
voir comme la conséquence du lemme \ref{lem:existV} qui sera démontré
dans la suite\footnote{Le même lemme assure également que la condition
$\varphi(\t) = c^{-1} E(u) \t$ détermine $\t$ à multiplication près par un
élément inversible de $\Z_p$.}.
Avec ces notations, le $(\varphi,\tau)$-module associé à la
représentation $\Z_p(n)$ (avec $n \in \Z$), c'est-à-dire la
représentation $T = \Z_p w$ où l'action de Galois est donnée par $g w =
\chi(g)^n w$, est engendré par la fonction $f : T \to \E^{\ent,\ur}$, $w
\mapsto \t^n$. Concrètement, il est décrit par les formules suivantes :
$$M = \E^\ent \cdot f
\quad ; \quad
\varphi(f) = c^{-n} E(u)^n \cdot f
\quad \text{et} \quad
\tau(f) = \Big(\frac{\tau(\t)} \t \Big)^n \cdot f.$$
Les éléments $\t$, $\tau(\t)$ et $E(u)$ sont inversibles dans $\E^\ent$,
de sorte que les égalités précédentes ont bien un sens, même lorsque
$n$ est négatif. Les $(\varphi,\tau)$-modules correspondant à $\F_p(n)$
et $\Q_p(n)$ sont donnés par des formules analogues.

\section{Réseaux dans les $(\varphi,\tau)$-modules}
\label{sec:reseaux}

Nous avons démontré dans la section précédente que la catégorie des
$\Z_p$-représentations galoisiennes de $G_K$ est équivalente à celle des
$(\varphi,\tau)$-modules sur $(\E^\ent, {\E^\ent_\tau})$. Ce résultat
peut paraître satisfaisant, mais il se heurte néanmoins à un problème
pratique important lié au fait que l'anneau ${\E^\ent_\tau}$ est
difficile à manipuler concrètement. On peut certes écrire ses éléments
comme des séries, mais celles-ci requièrent une infinité de variables et
des conditions de convergence subtiles.

Dans cette partie, nous aimerions montrer en quoi l'introduction de
réseaux à l'intérieur des $(\varphi,\tau)$-modules permet d'apporter des
éléments de réponse au problème précédent. Nous commençons par
introduire la notion de réseau dans les $\varphi$-modules étales dans le
\S \ref{subsec:phireseaux}, puis montrons dans le \S \ref{subsec:bornes}
comment celle-ci peut être utilisée pour établir des bornes explicites
portant sur la ramification des représentations galoisiennes. Forts de
ces résultats préliminaires, nous en arrivons ensuite, dans le \S
\ref{subsec:phitaures}, au c\oe ur de notre problème en introduisant la
notion de $(\varphi,\tau)$-réseau puis en démontrant le théorème
\ref{theo:congr} qui donne des contraintes fortes sur la forme des
éléments de $\E^\ent_\tau$ --- et notamment de leur représentation sous
forme de séries --- qui interviennent dans l'expression de
l'automorphisme $\tau$ (par exemple sous forme matricielle dans le cas
d'un module libre).

\subsection{Réseaux dans les $\varphi$-modules étales}
\label{subsec:phireseaux}

\subsubsection{Définitions et rappels}

Comme dans \cite{kisin}, on pose $\Sk = W[[u]]$. Cet anneau se plonge
naturellement dans $\E^\ent$ et dans $W(R)$ en envoyant comme d'habitude
$u$ sur le représentant de Teichmüller de $\upi$. En particulier, $\Sk$
apparaît comme un sous-anneau de l'intersection $\E^\ent \cap W(R)$, et
on démontre en fait facilement qu'il s'identifie à cette intersection.
Il est clair par ailleurs que le Frobenius (agissant sur $W(L)$ par
exemple) stabilise $\Sk$, et que ce dernier anneau est également stable
par l'action de $G_\infty$.

\begin{deftn}
\label{def:phires}
Soit $M$ un $\varphi$-module étale défini sur $\E^\ent$. Un
\emph{$\varphi$-réseau} dans $M$ est la donnée d'un sous-$\Sk$-module de
type fini $\frakM$ de $M$ qui est stable par $\varphi$ et qui est tel
que $\E^\ent \otimes_{\Sk} \frakM = M$.
\end{deftn}

\begin{rem}
Dans le cas où $M$ est un module libre sur $\E^\ent$, on se restreindra
souvent aux $\varphi$-réseaux qui sont eux-même libres (de type fini et
même rang) sur $\Sk$. Ce n'est en fait pas une véritable restriction car
il suit du théorème de structure des modules sur $\Sk$ (voir, par exemple,
théorème 3.1, chap. 5 de \cite{lang} pour un énoncé de ce théorème) que
si $\frakM$ est un $\varphi$-réseau quelconque dans $M$, alors $(\E^\ent
\otimes_\Sk \frakM) \cap \frakM[1/p]$ est un $\varphi$-réseau libre.
\end{rem}

Voici une autre définition importante qui est essentiellement dûe à
Fontaine et qui, en un certain sens, mesure la complexité d'un réseau.

\begin{deftn}
Soit $\frakM$ un $\varphi$-réseau à l'intérieur d'un $\varphi$-module
étale sur $\E^\ent$. Soit encore $U$ un élément de $W(R)$ et $h$ un
nombre entier positif ou nul.

On dit que $\frakM$ est de \emph{hauteur} divisant $U$ (resp. de
$U$-hauteur $\leq h$ pour un certain entier $h$) si le conoyau de
l'application $\id \otimes \varphi : W(R) \otimes_{\varphi,\Sk} \frakM
\to W(R) \otimes_{\Sk} \frakM$ est annulé par $U$ (resp.
$U^h$)\footnote{On notera en particulier que les locutions \og de
$U$-hauteur $\leq 1$ \fg\ et \og de hauteur divisant $U$ \fg\ sont
synonymes.}.

Finalement, on dit que $\frakM$ est de $U$-hauteur finie s'il est de
$U$-hauteur $\leq h$ pour un certain entier $h$.
\end{deftn}

\begin{rem}
La définition précédente de la hauteur a, en réalité, surtout un intérêt
lorsque $U \in \Sk$, auquel cas on peut se contenter de vérifier que le
conoyau de $\id \otimes \varphi : \Sk \otimes_{\varphi,\Sk} \frakM \to
\frakM$ est annulé par $U$ (sans aller jusqu'à étendre les scalaires à
$W(R)$). Toutefois, lorsque, dans la suite, nous manipulerons les
$(\varphi,\tau)$-réseaux, nous aurons besoin à plusieurs reprises de
considérer des $U \not\in \Sk$, et c'est pourquoi nous avons préféré 
donner directement la définition générale précédente.
\end{rem}

Dans le cas des $\varphi$-modules de $p$-torsion, on a le lemme suivant
qui dénote un comportement exemplaire des réseaux.

\begin{lemme}
\label{lem:resmodpn}
Tout $\varphi$-module étale sur $\E^\ent$ annulé par une puissance de 
$p$ admet un $\varphi$-réseau.

Tout $\varphi$-réseau dans un $\varphi$-module étale sur $\E^\ent$ qui
est annulé par une puissance de $p$ est de $u$-hauteur finie.
\end{lemme}

\begin{proof}
La première assertion résulte de la remarque suivante qui découle
directement de la définition de $\E^\ent$ : si $\frakM$ est un réseau
quelconque dans un $\varphi$-module étale sur $\E^\ent$ annulé par une
puissance de $p$, alors il existe un entier $n$ tel que $u^n \frakM$
soit stable par $\varphi$.

La seconde assertion, quant à elle, peut se voir comme une conséquence
du théorème de classification des modules sur $\Sk = W[[u]]$. En effet,
étant donné un $(\varphi,\tau)$-réseau $\frakM$ comme dans la
définition, le théorème en question assure que le conoyau de $\id
\otimes \varphi : \Sk \otimes_{\varphi,\Sk} \frakM \to \frakM$ est de
longueur finie comme $\Sk$-module. Il suffit alors pour conclure de
remarquer que $u$ n'est pas inversible dans $\Sk$.
\end{proof}

Les résultats du lemme précédent ne s'étendent pas au cas des
$\varphi$-modules libres sur $\E^\ent$. Dans le \S \ref{subsec:critere}
ci-après, nous examinerons l'équivalent de la première partie du lemme.
En ce qui concerne la deuxième partie du lemme, nous attirons
l'attention du lecteur sur le fait trivial suivant : tout
$\varphi$-réseau vivant dans un $\varphi$-module étale libre est de
hauteur divisant $U$ pour un certain $U \in \Sk$ (il suffit de prendre
pour $U$ le déterminant de $\varphi$ agissant sur le réseau). Par
contre, il n'est pas vrai que $U$ peut toujours être choisi de la forme
$u^n$ pour un certain entier $n$.

\subsubsection{Un critère pour l'existence de $\varphi$-réseaux}
\label{subsec:critere}

Soit $M$ un $\varphi$-réseau étale libre sur $\E^\ent$. Soit $T$ la
$\Z_p$-représentation de $G_\infty$ qui lui correspond. Les données $T$
et $M$ sont alors liées par la formule $M = \Hom_{\Z_p[G_\infty]}(T,
\E^{\ent,\ur})$. Une structure entière naturelle à l'intérieur de $M$
(et donc un candidat potentiel pour être un réseau) est l'espace $\frakM
= \Hom_{\Z_p[G_\infty]}(T, \Sk^\ur)$ où $\Sk^\ur = W(R) \cap
\E^{\ent,\ur}$.

\begin{prop}
\label{prop:critereres}
En reprenant les notations précédentes, les conditions suivantes sont 
équivalentes :
\begin{enumerate}
\item le sous-module $\frakM$ de $M$ est un $\varphi$-réseau ;
\item le $\varphi$-module $M$ admet un $\varphi$-réseau qui est 
libre sur $\Sk$ ;
\item le $\varphi$-module $M$ admet un $\varphi$-réseau.
\end{enumerate}
\end{prop}

\begin{proof}
Le lemme 2.1.10 de \cite{kisin} montre que la première condition
implique la deuxième. Comme cette dernière implique clairement la
troisième, il suffit de montrer que la troisième condition implique la
première.
Soit $\frakM'$ un $\varphi$-réseau dans $M$ (qui existe, par hypothèse).
La première étape consiste à démontrer que $\frakM'$ est inclus dans
$\frakM$. On considère les éléments de $\frakM'$ comme des morphismes
$\Z_p$-linéaires et $G_\infty$-équivariants de $T$ dans $\E^{\ent,\ur}$.
Soit $X$ le sous-$\Sk$-module de $\E^{\ent,\ur}$ engendré par les images
des éléments de $\frakM'$. Il est stable par $\varphi$ (puisque
$\frakM'$ l'est) et de type fini sur $\Sk$ (puisque $\frakM'$ l'est). On
en déduit qu'il est inclus dans $W(R)$, ce qui est bien ce que l'on
avait annoncé. La conclusion s'obtient maintenant facilement. En effet,
de l'inclusion $\frakM' \subset \frakM$, on déduit que $\E^{\ent,\ur}
\otimes_\Sk \frakM$ contient $\E^{\ent,\ur} \otimes_\Sk \frakM' = M$
et, par suite, que $\frakM$ est un $\varphi$-réseau de $M$.
\end{proof}

La proposition précédente permet en particulier de montrer que certaines
représentations $T$ correspondent à des $\varphi$-modules n'admettant
pas de $\varphi$-réseau. C'est par exemple le cas de la représentation
$\Z_p(-1)$, comme nous nous proposons de le vérifier pour conclure ce
numéro. Soit $w$ un générateur de $\Z_p(-1)$. D'après le deuxième
exemple traité dans le \S \ref{subsec:exemples}, le $\varphi$-module
associé $M$ est engendré par la fonction $f : \Z_p(-1) \to
\E^{\ent,\ur}$, $w \mapsto V^{-1}$ où $V$ vérifie l'équation $\varphi(V)
= c^{-1} E(u) V$ et où, dans cette dernière égalité, $E(u)$ désigne le
polynôme minimal de l'uniformisante $\pi$ et $c = \frac{E(0)} p \in
\Z_p^\times$. D'après la proposition \ref{prop:critereres}, pour montrer
que $M$ n'admet pas de $\varphi$-réseau, il suffit de montrer qu'aucun
élément non nul de $M$ ne prend ses valeurs dans $W(R)$. Autrement dit,
il suffit de montrer que $\E^\ent \cap V W(R)$ est réduit à $0$. 

On pose $\Sk' = \E^\ent \cap V W(R)$. De $V W(R) \subset W(R)$, on
déduit facilement que $\Sk' \subset \Sk$. Il est clair par ailleurs que
cette inclusion induit un morphisme injectif $\Sk'/p\Sk' \to \Sk/p\Sk =
k[[u]]$. On en déduit que $\Sk'/p\Sk'$ est un $k[[u]]$-module libre de
rang $\leq 1$, et donc que $\Sk'$ est un $\Sk$-module libre de rang
$\leq 1$. Soit $a$ un élément de $\Sk$, éventuellement nul, qui engendre
$\Sk'$. D'après une variante du théorème de préparation de Weierstrass
(voir, par exemple, théorème 2.1, chap. 5 de \cite{lang}), on peut
supposer que $a$ est un polynôme ; on le note à partir de maintenant
$A(u)$. Soit $A^\sigma$ le polynôme obtenu à partir de $A$ en appliquant
le Frobenius $\sigma$ à chacun de ses coefficients. On a $\varphi(A(u))
= A^\sigma(u^p)$. D'autre part, on vérifie que le quotient
$\frac{\varphi(A(u))}{E(u)}$ appartient encore à $\Sk'$. Il en résulte
que $A(u) E(u)$ divise $A^\sigma(u^p)$ dans $\Sk$. Les éléments de $\Sk$
définissant des séries convergentes sur le disque de centre $0$ et de
rayon $1$, la divisibilité trouvée implique que le polynôme $A^\sigma$
s'annule en $\pi^p$. Par la théorie des polygones de Newton, il en
résulte que $A$ s'annule en un élément de valuation $p$, à partir de
quoi on trouve que $A^\sigma$ admet une racine de valuation $p^2$. Par
récurrence, on démontre que $A^\sigma$ admet une racine de valuation
$p^n$ pour tout entier $n$. Ceci n'est évidemment possible que si $A(u)$
est le polynôme nul, c'est-à-dire si $\Sk' = 0$. On a donc finalement
bien demontré ce que l'on souhaitait.

\subsubsection{Calcul de la représentation galoisienne associée}

On fixe un $\varphi$-module $M$ sur $\E^\ent$, ainsi qu'un $\varphi$-réseau
$\frakM$ à l'intérieur de $M$. On suppose que l'on est dans
l'alternative suivante : soit $M$ est annulé par $p^n$ pour un certain
$n$, soit $M$ est libre comme module sur $\E^\ent$. Dans le deuxième cas,
on pose $n = \infty$, et on convient que $W_\infty (R) = W(R)$. Soit
$\calT(M)$ la représentation galoisienne associée à $M$. Dans ce
paragraphe, nous donnons plusieurs formules permettant de calculer
$\calT(M)$ directement à partir de $\frakM$. La plus simple consiste
évidemment à commencer par retrouver $M$ en étendant les scalaires à
$\E^\ent$, ce qui conduit à l'expression suivante :
$$\calT(M) = \Hom_{\E^\ent, \varphi} (\E^\ent \otimes_\Sk \frakM,
{\E^{\ent,\ur}} / p^n {\E^{\ent,\ur}}).$$
Toutefois, on aimerait justement éviter cette solution, car l'un des
intérêts d'utiliser des réseaux est bien sûr de travailler avec $\Sk$ 
à la place de $\E^\ent$. La proposition B.1.8.3 
de \cite{fontaine} permet de faire cela. Elle implique par exemple
que :
\begin{equation}
\label{eq:fontaine}
\calT(M) = \Hom_{\Sk, \varphi} (\frakM, W_n(R) \cap {\E^{\ent,\ur}} 
/ p^n {\E^{\ent,\ur}}) = \Hom_{\Sk, \varphi} (\frakM, W_n(R)).
\end{equation}
Cette formule implique en particulier que tout morphisme de $\frakM$
dans $W_n(R)$ qui est compatible à $\varphi$ prend nécessairement ses
valeurs dans ${\E^{\ent,\ur}} / p^n {\E^{\ent,\ur}}$ (ce que l'on peut
démontrer directement sans difficulté).

\medskip

En s'inspirant de \cite{carliu}, on peut enfin donner une troisième
description de $\calT(M)$ qui fait intervenir non pas $W_n(R)$, mais
plutôt certains de ces quotients. Cette idée, qui pourrait sembler
inutilement complexe à première vue, va en fait s'avérer très fructueuse
tout au long de ce chapitre (comme elle l'a déjà d'ailleurs été dans
\cite{carliu}) : elle sera la clé pour obtenir des bornes sur la
ramification dans le \S \ref{subsec:bornes}, mais aussi dans le \S
\ref{subsec:formetau} lorsque l'on s'évertuera à établir des congruences
afin de préciser la forme de l'opérateur $\tau$. \emph{On se donne à
partir de maintenant un élément $U \in W(R)$ qui n'est pas multiple de
$p$ et on suppose que $\frakM$ est de hauteur divisant $U$.}

\begin{lemme}
\label{lem:existV}
Il existe un élément $V \in W(R)$ qui n'est pas multiple de $p$ et qui
vérifie $\varphi(V) = UV$.
\end{lemme}

\begin{proof}
On construit $V$ par approximations successives : par récurrence, on
construit une suite $(V_n)_{n \geq 1}$ d'éléments de $W(R)$ telle que
$V_{n+1} \equiv V_n \pmod {p^n}$ et $\varphi(V_n) = U V_n$ pour tout
$n$. Pour construire $V_1$, il suffit d'extraire une racine $(p-1)$-ième
de $U$ dans $R$, ce qui est possible puisque $\Frac R$ est
algébriquement clos et que $R$ est intégralement clos. Si maintenant
$V_n$ est construit, on cherche $V_{n+1}$ sous la forme $V_n + p^n X$
avec $X \in W(R)$. La condition que doit vérifier $X$ s'écrit :
$$\varphi(X) - U X \equiv \frac{U V_n - \varphi(V_n)}{p^n} \pmod p.$$ 
Si $x$ et $a$ désignent respectivement la réduction de $X$ et $\frac{U
V_n - \varphi{V_n}}{p^n}$ modulo $p$, trouver $X$ revient à résoudre
l'équation $x^p - U x = a$ dans $R$. Or, par le même argument que
précédemment, cette équation a bien une solution dans $R$, et la
récurrence se termine ainsi. Enfin $V = \lim_{n \to \infty} V_n$ existe
et vérifie $\varphi(V) = UV$.
\end{proof}

\begin{rem}
\label{rem:existV}
Il est évident que si $V$ vérifie $\varphi(V) = UV$ et si $a \in
\Z_p$, alors $aV$ vérifie la même équation. Un examen de la preuve
précédente montre que ce sont les seuls. Autrement dit, l'ensemble
des solutions de l'équation $\varphi(V) = UV$ est un $\Z_p$-module
libre de rang $1$. On notera en particulier que l'idéal engendré
par $V$ est uniquement déterminé en fonction de $U$.
\end{rem}

On rappelle que l'on a noté $\m_R$ l'idéal maximal de $R$, c'est-à-dire 
l'idéal formé des éléments de $R$ de valuation strictement positive.
Pour tout élément $X \in W(R)$, on pose
$$\calT_X(\frakM) = \Hom_{\Sk, \varphi} (\frakM, W_n(R)/
(X \cdot W_n(\m_R))).$$
La réduction modulo $X \cdot W_n(\m_R)$ définit des morphismes
canoniques $\rho_X : \calT(M) \to \calT_X(\frakM)$ et $\rho_{Y,X} :
\calT_Y(\frakM) \to \calT_X(\frakM)$ pour tout $Y \in W(R)$ multiple
de $X$. En particulier, on a le diagramme suivant :
$$\xymatrix {
& \calT(M) \ar[ld]_-{\rho_{UV}} \ar[rd]^-{\rho_V} \\
\calT_{UV}(\frakM) \ar[rr]^-{\rho_{UV,V}} & &\calT_V (\frakM) 
}$$
qui est manifestement commutatif.

\begin{prop}
\label{prop:rhoUVV}
Le morphisme $\rho_V$ est injectif, et son image s'identifie dans
$\calT_V(M)$ à l'image de $\rho_{UV,V}$ 
\end{prop}

\begin{proof}
On démontre d'abord l'injectivité. Soient $f$ et $g$ deux éléments de 
$\calT(M)$ tels que $f \equiv g \pmod{V \cdot W_n(\m_R)}$. Étant donné 
que $\frakM$ est de type fini, il existe un nombre réel $\nu > 0$ tel 
que la congruence $f \equiv g$ ait lieu modulo $V \cdot W_n(\a_R^{\geq 
\nu}$ où $\a_R^{\geq \nu}$ désigne l'idéal des éléments de $R$ de 
valuation $\geq \nu$. On pose $I = W_n(\a_R^{\geq \nu}$ ; on a alors 
$\bigcap_{i \geq 0} \varphi^i(I) = 0$ et $\varphi(I) \subset I$. Pour 
tout entier $i$, on pose $I_i = V\varphi^i(I)$ et on note $f_i$ (resp. 
$g_i$) la réduction de $f$ (resp. de $g$) modulo $I_i$. On va montrer, 
par récurrence, que $f_i = g_i$ pour tout entier $n$. Puisque 
$\bigcap_{i \geq 0} I_i = 0$, il en résultera que $f = g$, et donc 
l'injectivité souhaitée.
L'égalité $f_0 = g_0$ ayant déjà été justifiée, il suffit de traiter
l'hérédité. Soit $\psi : \frakM \to \Sk \otimes_{\varphi,\Sk} \frakM$,
$x \mapsto (\id \otimes \varphi)^{-1}(U x)$. Le morphisme $f$ induit
alors une application linéaire :
$$\id \otimes f_i : \Sk \otimes_{\varphi,\Sk} \frakM \to \Sk
\otimes_{\varphi,\Sk} \frac{W_n(R)}{I_iW_n(R)} = \frac{\Sk \otimes_{\Sk,
\varphi} W_n(R)}{\varphi(I_i)W_n(R)} = \frac{\Sk \otimes_{\varphi,\Sk}
W_n(R)}{U I_{i+1} W_n(R)}.$$
Le fait que $f$ commute à $\varphi$ implique la commutativité du
diagramme suivant :
\begin{equation}
\label{eq:diagfn}
\raisebox{0.5\depth}{\xymatrix @C=40pt {
\Sk \otimes_{\varphi,\Sk} \frakM \ar[d]^-{\id \otimes f_i}
& \frakM \ar[l]_-{\psi} \ar[r]^-{x \mapsto Ux}
& \frakM \ar[d]^-{f \mod U I_{i+1}} \\
\frac{\Sk \otimes_{\varphi,\Sk} W_n(R)}{U I_{i+1} W_n(R)}
\ar[rr]^-{\id \otimes \varphi} & &
\frac{W_n(R)}{U I_{i+1} W_n(R)} }}
\end{equation}
Comme on a bien sûr un diagramme analogue pour $g$, on déduit que
l'égalité $f_i = g_i$ implique 
$$U \cdot (f \mod U I_{i+1}) = U \cdot (g \mod U I_{i+1})$$
et donc finalement $f_{i+1} = g_{i+1}$ en divisant par $U$ (qui
n'est pas diviseur de zéro dans $W_n(R)$).

On en vient maintenant à la preuve de la propriété concernant les 
images. Soit $f_{UV}$ un élément de $\calT_{UV}(\frakM)$. Il s'agit de 
démontrer que $\rho_{UV,V}(f_{UV})$, qui est une application de $\frakM$ 
dans $W_n(R)/V \cdot W_n(\m_R)$ compatible au Frobenius, se relève en un 
morphisme $\frakM \to W_n(R)$ encore compatible au Frobenius. Étant 
donné que $\frakM$ est de type fini, il existe un idéal $I$ comme 
précédemment tel que $f_{UV}$ se relève en un morphisme $f_{UVI} : 
\frakM \to W(R)/UVI$. On pose $f_0 = f_{UVI} \mod I_0$. On va construire 
par récurrence sur $i$, une suite d'applications $f_i : \frakM \to 
W_n(R)/I_iW_n(R)$ telles que $f_{i+1} \equiv f_i \pmod {I_i}$ pour tout 
$i$. Le diagramme \eqref{eq:diagfn} assure que, si l'on pose
$$\alpha_i = (\id \otimes \varphi) \circ (\id \otimes f_i) \circ
\psi : \frakM \to W(R)/UI_{i+1}$$
un bon candidat pour $f_{i+1}$ est $\frac{\alpha_i} U$. Mais pour
pouvoir le définir ainsi, il faut montrer au préalable que $U$
divise $\alpha_i$. Or, ce dernier fait est vrai et suit des deux
remarques suivantes : \emph{primo}, commee $f_0$ se rélève 
modulo $UVI$ (le relève étant donné par $f_{UVI}$), l'élément $U$ (et
même en fait $UV$) divise $\alpha_0$, et \emph{secundo}, comme
$f_i \equiv f_0 \pmod V$, on a $\alpha_i \equiv \alpha_0 \pmod
{UV}$. On peut donc bien considérer l'application $\frac{\alpha_i} U$
qui est défini sur $\frakM$ mais, à cause de la division par $U$, prend
ses valeurs dans $W(R)/I_{i+1}$ (et pas $W(R)/U I_{i+1}$ comme c'était
le cas pour $\alpha_i$). Enfin, de la congruence $f_{i-1} \equiv f_i
\pmod{I_{i-1}}$, on déduit $\alpha_{i-1} \equiv \alpha_i
\pmod{\varphi(I_{i-1})}$, puis $f_i \equiv f_{i+1} \pmod{I_i}$ étant 
donné que $\varphi(I_{i-1}) = UI_i$. 
Enfin, en passant à la limite, on obtient une application $f : \frakM
\to W(R)$ qui relève $f_0$ et qui commute à $\varphi$.
\end{proof}

Il suit de la proposition précédente que $\calT(M)$ s'identifie à
l'image de $\rho_{UV,V}$ et donc, comme nous l'avions annoncé, nous
avons bien obtenu une description de cette représentation galoisienne
qui ne fait pas intervenir $W_n(R)$ lui-même mais seulement deux de ces
quotients.

\begin{rem}
La proposition \ref{prop:rhoUVV} vaut encore avec $n = \infty$ ; pour
l'établit, il suffit de passer à la limite sur $n$.
\end{rem}

\subsection{Bornes pour la ramification}
\label{subsec:bornes}

Sans surprise, le but de ce numéro est de démontrer le théorème
\ref{introtheo:bornes} de l'introduction. On commence par quelques brefs
rappels au sujet des filtrations de ramification dans le \S
\ref{subsec:ramif}. Les deux paragraphes suivants sont consacrés à la
preuve du théorème, tandis que dans le \S \ref{subsec:reciproque}, on
établit une réciproque partielle.

\subsubsection{Rappels sur les filtrations de ramification}
\label{subsec:ramif}

Soit $\kappa$ un corps complet pour une valuation discrète dont le corps
résiduel est parfait de caractéristique $p$ (dans les applications qui
nous intéressent, on prendra $\kappa = K$ ou $\kappa = F_0 = k((u))$). 
Pour toute extension finie $\kappa'$ de $\kappa$, on appelle $v_{\kappa'}$
la valuation sur $\kappa'$ normalisée par $v_{\kappa'}(\kappa'{}^\star) = \Z$.
Si $\kappa'$ est une extension galoisienne finie de $\kappa$ de groupe
de Galois $G$, la \emph{filtration de ramification en numérotation
inférieure} de $G$ est la filtration $(G_{(\lambda)})_{\lambda \in 
\R^+}$ définie comme suit :
$$G_{(\lambda)} = \big\{ \, \sigma \in G \,\, | \,\, v_{\kappa'} 
(\sigma(x) - x) \geq \lambda ,\, \forall x \in \O_{\kappa'} \, \big\}$$
où $\O_{\kappa'}$ est l'anneau des entiers de $\kappa'$. Les
$G_{(\lambda)}$ sont des sous-groupes distingués de $G$, et la
filtration qu'ils forment est décroissante, exhaustive et séparée. On
introduit encore la fonction $\varphi_{\kappa'/\kappa} : \R^+ \to
\R^+$ définie par :
$$\varphi_{\kappa'/\kappa}(\lambda) = \int_0^\lambda \frac{\Card
G_{(t)}}{\Card G_{(1)}} dt.$$
C'est une fonction affine par morceaux, concave et bijective, dont on
note $\psi_{\kappa'/\kappa}$ l'inverse. La \emph{filtration de 
ramification en numérotation supérieure} est définie par l'égalité
$G^{(\mu)} = G_{(\psi_{\kappa'/\kappa}(\mu))}$ pour tout réel 
$\mu \geq 0$. On renvoie à \cite{serre}, chap. IV pour les propriétés 
usuelles la concernant, et en particulier le théorème d'Herbrand. La filtration de
ramification en numérotation supérieure s'étend à une extension
galoisienne $\kappa'/\kappa$ non nécessairement finie en posant :
$$\Gal(\kappa'/\kappa)_{(\mu)} = \varprojlim_{\kappa''} 
\Gal(\kappa''/\kappa)_{(\mu)}$$
où la limite projective est prise sur toutes les extensions finies
galoisiennes $\kappa''$ de $\kappa$ incluses dans $\kappa'$. Dans le 
cas où $\kappa'/\kappa$ est une extension algébrique séparable non
galoisienne, on ne peut certes pas définir de filtration sur le groupe
de Galois puisque celui-ci n'existe pas mais les fonctions
$\varphi_{\kappa'/\kappa}$ et $\psi_{\kappa'/\kappa}$, elles, ont
encore un sens ; on peut les définir, par exemple, grâce aux
formules :
$$\psi_{\kappa'/\kappa} (\mu) = \int_0^\mu [I_\kappa : I_{\kappa'}
G_\kappa^{(t)}] dt \quad \text{et} \quad \varphi_{\kappa'/\kappa} =
\psi_{\kappa'/\kappa}^{-1}$$
où $I_\kappa$ et $I_{\kappa'}$ sont respectivement les sous-groupes
d'inertie des groupes de Galois absolus de $\kappa$ et $\kappa'$ 
(voir \cite{wintenberger}, \S 1.2.1).
La théorie qui vient d'être rappelée brièvement s'applique en
particulier aux corps $\kappa = K$ et $\kappa = F_0$. Les groupes de
Galois absolus $G_K$ et $\Gal(F_0^\sep/F_0) \simeq G_\infty$ héritent 
ainsi d'une filtration de ramification en numération supérieure.

\subsubsection{Bornes pour les représentations de $G_\infty$}

On démontre dans ce paragraphe l'assertion 1 du théorème
\ref{introtheo:bornes}. Soient $M$ un $\varphi$-module $M$ sur $\E^\ent$
annulé par $p^n$, et $\frakM$ un $\varphi$-réseau de $M$ de hauteur
divisant $U$, pour un élément $U \in \Sk$ qui n'est pas multiple de $p$.
Comme dans l'énoncé du théorème, on note $T$ la représentation
galoisienne associée à ces données et on pose $h = v_R(U \mod p)$. Pour
tout nombre réel $v \geq 0$ et tout anneau $A$ muni d'une valuation (par
exemple $A \subset R$), on note encore $\a_A^{>v}$ l'idéal des éléments
de $A$ de valuation strictement plus grande que $v$. Enfin, si $F$ est
une extension algébrique de $F_0$, on désigne par $\O_F$ son anneau des
entiers. Dans le cas où $F$ est inclus dans $F_0^\sep$, on a simplement
$\O_F = F \cap R$.

Comme cela se fait usuellement dans ce genre de situations, on va
utiliser la propriété $(P_m)$ de Fontaine (introduite dans
\cite{fontaine-varab}, proposition 1.5). Rappelons qu'ici $m$ est un
nombre réel positif ou nul et que, par définition, une extension $F$ de
$F_0$ vérifie $(P_m)$ si, et seulement si pour toute extension
algébrique $E$ de $F_0$, s'il existe un morphisme de $\O_{F_0}$-algèbres
$\O_F \to \O_E/\a_E^{>m}$, alors il existe un $F_0$-plongement de $F$
dans $E$. Le lien avec la filtration de ramification en numérotation
supérieure est donnée par la proposition suivante qui, dans cette
formulation, est dûe à Yoshida (voir \cite{yoshida}).

\begin{prop}
Soit $F$ une extension finie galoisienne de $F_0$ de groupe de
Galois $G$. On définit :
\begin{itemize}
\item l'entier $m_0$ comme la borne inférieure des réels $m$ 
tels que $(P_m)$ soit satisfaite, et
\item l'entier $\mu_0$ comme la borne inférieure des réels $\mu$
tel que $G^{(\mu)} = \{\id_F\}$.
\end{itemize}
Alors $m_0 = \mu_0$.
\end{prop}

Notons $\rho : G_\infty \to \GL(T)$ le morphisme donnant l'action de
$G_\infty$ sur $T$, et $F$ l'extension finie galoisienne de $F_0$ qui
est en correspondance avec le sous-groupe distingué $\ker \rho$ de
$G_\infty$. D'après la proposition précédente, pour démontrer la
première assertion du théorème \ref{introtheo:bornes}, il suffit de
prouver que l'extension $F$ vérifie la propriété $(P_m)$ pour $m =
\frac{h p^n}{p-1}$.

\begin{lemme}
\label{lem:Fcut}
Soit $E$ une extension de $F_0$ incluse dans $F_0^\sep$. Alors
l'application injective naturelle 
$$\Hom_{\Sk,\varphi} (\frakM, W_n(\O_E)) \to 
\Hom_{\Sk,\varphi} (\frakM, W_n(R)) = T$$
(déduite de l'inclusion $W_n(\O_E) \to W_n(R)$) est un isomorphisme
si, et seulement si $E$ contient $F$.
\end{lemme}

\begin{proof}
Il est clair que si l'application du lemme est un isomorphisme, alors le
groupe de Galois absolu de $E$ agit trivialement sur $T$. D'où $F
\subset E$. Pour la réciproque, on remarque qu'en vertu de
l'égalité \eqref{eq:fontaine}, on sait que tous les morphismes $f :
\frakM \to W_n(R)$ compatibles à $\varphi$ prennent leurs valeurs dans
${\E^{\ent,\ur}}$ et donc, en particulier, dans $W_n(\O_{F_0^\sep})$. Par
ailleurs, par définition de $F$, son groupe de Galois absolu
$G_F$ agit trivialement sur $T$, ce qui implique que tous les
morphismes $f$ comme précédemment prennent leurs valeurs dans $W_n
(\O_{F_0^\sep})^{G_F} = W_n(\O_F)$. Il en résulte que, 
si $F \subset E$, l'application du lemme est bien bijective.
\end{proof}

Pour tout réel $v \geq 0$, on introduit à présent le sous-ensemble
$W_n(\a_R^{>v})$ de $W_n(R)$ formé des éléments dont toutes les
coordonnées sont dans $\a_R^{>v}$. C'est un idéal de $W_n(R)$, et le
quotient de $W_n(R)/ W_n(\a_R^{>v})$ s'identifie à $W_n(R/\a_R^{>v})$.
D'après le lemme \ref{lem:existV}, il existe $V \in W(R)$, non multiple
de $p$, tel que $\varphi (V) = UV$. on fixe un tel élément $V$ ; il
vérifie $V^{p-1} = U \pmod p$, ce qui implique que $v_R(V \mod p) =
\frac h {p-1}$.

\begin{lemme}
\label{lem:incWitt}
On a $W_n(\a_R^{>m}) \subset UV \cdot W_n(\m_R)$ pour $m =
\frac{hp^n}{p-1}$.
\end{lemme}

\begin{proof}
On raisonne par récurrence sur $n$ et donc, pour éviter les confusions,
on notera $m(n)$ à la place de $m$ tout au long de la démonstration.
Pour $n = 1$, on remarque que $v_R(UV \mod p) = m(1)$ et donc que les
deux idéaux considérés sont égaux. On suppose à présent que l'inclusion
$W_n(\a_R^{>m(n)}) \subset UV \cdot W_n(\m_R)$ est satisfaite et
on considère un élément $X = (x_1, \ldots, x_{n+1}) \in
W_{n+1}(\a_R^{>m(n+1)})$. On veut montrer que $X \in UV \cdot
W_{n+1}(\m_R)$. Soient $\lambda \in R$ un élément quelconque de
valuation $hp$ et $[\lambda] \in W(R)$ son représentant de Teichmüller.
Un calcul immédiat sur les valuations montre que les composantes du
vecteur de Witt de longueur $n$ suivant :
$$\frac 1 {[\lambda]} (x_1, \ldots, x_n) = 
\Big(\frac{x_1} \lambda, \frac{x_2}{\lambda^p}, \ldots, 
\frac{x_n}{\lambda^{p^{n-1}}}\Big)$$
sont toutes dans $\a_R^{>m(n)}$. L'hypothèse de récurrence 
s'applique donc et assure qu'il existe $Y \in W_n(\m_R)$ tel que
$\frac 1 {[\lambda]} (x_1, \ldots, x_n) = UVY$. On note encore $Y$ un
élément de $W_{n+1}(\m_R)$ qui relève $Y$ et on pose
$$\Delta = (x_1, \ldots, x_n, x_{n+1}) - [\lambda] UVY.$$
Les $n$ premières coordonnées de $\Delta$ sont nulles, tandis que, si on
pose $[\lambda] UVY = (z_1, \ldots, z_{n+1})$, la dernière coordonnée de
$\Delta$ s'exprime comme un polynôme homogène de degré $p^n$ en les
$x_1, \ldots, x_{n+1}, z_1, \ldots, z_{n+1}$, à condition de donner le
poids $p^{i-1}$ aux variables $x_i$ et $z_i$. Comme, en outre, $z_i$ 
s'écrit comme le produit de $\lambda^{p^i}$ par un élément de $\m_R$, 
on en déduit que la dernière coordonnée de $\Delta$ est de valuation 
strictement supérieure à $m(n+1)$. Ainsi, $\Delta$ appartient à $p^n UV 
\cdot W_{n+1}(\m_R)$ et donc \emph{a fortiori} de $UV \cdot 
W_{n+1}(\m_R)$. Il s'ensuit enfin que $X \in UV \cdot W_{n+1}(\m_R)$ 
comme voulu.
\end{proof}

Nous sommes prêts à vérifier la propriété $(P_m)$ pour le corps $F$ et
le nombre $m = \max(1,\frac{hp^n}{p-1})$. Soit $E$ une extension
algébrique de $F_0$ incluse dans $F_0^\sep$. Soit également $f : \O_F
\to \O_E/\a_E^{> m}$ un morphisme de $\O_{F_0}$-algèbres. On considère
la composée suivante :
$$\begin{array}{ll}\Psi_V : &
T = \Hom(\frakM, W_n(\O_F)) \to \Hom(\frakM, W_n(\O_F/\a_F^{>m}))
\to \Hom(\frakM, W_n(\O_E/\a_E^{>m})) \medskip \\
& \displaystyle \hspace{1cm}
\to \Hom\Big(\frakM, \frac{W_n(\O_E)}{UVW_n(\m_R) 
\cap W_n(\O_E)}\Big)
\to \Hom\Big(\frakM, \frac{W_n(\O_E)}{VW_n(\m_R) 
\cap W_n(\O_E)}\Big)
\end{array}$$
où la deuxième flèche est induite par $f$ et l'existence de la troisième
résulte du lemme \ref{lem:incWitt}. Une adaptation de la proposition
\ref{prop:rhoUVV} assure que, pour tout $\psi_F \in T$, il existe une
unique application $\psi_E \in \Hom(\frakM, W_n(\O_E))$ telle que
$\Psi_V(\psi_F)$ s'identifie à $\psi_E \mod V W_n(\m_R) \cap W_n(\O_E)$.
Ceci permet de construire un morphisme $\Psi : T \to \Hom(\frakM,
W_n(\O_E))$ relevant $\Psi_V$.

\begin{lemme}
\label{lem:Psiinj}
Le morphisme $\Psi$ précédent est injectif.
\end{lemme}

\begin{proof}
On raisonne par récurrence sur $n$. Pour $n = 1$, on considère
$u_F$ une uniformisante de $F$. Son polynôme minimal sur $F_0$ est un
polynôme d'Eisenstein que l'on note $P$. L'élement $x = f(u_F) \in
\O_E/\a_E^{>m}$ est alors une racine de $P$. Comme $m \geq 1$, le
coefficient constant de $P$ ne s'annule pas dans $\O_E/\a_E^{>m}$. On en
déduit que $x$ a la même valuation que $u_F$ puis que $f$ induit une
application injective $f_{h/(p-1)} : \O_F/\a_E^{> h/(p-1)} \to
\O_E/\a_E^{> h/(p-1)}$. Ainsi le morphisme
$$\textstyle \Hom\big(\frakM, \frac{\O_E}{V \m_R \cap \O_E}\big) = 
\Hom(\frakM, \O_E/\a_E^{>h/(p-1)}) \to \Hom(\frakM, \O_F/\a_F^{>h/(p-1)}) 
= \Hom\big(\frakM, \frac{\O_F}{V \m_R \cap \O_F}\big)$$
est, lui aussi, injectif, ce qui permet de conclure.

Pour l'hérédité, on suppose que $\frakM$ est annulé par $p^{n+1}$ et on
considère la suite exacte $0 \to p\frakM \to \frakM \to \frakM / p\frakM
\to 0$. Les modules $p\frakM$ et $\frakM/p\frakM$ sont encore des
$\varphi$-réseaux à l'intérieur, respectivement de $pM$ et $M/pM$. En
outre, ils sont tous les deux de hauteur divisible par $U$. En effet,
c'est évident pour $\frakM/p\frakM$ et, pour $p\frakM$, on raisonne
comme suit. On note tout d'abord que l'application $\id \otimes \varphi :
W(R) \otimes_{\varphi,\Sk} \frakM/p\frakM \to W(R) \otimes_\Sk
\frakM/p\frakM$ est injective, ce qui implique que le morphisme suivant induit par
l'inclusion naturelle de $p\frakM$ dans $\frakM$ :
$$\coker\big( \id \otimes \varphi : W(R) \otimes_{\varphi,\Sk} (p\frakM)
\to W(R) \otimes_\Sk (p\frakM) \big) \longrightarrow \coker \big( W(R)
\otimes_{\varphi,\Sk} \frakM \to W(R) \otimes_\Sk \frakM \big)$$
est, lui aussi, injectif. Comme l'espace d'arrivée est annulé par $U$,
on en déduit qu'il en est de même de l'espace de départ, ce qui veut
bien dire que $p\frakM$ est de hauteur divisant $U$. On considère
maintenant le diagramme suivant :
$$\xymatrix @R=20pt {
0 \ar[r] &
\calT(\frakM/p\frakM) \ar[r] \ar[d] & \calT(\frakM) \ar[r] \ar[d] & 
\calT(p\frakM) \ar[r] \ar[d] & 0 \\
& \Hom(\frakM/p\frakM, \O_E) \ar@{^(->}[d]
& \Hom(\frakM, W_{n+1}(\O_E)) \ar@{^(->}[d]
& \Hom(p\frakM, W_n(\O_E)) \ar@{^(->}[d] \\
0 \ar[r] &
\calT(\frakM/p\frakM) \ar[r] & \calT(\frakM) \ar[r] & \calT(p\frakM)
\ar[r] & 0
}$$
où les flèches verticales du haut sont les morphismes $\Psi$
correspondant respectivement à $\frakM/p\frakM$, $\frakM$ et $p\frakM$,
et les flèches verticales du bas sont les inclusions canoniques. Par
hypothèse de récurrence, on sait que les flèches verticales en haut à
gauche et en haut à droite sont injectives. Un chasse au diagramme
termine alors la récurrence.
\end{proof}

Il est maintenant aisé de conclure. Du lemme \ref{lem:Psiinj}, il
découle que l'ensemble $\Hom(\frakM, W_n(\O_E))$ a au moins autant
d'éléments que $T$. L'inclusion naturelle $\Hom(\frakM, W_n(\O_E)) \to
T$ est donc nécessairement une bijection. Le lemme \ref{lem:Fcut} assure
que $E$ contient $F$, ce qui est exactement ce qu'il fallait vérifier
pour établir la propriété $(P_m)$.

\subsubsection{Bornes pour les représentations de $G_K$}
\label{subsec:borneGK}

Expliquons à présent comment la deuxième partie du théorème
\ref{introtheo:bornes} se déduit de la première. On rappelle tout
d'abord que la théorie du corps des normes ne se contente pas de fournir
un isomorphisme canonique entre les groupes $G_\infty$ et $G_{F_0} =
\Gal(F_0^\sep/F_0)$, mais qu'elle compare aussi les filtrations de
ramification en numérotation supérieure. Précisément, le corollaire
3.3.6 de \cite{wintenberger} affirme que pour tout réel $\mu \geq 0$,
les groupes $G_{F_0}^{(\mu)}$ et $G_\infty \cap
G_K^{(\varphi_{K_\infty/K}(\mu))}$ se correspondent \emph{via}
l'isomorphisme du corps des normes. De surcroît, la fonction
$\varphi_{K_\infty/K}$ a été calculée dans \cite{carliu}, \S 4.3. Voici
sa représentation graphique :

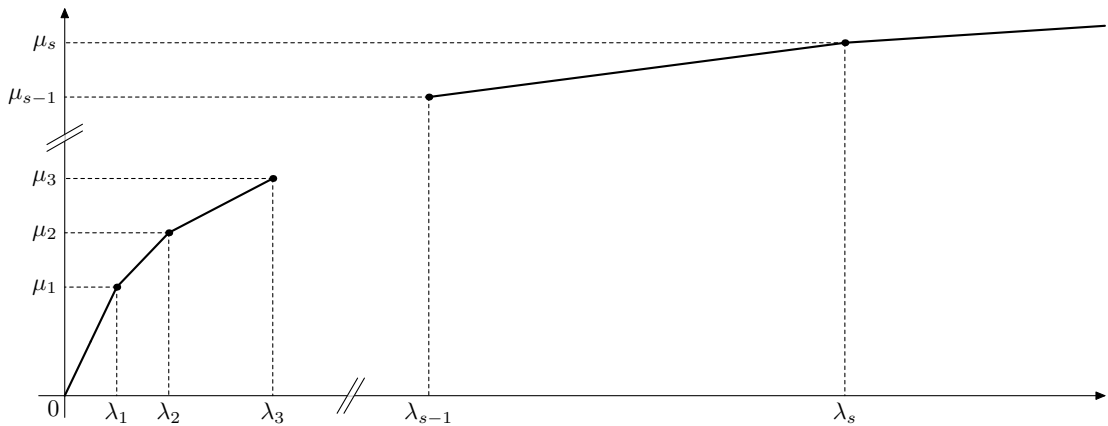
\begin{center}
\includegraphics[scale=0.9]{figure.1}
\end{center}

\noindent
où $\lambda_s = 1 + \frac{ep^s}{p-1}$ et $\mu_s = 1 + e(s + \frac
1{p-1})$. Si l'on souhaite une formule explicite, on peut écrire
pour tout $\lambda \geq \lambda_1$ :
\begin{equation}
\label{eq:phiKinfK}
\varphi_{K_\infty/K}(\lambda) = 1 + es + e \cdot \Bigg(
\frac{p^{\{s\}}}{p-1} - \{s\}\Bigg)
\end{equation}
où $s = \log_p(\frac{(p-1)(\lambda-1)} e) \geq 1$ et où $\{s\}$ désigne
sa partie décimale. On voit en particulier que pour sur l'intervalle
$[\lambda_1, +\infty[$, la fonction $\varphi_{K_\infty/K}$ s'écrit comme
la somme de la fonction $\lambda \mapsto e \cdot \log_p \lambda$ et d'une 
fonction bornée.

\medskip

Rappelons qu'au début de l'introduction, nous avions défini l'extension
$K(\zeta_{p^\infty})/K$ obtenue en ajoutant toutes les racines
primitives $p^n$-ième de l'unité et que nous avions posé $\Gamma =
\Gal(K(\zeta_{p^\infty})/K)$. Le caractère cyclotomique $\chi$ identifie
$\Gamma$ à un sous-groupe d'indice fini de $\Z_p^\times$. Ce groupe est,
par ailleurs, muni de la filtration de ramification en numérotation
supérieure, notée $(\Gamma^{(s)})_{s \geq 0}$ dans la suite.

\begin{lemme}
\label{lem:ramifGamma}
Il existe une constante $c_0(K) \geq 0$ ne dépendant que de $K$ telle
que pour tout $\gamma \in \Gamma$, on ait $\gamma \in
\Gamma^{(1+es-c_0(K))}$ où $s$ est la valuation $p$-adique de
$\chi(\gamma)-1$.
\end{lemme}

\begin{proof}
Le lemme est une conséquence directe d'un théorème très général de Sen
(voir \cite{sen}). Nous donnons toutefois ci-dessous une autre
démonstration, plus élémentaire. Si le corps $K$ n'est pas absolument
ramifiée, le lemme résulte d'un calcul classique que l'on peut trouver
par exemple dans \cite{serre}, chap. IV, \S 4 ; on peut alors même choisir
$c_0(K) = 0$ (attention au décalage d'indice dans la numérotation des
groupes de ramification entre notre convention et celle de \emph{loc.
cit.}).

Pour le cas général, on pose $K' = W[1/p]$ (c'est la sous-extension
maximale non absolument ramifiée de $K$) et $\Gamma' = \Gal(K'(\zeta
_{p^\infty})/K') \simeq \Z_p^\times$. Le groupe $\Gamma$ apparaît
naturellement comme un sous-groupe de $\Gamma'$ et d'après le
théorème d'Herbrand, on a l'égalité :
$$\Gamma^{(\mu)} = \Gamma \cap (\Gamma')^{(\varphi_{K'/K}(\mu))}.$$
Ainsi, en utilisant le cas déjà traité d'un corps de base non absolument
ramifié, on obtient $\gamma \in \Gamma^{(\psi_{K'/K}(1+s))}$. Or à 
partir d'un certain moment, la fonction $\psi_{K'/K}$ est affine de
pente $e$. On en déduit qu'il existe une constante $c_0(K) \geq 0$
telle que $\psi_{K'/K}(1+s) \geq 1 + es - c_0(K)$. Le lemme s'ensuit.
\end{proof}

\begin{rem}
\label{rem:Kmodram}
Dans le cas où l'extension $K/K'$ est modérément ramifiée ---
c'est-à-dire dans le cas où $e$ est premier avec $p$ ---, la fonction
$\psi_{K'/K}$ vaut l'identité sur $[0,1]$, et est tout de suite après
affine de pente $e$. On en déduit que, dans ce cas, on peut prendre
$c_0(K) = 0$, c'est-à-dire que $\gamma \in \Gamma^{(1+es)}$.
\end{rem}

Étant donné que l'image de $\chi : \Gamma \to \Z_p^\times$ est d'indice
finie dans $\Z_p^\times$, il existe un entier $s_0(K)$ --- ne dépendant
bien sûr que de $K$ --- tel que $1 + p^{s_0(K)} \Z_p \subset \chi(\Gamma)$.
Soit $s$ un nombre entier vérifiant :
\begin{equation}
\label{eq:inegs}
s \geq s_0(K) \quad \text{et} \quad s > \frac{c_0(K)} e +
\textstyle \max\big(\frac 1{p-1}, n + \log_p(\frac h e)\big).
\end{equation}
On choisit en outre un élément $\tau \in G_K$ tel que $c(\tau) = 1$ et 
$\chi(\tau) = 1$. Soient $\mu$ un nombre réel $> 1 + e(s + \frac 
1{p-1})$ et $g \in G_K^{(\mu)}$. D'après ce que l'on sait à propos de la 
ramification de l'extension $K_\infty/K$, on a $g \in G_s$. D'après le 
lemme \ref{lem:strucGK} et la discussion menée au \S 
\ref{subsec:quotHinf}, quitte à modifier un peu $\tau$, on peut écrire 
$g$ sous la forme $g = \tau^a g'$ pour un certain $g' \in G_\infty$. Du 
fait que $g$ fixe $K_s$ point par point, on déduit qu'il en est de même 
pour $\tau^a$, et donc que $p^s$ divise $a$. Par ailleurs, puisque $s$ 
est choisi supérieur ou égal à $s_0(K)$, il existe un élément $\gamma 
\in \Gamma$ tel que $\chi(\gamma) = 1 + a$. L'autre inégalité 
supposée sur $s$ (voir formule \eqref{eq:inegs}),
quant à elle, implique, avec le lemme \ref{lem:ramifGamma}, que 
$\gamma \in \Gamma^{(\mu_0)}$ avec $\mu_0 = 1 + e \cdot \max(\frac 
1{p-1}, n + \log_p(\frac h e))$.
Par ailleurs, étant donné que les extensions $K_\infty/K$ et
$K(\zeta_{p^\infty})/K$ sont linéairement disjointes, le groupe de
Galois $\Gal(K_\infty (\zeta_{p^\infty})/K_\infty)$ s'identifie
canoniquement à $\Gamma$. \emph{Via} cette identification, on a en 
outre l'égalité :
$$\Gal(K_\infty (\zeta_{p^\infty})/K_\infty) \cap G_K^{(\mu_0)} =
\Gamma^{(\mu_0)}.$$
On en déduit que $\gamma$ (qui, pour l'instant, vit dans $\Gamma$)
se relève en un
élement de $G_K$ qui appartient à $G_\infty \cap G_K^{(\mu_0)}$. À
partir de maintenant, on fixe $\gamma$ un tel relevé ; il vérifie
encore $\chi(\gamma) = 1+a$. La relation de commutation 
$\gamma \tau \equiv \tau^{\chi(\gamma)} \gamma \pmod {H_\infty}$ 
entraîne donc l'existence d'un élément $g'' \in H_\infty$ tel que $\tau^a = 
\gamma \tau \gamma^{-1} \tau^{-1} g''$. On a alors l'écriture $g'' g' 
= (\tau \gamma \tau^{-1}) \cdot \gamma^{-1} \cdot g$, de laquelle on
déduit que $g'' g' \in G_K^{(\mu_0)}$ (puisque $\mu > 1 + e(s+\frac 1
{p-1}) > \mu_0$). Ainsi $g'' g' \in G_\infty \cap G_K^{(\mu_0)} = 
G_\infty^{\psi_{K_\infty/K}(\mu_0)}$. Or, un calcul facile basé sur
la formule \eqref{eq:phiKinfK} montre que $\psi_{K_\infty/K}(\mu_0)
> \max(1, \frac{h p^n}{p-1})$. 
L'assertion 1 du théorème \ref{introtheo:bornes} implique alors que $g'' 
g'$ agit trivialement sur $T$. Comme, de même, $\gamma$ était élément de 
$G_\infty \cap G_K^{(\mu_0)}$, il agit aussi trivialement sur $T$. Ainsi, 
de l'écriture $g = \gamma \tau \gamma^{-1} \tau^{-1} \cdot (g'' g')$, on 
déduit que $g$, à son tour, agit trivialement sur $T$.
Comme ceci est valable pour tout $g \in G_K^{(\mu)}$, tout $\mu > 1 + 
e(s + \frac 1 {p-1})$ et tout $s$ vérifiant les inégalités 
\eqref{eq:inegs}, on a finalement démontré qu'en posant\footnote{La 
constante peut être encore légèrement améliorée, mais la forme que nous 
donnons nous a semblé un peu plus agréable.} $$c(K) = 1 + \frac e{p-1} + 
e \cdot s_0(K) + c_0(K),$$ le sous-groupe $G_K^{(\mu)}$ agit 
trivialement sur $T$ pour tout $\mu > c(K) + e\cdot \max(\frac 1{p-1}, n 
+ \log_p (\frac h e))$. Le théorème \ref{introtheo:bornes} est démontré.

\medskip

Lorsque $K$ est absolument modérément ramifiée, on a déjà vu, dans la
remarque \ref{rem:Kmodram}, que l'on pouvait prendre $c_0(K) = 0$. Par
ailleurs, il est facile de voir que l'entier $s_0(K)$ peut être, dans ce
cas, choisi égal à $1$. La formule qui a été obtenue précédemment
donne alors la valeur $1 + e + \frac e{p-1}$ pour $c(K)$. En réalité, un
examen de la preuve précédente montre que, dès que $h \geq e$, on peut
même prendre $c(K) = 1 + \frac e{p-1}$.

\subsubsection{Une réciproque partielle au théorème
\ref{introtheo:bornes}}
\label{subsec:reciproque}

Le théorème \ref{introtheo:bornes} que nous venons de démontrer permet
de contrôler la ramification d'une représentation de torsion de
$G_\infty$ ou $G_K$ en fonction de la hauteur du $\varphi$-module
associé. Dans le cas des objets annulées par $p$, il se trouve qu'être
de hauteur divisant $U$ (pour un élément $U \in \Sk$) équivaut à être de
$u$-hauteur $\leq h = v_R (U \mod p)$. Ainsi, la ramification d'une
$\F_p$-représentation est contrôlée par la $u$-hauteur du
$\varphi$-module associé, qui est simplement un nombre entier : dans le
cas des représentations de $G_\infty$ par exemple, si cette $u$-hauteur
est $\leq h$, alors le sous-groupe $G_\infty^{(\mu)}$ agit trivialement
pour tout $\mu > \frac{hp}{p-1}$. 

On peut alors se demander si, réciproquement, la ramification de la
représentation contrôle la $u$-hauteur du $\varphi$-module associé. La
proposition suivante apporte une réponse affirmative à cette question.

\begin{prop}
\label{prop:reciproque}
Soit $T$ une $\F_p$-représentation de dimension finie de $G_\infty$,
et soit $h$ un entier tel que $G_\infty^{(\mu)}$ agisse trivialement
sur $T$ pour tout $\mu > \frac{hp}{p-1}$. Alors le $\varphi$-module
étale associé à $T$ admet un $\varphi$-réseau de $u$-hauteur $\leq
hp$.

\smallskip

Soit $T$ une $\F_p$-représentation de dimension finie de $G_K$, et soit
$h$ un entier strictement positif tel que $G_K^{(\mu)}$ agisse
trivialement sur $T$ pour tout $\mu > 1 + \frac e{p-1} + e \cdot (1 +
\log_p(\frac h e))$. Alors le $(\varphi,\tau)$-module associé à $T$
admet un $(\varphi,\tau)$-réseau de $u$-hauteur $\leq hp$.
\end{prop}

\begin{rem}
On prendra garde d'une part à ce que la borne portant sur $\mu$ ne fait 
pas intervenir la constante $c(K)$ comme c'était le cas dans le théorème
\ref{introtheo:bornes} et d'autre part à ce que la borne que l'on 
obtient au final pour la $u$-hauteur est bien $hp$, et non pas $h$.
\end{rem}

\begin{proof}
La deuxième partie de la proposition se déduit de la première en 
utilisant les liens rappelés précédemment entre les filtrations de 
ramification sur $G_\infty$ et $G_K$. Nous laissons cet exercice au 
lecteur et nous concentrons à partir de maintenant sur la preuve de la 
première affirmation.

Le $\varphi$-module étale associé à la représentation $T$ est $M =
\Hom_{\F_p[G_\infty]} (T, F_0^\sep)$ et admet un $\varphi$-réseau
canonique donné par $\frakM = \Hom_{\F_p[G_\infty]} (T, \O_{F_0^\sep})$.
Montrons que $\frakM$ est de $u$-hauteur $\leq hp$. Étant donné notre
hypothèse, il existe une extension galoisienne finie $F$ sur $F_0$ telle
que, d'une part, $\frakM = \Hom_{\F_p[G_\infty]} (T, \O_F)$ et, d'autre
part, le sous-groupe de ramification $\Gal(F/F_0) ^{(\mu)}$ soit trivial
pour tout $\mu > \frac{hp}{p-1}$. Le corps $F$ est à l'évidence un
$F_0$-espace vectoriel de dimension finie muni d'un opérateur de
Frobenius, qui est tout simplement l'élévation à la puissance $p$ :
c'est donc un $\varphi$-module. Son anneau des entiers $\O_F$ définit en
outre un $\varphi$-réseau à l'intérieur de $F$. Le lemme 3.2 de
\cite{leborgne} implique que ce $\varphi$-réseau est de $u$-hauteur
$\leq \lceil (p-1) v_R( \d_{F/F_0}) \rceil$ où $\d_{F/F_0}$ est la
différente de l'extension $F/F_0$ et où $\lceil x \rceil$ désigne la
partie entière supérieure du nombre réel $x$. De la proposition 1.3 de
\cite{fontaine-varab} et de la borne sur la ramification de $F$, il suit
$v_R(\d(F/F_0)) \leq \frac{hp} {p-1}$, d'où on déduit que le
$\varphi$-réseau $\O_F$ est de $u$-hauteur $\leq hp$.

Reste à montrer comment cela implique que $\frakM$ est, lui-même, de 
$u$-hauteur $\leq hp$.  Soit $f$ un élément de $\frakM$ ; c'est un 
morphisme de $T$ dans $\O_F$ compatible à l'action de Galois. On veut 
montrer que $u^{hp} f$ est dans l'image de $\id \otimes \varphi_\frakM : 
\Sk \otimes_{\varphi,\Sk} \frakM \to \frakM$. Or, l'application $u^{hp} 
f$ prend manifestement ses valeurs dans le sous-ensemble $u^{hp} \O_F$ 
qui est inclus dans l'image de $\id \otimes \varphi_{\O_F} : \Sk 
\otimes_{\varphi,\Sk} \O_F \to \O_F$ puisque $\O_F$ est de $u$-hauteur 
$\leq hp$. Par ailleurs, on vérifie facilement que cette dernière 
application est injective. Il existe par suite un morphisme $g : T \to 
\Sk \otimes_{\varphi,\Sk} \O_F$ qui est compatible à l'action de Galois, 
et permet de factoriser $u^{hp}f$ comme suit : $u^{hp} f = (\id \otimes 
\varphi_{\O_F}) \circ g$. Ce morphisme $g$ définit un élément de $\Sk 
\otimes_{\Sk, \varphi} \frakM$ qui s'envoie sur $u^{hp} f$ par $\id 
\otimes \varphi_\frakM$. Ceci montre bien que $u^{hp} f$ est dans 
l'image de $\id \otimes \varphi_\frakM$.
\end{proof}

Le résultat précédent s'étend aux représentations annulées par une
puissance de $p$ comme suit.

\begin{prop}
\label{prop:reciproque2}
Soit $T$ une $\Z_p$-représentation de longueur finie de $G_\infty$
(resp. de $G_K$), soit $n$ un entier tel que $p^n T = 0$, et soit $h$ 
un entier tel que $G_\infty^{(\mu)}$ (resp. $G_K^{(\mu)}$) agisse
trivialement sur $T$ pour tout $\mu > \frac{hp}{p-1}$ (resp pour tout
$\mu > 1 + \frac e{p-1} + e \cdot (1 + \log_p(\frac h e))$).  Alors,
pour tout élément $U \in W(R)$ tel que $v_R(U \mod p) \geq hp$, le
$\varphi$-module étale (resp. le $(\varphi,\tau)$-module) associé à $T$
admet un $\varphi$-réseau (resp. un $(\varphi,\tau)$-réseau) de
$U$-hauteur $\leq n$.
\end{prop}

\begin{proof}

De même que pour la proposition \ref{prop:reciproque}, il suffit de 
traiter le cas \og $G_\infty$ \fg. Si $T$ une représentation vérifiant 
les hypothèses de l'énoncé, il existe une extension finie galoisienne 
$F/F_0$ telle que $\Gal(F_0^\sep/F)$ agisse trivialement sur $T$ et 
$\Gal(F/F_0)^{(\mu)}$ soit trivial pour tout $\mu > \frac{hp}{p-1}$. 
Soit $\E^{\ent,F}$ l'unique extension étale de $\E^\ent$ incluse dans 
$\E^{\ent,\ur}$ dont le corps résiduel est égal à $F$. 
Le $\varphi$-module associé à $T$ s'écrit alors $M = 
\Hom_{\Z_p[G_\infty]}(T, \E^{\ent,F})$ et un $\varphi$-réseau naturel à 
l'intérieur de $M$ est $\frakM \Hom_{\Z_p[G_\infty]}(T, \Sk^F/p^n\Sk^F)$ 
avec $\Sk^F = \E^{\ent,F} \cap W(R)$. Pour conclure, il suffit de 
démontrer que $\Sk^F/p^n\Sk^F$, vu comme $\varphi$-module sur $\Sk$, est 
de $U$-hauteur $\leq n$. Or, par la démonstration de la proposition 
\ref{prop:reciproque}, on sait que cette propriété vaut pour $n=1$ 
(puisque $\Sk^F/p\Sk^F$ s'identifie à l'anneau des entiers de $F$) et, 
par ailleurs, on a la suite exacte
$$0 \to \Sk^F/p^{n-1} \Sk^F \to \Sk^F/p^n \Sk^F \to \Sk^F/p\Sk^F
\to 0$$
à partir de laquelle on conclut par récurrence.
\end{proof}

On remarquera que, si pour $n = 1$, l'écart entre les bornes appraissant 
dans le théorème \ref{introtheo:bornes} d'une part et la proposition 
\ref{prop:reciproque} ne différaient que par la constante $c(K)$ dans la 
borne sur $\mu$ et par un facteur $p$ dans celle sur la $u$-hauteur, ce 
même écart se creuse de façon drastique lorsque $n$ augmente : dans la 
proposition \ref{prop:reciproque2}, la présence de $n$ n'apparaît plus 
dans la borne sur $\mu$ mais dans celle sur la hauteur !

\subsection{Les $(\varphi,\tau)$-réseaux}
\label{subsec:phitaures}

Avant de pouvoir définir des réseaux dans les $(\varphi,\tau)$-modules
sur $(\E^\ent, {\E^\ent_\tau})$, il est nécessaire introduire des
sous-anneaux de $\E^\ent$ et ${\E^\ent_\tau}$ sur lesquels ces réseaux
seront définis. Le sous-anneau de $\E^\ent$ que l'on va considérer est
bien entendu $\Sk$, qui est déjà l'anneau qui intervenait dans la
définition de $\varphi$-réseaux. L'égalité $\Sk = \E^\ent \cap W(R)$
nous conduit à choisir comme sous-anneau de ${\E^\ent_\tau}$,
l'intersection ${\E^\ent_\tau} \cap W(R)$ que l'on note $\Sk_\tau$.
Celle-ci est stable par le Frobenius et l'action de $G_K$.

\begin{deftn}
\label{def:phitaures}
Soit $M$ un $(\varphi,\tau)$-module sur $(\E^\ent, {\E^\ent_\tau})$. Un
\emph{$(\varphi,\tau)$-réseau} dans $M$ est la donnée d'un
$\varphi$-réseau $\frakM$ dans $M$ tel que le sous-module $\Sk_\tau
\otimes_{\Sk} \frakM$ de ${\E^\ent_\tau} \otimes_{\E^\ent} \frakM$ soit
stable par l'opérateur $\tau$.

Pour $U \in W(R)$ et $h \in \N$, on dit que $\frakM$ est de hauteur
divisant $U$ (resp. de $U$-hauteur $\leq h$) s'il l'est en tant que
$\varphi$-module.
\end{deftn}

\begin{rem}
Comme $-1$ est la limite dans $\Z_p$ d'une suite d'entiers positifs, la
définition ci-dessus implique que tout $(\varphi,\tau)$-réseau $\frakM$
de $M$ est tel que $\Sk_\tau \otimes_\Sk \frakM$ est stable par
$\tau^{-1}$.
\end{rem}

\subsubsection{Un théorème d'existence dans le cas de torsion}
\label{subsec:existres}

Le but de ce paragraphe est de démontrer la proposition suivante.

\begin{prop}
\label{prop:existres}
Tout $(\varphi,\tau)$-module sur $(\E^\ent, \E^\ent_\tau)$ qui est
annulé par une puissance de $p$ admet un $(\varphi,\tau)$-réseau.
\end{prop}

Pour cela, on fixe $M$ un $(\varphi,\tau)$-module sur $(\E^\ent,
\E^\ent_\tau)$ annulé par une puissance de $p$ et on considère
l'ensemble $\calF_{\varphi}(M)$ (resp. $\calF_{\varphi,\tau}(M)$) des
$\varphi$-réseaux (resp. $(\varphi, \tau)$-réseaux à l'intérieur de $M$.
On veut démontrer que $\calF_{\varphi,\tau}(M)$ est non vide. On sait
déjà (voir lemme \ref{lem:resmodpn}) que $\calF_{\varphi}(M)$ n'est pas
vide. Pour tout entier $h$, soit $\calF_{\varphi}^{\leq h} (M)$ le
sous-ensemble de $\calF_{\varphi}(M)$ formé des $\varphi$-réseaux de
$u$-hauteur $\leq h$. Étant donné que $M$ est annulé par une puissance
de $p$, on sait par la deuxième partie du lemme \ref{lem:resmodpn} que
$\calF_{\varphi}(M)$ est la réunion des $\calF_{\varphi}^{\leq h} (M)$.
On en déduit qu'il existe un entier $h$ tel que $\calF_{\varphi}^{\leq
h}(M)$, lui-même, est non vide. Par ailleurs, la relation d'inclusion
fait de $\calF_{\varphi}^{\leq h}(M)$ un ensemble partiellement ordonné
dont les propriétés essentielles sont dégagées dans \cite{carliu-smf},
\S 3.2. En particulier, le corollaire 3.2.6 nous apprend que
$\calF_{\varphi} ^{\leq h}(M)$ admet un plus petit et un plus grand
élément. La proposition \ref{prop:existres} que nous voulons démontrer
résulte ainsi directement du lemme suivant.

\begin{lemme}
\label{lem:stabtau}
Le plus petit élément de $\calF_{\varphi}^{\leq h}(M)$ appartient à
$\calF_{\varphi,\tau}^{\leq h}(M)$.
\end{lemme}

\begin{proof}
Pour cette démonstration, nous allons avoir besoin de travailler avec
une version partiellement déperfectisée des $(\varphi,\tau)$-modules
(voir \S \ref{subsec:variante2}). Plus précisément, on considère
l'anneau $\E_{u\tnp,\tau}^\ent$ introduit dans ce numéro et on pose
$\Sk_{u\tnp,\tau} = W(R) \cap \E_{u\tnp,\tau}^\ent$.

Soit $\frakM$ le plus petit élément de $\calF_{\varphi}^{\leq h}(M)$.
C'est en particulier un $\varphi$-réseau dans $M$ de $u$-hauteur $\leq
h$. Nous allons montrer que $\Sk_{u\tnp,\tau} \otimes_\Sk \frakM$ est
stable par $\tau$ ce qui impliquera, à l'évidence, que $\Sk_\tau
\otimes_\Sk \frakM$ est, lui aussi, stable par $\tau$ comme voulu. On
pose : $$\frakM' = M \cap \:\tau^{-1}(\Sk_{u\tnp,\tau} \otimes_\Sk
\frakM).$$ Du fait que $\varphi$ et $\tau$ commutent, on déduit
facilement que $\frakM'$ est un $\varphi$-réseau dans $M$. Pour
conclure, il suffit de montrer qu'il est de $u$-hauteur $\leq h$. En
effet, par minimalité de $\frakM$, on aura alors $\frakM \subset
\frakM'$ d'où, en appliquant $\tau$, on trouvera $\tau(\frakM) \subset
\tau(\frakM') \subset \Sk_{\np,\tau} \otimes_\Sk \frakM$. Reste donc à
montrer que $\frakM'$ est de $u$-hauteur $\leq h$. Soit $x \in \frakM'$.
Comme $M$ et $\Sk_{u\tnp,\tau} \otimes_\Sk \frakM$ sont de $u$-hauteur
$\leq h$, on a les écritures suivantes :
\begin{equation}
\label{eq:uhx1}
u^h x = \sum_{i=0}^{p-1} u^i \varphi(x_i) \quad \text{et} \quad
u^h \tau(x) = \sum_{i=0}^{p-1} u^i \varphi(y_i)
\end{equation}
avec $x_i \in M$ et $y_i \in \Sk_{u\tnp,\tau} \otimes_\Sk
\frakM$. On déduit de la deuxième de ces égalités et de la formule
$\tau(u) = [\ueps] u$, que :
\begin{equation}
\label{eq:uhx2}
u^h x = \sum_{i=0}^{p-1} u^i \varphi([\ueps]^{(h-i)/p} y_i).
\end{equation}
avec, à nouveau, $[\ueps]^{(h-i)/p} y_i \in \Sk_{u\tnp,\tau} \otimes_\Sk
\frakM$. Or, du fait que $\Sk_{u\tnp,\tau}$ est un
$\varphi(\Sk_{u\tnp,\tau})$-module libre de base $1, u, \ldots, u^{p-1}$
(car c'est le cas après réduction modulo $p$), on déduit des deux
expressions de $u^h x$ données respectivement par \eqref{eq:uhx2} et la
première égalité de \eqref{eq:uhx1}, que $x_i = [\ueps]^{(h-i)/p} y_i$
pour tout $i \in \{0, \ldots, p-1\}$. Il s'ensuit que $x_i \in \frakM'$,
ce qui démontre que $\frakM'$ est de $u$-hauteur $\leq h$ comme annoncé.
\end{proof}

\subsubsection{Réduction modulo $UV$ des $(\varphi,\tau)$-réseaux}
\label{subsec:phitaumodUV}

On fixe un élément $U \in \Sk_\tau$ qui n'est pas multiple de $p$ et on 
se donne un $\varphi$-réseau $\frakM$ de hauteur divisant $U$ dans un 
$\varphi$-module étale sur $\E^\ent$ (qui peut, oui ou non, avoir de la 
torsion). D'après le lemme \ref{lem:existV}, il existe $V \in W(R)$ 
tel que $\varphi(V) = UV$. En réalité, un examen de la démontration de
ce lemme montre que $V$ peut être choisir dans $\Sk_\tau$ (et même, en
fait, qu'il ne peut être choisi en dehors de $\Sk_\tau$). On suppose,
dans ce numéro, que $V \in \Sk_\tau$.
On pose en outre $\Sk_\tau^+ = \Sk_\tau \cap W(\m_R)$ ; c'est un idéal 
de $\Sk_\tau$ que l'on peut encore définir comme l'idéal noyau du 
morphisme $\Sk_\tau \to W(\bar k)$ déduit de la projection canonique 
$W(R) \to W(\bar k)$.
Étant donné un entier $p$-adique $a$, ainsi qu'un élément $X \in W(R)$, 
on définit $\Tau_X(\frakM)$ comme l'ensemble des automorphismes 
$\tau^a$-semi-linéaires
$$\tau^{(a)}_X : (\Sk_\tau/X \Sk_\tau^+) \otimes_\Sk \frakM 
\to (\Sk_\tau/X \Sk_\tau^+) \otimes_\Sk \frakM$$
commutant à $\varphi$ et tels que pour tout $g \in G_\infty/H_\infty$, 
on ait
\begin{equation}
\label{eq:commuttauX}
(g \otimes \id) \circ \tau^{(a)}_X = (\tau^{(a)}_X)^b \circ 
((\tau^{ab} g \tau) \otimes \id).
\end{equation}
où $b$ est l'unique élément de $\Z_p$ tel que $\chi(g) \cdot [b]
_{\chi(\tau)} = [a]_{\chi(\tau)}$.
Si $X$ divise $Y$, la réduction modulo $X \Sk_\tau^+$ définit une
application canonique $\rho_{Y,X} : \Tau_Y(\frakM) \to \Tau_X(\frakM)$.
Par ailleurs, $\Tau_0(\frakM)$ est exactement l'ensemble des $\tau$ 
qui font de $\frakM$ un $(\varphi,\tau)$-réseau. La proposition
suivante, qui est très semblable à la proposition \ref{prop:rhoUVV},
montre que $\Tau_0(\frakM)$ se retrouve entièrement à partir du
morphisme $\rho_{UV,V} : \Tau_{UV}(\frakM) \to \Tau_{V}(\frakM)$.

\begin{prop}
\label{prop:relevement}
On se donne $\frakM$ un $\varphi$-réseau de hauteur divisant $U$ (avec 
$U$ comme précédemment). L'application $\rho_{0,V} : \Tau_0(\frakM) \to 
\Tau_V(\frakM)$ est injective et son image s'identifie à celle de 
$\rho_{UV,V} : \Tau_{UV}(\frakM) \to \Tau_{V}(\frakM)$.
\end{prop}

\begin{proof}
En utilisant la commutation de $\tau^{-a}$ et $\varphi$, on démontre tout
de suite que $\frakM$ est de hauteur divisant $\tau^{-a}(U)$ pour tout
$a \in \Z_p$. En reprenant la preuve de la proposition
\ref{prop:rhoUVV}, on démontre alors qu'étant donné un élément
$\tau^{(a)}_{UV} \in \Tau_{UV}(\frakM)$, il existe une unique
application $\tau^a$-semi-linéaire $\tau^{(a)} : \Sk_\tau \otimes_\Sk
\frakM \to \Sk_\tau \otimes_\Sk \frakM$ qui commute à $\varphi$ et est
congrue à $\tau^{(a)}_{UV}$ modulo $V \Sk_\tau^+$. L'injectivité de
$\rho_{0,V}$ résulte de cela, tandis pour établir l'égalité annoncée
entre les images, il suffit de prouver que l'application $\tau^{(a)}$
construite vérifie la relation \eqref{eq:commuttauX}.

Pour cela, le plus rapide est sans doute de remarquer qu'étant donné un
élément $b \in \Z_p$, on peut appliquer ce qui précède avec $ab$ à la
place de $a$. On obtient comme cela en particulier l'existence d'une
unique application $\tau^{ab}$-semi-linéaire $\tau^{(ab)} : \Sk_\tau
\otimes_\Sk \frakM \to \Sk_\tau \otimes_\Sk \frakM$ qui est congrue à
$(\tau_0^{(a)})^b$ modulo $V \Sk_\tau^+$ et qui commute à $\varphi$. 
Or, pour $g \in G_\infty/H_\infty$ et $b$ tel que $\chi(g) \cdot
[b]_{\chi(\tau)} = [a]_{\chi(\tau)}$, les deux applications $\tau^b$ et
$(g \otimes \id) \circ \tau \circ ((\tau^{-b} g \tau)^{-1} \otimes \id)$
sont $\tau^{ab}$-semi-linéaires, congrues à $(\tau_0^{(a)})^b$ modulo
$I_0$ et commutent à $\varphi$ ; on en déduit donc bien qu'elles
coïncident.
\end{proof}

La proposition \ref{prop:relevement} implique que, si $\frakM$ et 
$\frakM'$ sont deux $(\varphi,\tau)$-réseaux de hauteur divisant $U$,
alors on a l'égalité :
\begin{equation}
\label{eq:homUV}
\Hom_{\varphi,\tau}(\frakM,\frakM') = \text{image}\big(
\Hom_{\varphi,\tau,UV}(\frakM,\frakM') \to
\Hom_{\varphi,\tau,V}(\frakM,\frakM')\big)
\end{equation}
où, si $X$ est un élément de $\Sk_\tau$, la notation $\Hom_{\varphi,
\tau,X}(\frakM,\frakM')$ désigne l'ensemble des applications $f :
\frakM \to \frakM'$ qui commutent à $\varphi$ et dont la réduction
modulo $X \Sk_\tau^+$ commute à $\tau \mod X \Sk_\tau^+$.
On en déduit que le foncteur de réduction modulo $UV \Sk_\tau^+$ 
induit un foncteur pleinement fidèle de la catégorie des $(\varphi, 
\tau)$-réseaux sur $(\Sk,\Sk_\tau)$ de hauteur divisant $U$ dans la 
catégorie des $(\varphi,\tau)$-réseaux sur $(\Sk,\Sk_\tau/UV\Sk_\tau^+)$
où les morphismes sont définis par la formule du membre de droite de
\eqref{eq:homUV}.

\subsubsection{L'action de $\tau$ sur un $(\varphi,\tau)$-réseau}
\label{subsec:formetau}

L'objectif de ce numéro est de démontrer le théorème suivant (qui jouera 
un rôle clé dans la suite) qui donne un contrôle sur l'action de $\tau$ 
agissant sur un $(\varphi, \tau)$-réseau en fonction de la hauteur de 
celui-ci.

\begin{theo}
\label{theo:congr}
Il existe une constante $c'(K)$ ne dépendant que de $K$ telle que
pour toute donnée de
\begin{itemize}
\item deux éléments $U$ et $V$ de $W(R)$ non multiples de $p$ tels que 
$\varphi(V) = UV$, et
\item un $(\varphi,\tau)$-réseau $\frakM$ de hauteur divisant $U$ à
l'intérieur d'un $(\varphi,\tau)$-module $M$,
\end{itemize}
on ait, en notant $s_0$ le plus petit entier $\geq \log_p(v_R(U \mod
p)) + c'(K)$, l'inclusion suivante :
\begin{equation}
\label{eq:controle}
(\tau^{p^s}- \id) (\frakM) \subset \big( V \Sk_\tau^+ +
p^{s-s_0} \Sk_\tau \big) \otimes_\Sk \frakM
\end{equation}
pour tout entier $s \geq s_0$.
\end{theo}

\begin{rem}
On sait par le lemme \ref{lem:existV} que, dès que l'on se donne un 
élément $U \in W(R)$, il existe $V$ vérifiant la condition requise. En 
outre, d'après la remarque \ref{rem:existV}, l'idéal engendré par $V$ 
--- et donc en particulier l'inclusion \eqref{eq:controle} --- ne dépend 
que de $U$. Ainsi, l'élément $V$ ne joue pas véritablement de rôle dans 
le théorème précédent, hormis le fait qu'il permet de l'énoncer 
relativement simplement.
\end{rem}

\begin{proof}
On se donne des éléments $U,V \in W(R)$ tels que $U$ ne soit pas 
multiple de $p$ et $\varphi(V) = UV$. On fixe un $(\varphi,\tau)$-module 
$M$ sur $(\E^\ent, {\E^\ent_\tau})$ ainsi qu'un $(\varphi,\tau)$-réseau 
$\frakM$ de $M$ de hauteur divisant $U$. Comme d'habitude, on désigne 
par la lettre $T$ la représentation de $G_K$ associée à $M$.  On suppose 
en outre pour commencer que $\frakM$ est annulé par $p^n$ pour un 
certain entier $n$, et on considère un nombre entier $t$ vérifiant les 
deux conditions suivantes :
\begin{itemize}
\item on a $\tau^{p^t} (u) \equiv u \pmod{p^n, U V \Sk_\tau^+}$ ;
\item il existe un élément $\gamma_t \in G_\infty$ agissant trivialement
sur $T$ et tel que $v_p(\chi(\gamma_t) - 1) = t$ (où $v_p$ est la
valuation $p$-adique normalisée par $v_p(p) = 1$).
\end{itemize}
On expliquera, dans la suite de la démonstration, comment déterminer
un tel $t$, mais pour l'instant on se contente de supposer qu'il
est donné. Afin de minimiser les confusions, on notera tout au long de
la preuve $\tau_\frakM$ l'automorphisme $\tau$ agissant sur $\Sk_\tau
\otimes_\Sk \frakM$.
De $\tau^{p^t} (u) \equiv u \pmod{p^n, U V \Sk_\tau^+}$, 
on déduit que l'application $\tau^{p^t}$-semi-linéaire
$$\tau_{\frakM,\triv}^{(p^t)} = \tau^{p^t} \otimes \id :
\frac{\Sk_\tau}{U V \Sk_\tau^+} \otimes_\Sk \frakM \to
\frac{\Sk_\tau}{U V \Sk_\tau^+} \otimes_\Sk \frakM$$
est bien définie. La proposition \ref{prop:relevement} s'applique 
et montre qu'il existe une
application $\tilde \tau_{\frakM,\triv}^{(p^t)} : \Sk_\tau \otimes_\Sk
\frakM \to \Sk_\tau \otimes_\Sk \frakM$ qui commute avec $\varphi$,
vérifie la relation \eqref{eq:commuttauX} et est congrue à
$\tau_{\frakM,\triv}^{(p^t)}$ modulo $V \Sk_\tau^+$.
L'automorphisme $\tilde \tau_{\frakM,\triv}^{(p^t)}$ prolongé à
${\E^\ent_\tau} \otimes_{\E^\ent} M$ définit une structure de
$(\varphi,\tau^{p^t})$-module sur $\frakM$ à laquelle il correspond une
représentation $T_\triv$ du sous-groupe $G_t$. Puisque
les représentations $T$ et $T_\triv$ sont liées au même
$\varphi$-module, leurs restrictions à $G_\infty$ coïncident.

On rappelle que l'on a supposé l'existence d'un élément $\gamma_t \in 
G_\infty$ qui agit trivialement sur $T$ et qui est tel que $v_p 
(\chi(\gamma_t) - 1) = t$. L'élément $\tau^{\chi_\tau(\gamma_t)-1}$ est 
égal au produit $\gamma_t \tau \gamma_t^{-1} \tau^{-1}$ dans le quotient 
$G_K / H_\infty$ ; autrement dit, il existe $h \in H_\infty$ tel 
que $\tau^{\chi_\tau(\gamma_t)-1} = \gamma_t \tau \gamma_t^{-1} 
\tau^{-1} h$. On en déduit que $\tau^{\chi_\tau(\gamma_t)-1}$ agit sur 
$T$ de la même manière que $h$. Puisque $\gamma_t$ appartient à 
$G_\infty$, il agit également trivialement sur $T_\triv$ et, en 
reprenant le raisonnement précédent, on trouve que 
$\tau^{\chi_\tau(\gamma_t)-1}$ agit sur $T_\triv$ de la même façon
que $h$. Au final, le sous-groupe de $G_K$ engendré par $G_\infty$
et $\tau^{\chi_\tau(\gamma_t)-1}$ agit donc de la même façon sur $T$
et $T_\triv$. Or, par le lemme \ref{lem:bijanalogue}, on sait que la 
valuation $p$-adique de $\chi_\tau(\gamma_t) - 1$ est égale à $t$ ; le
sous-groupe engendré par $G_\infty$ et $\tau^{\chi_\tau(\gamma_t)-1}$
n'est donc rien d'autre que $G_t$. Nous venons ainsi de démontrer que
$T_\triv$ s'identifie à la représentation $T$ restreinte au sous-groupe 
$G_t$.
Au niveau des $(\varphi,\tau)$-modules, cela nous donne l'égalité
$\tau_{\frakM}^{p^t} = \tilde \tau_{\frakM,\triv}^{(p^t)}$, d'où
il résulte notamment que $\tau_{\frakM}^{p^t} \equiv \tau^{p^t}
\otimes \id \pmod {V\Sk_\tau^+}$. Pour $x \in \frakM$, on obtient
en particulier $\tau_{\frakM}^{p^t} (x) \equiv x \pmod {V
\Sk_\tau^+}$. Cette congruence reste clairement vraie si $t$ est
remplacée par un entier $s \geq t$. Autrement dit, on a l'inclusion
\begin{equation}
\label{eq:inclpreuve}
(\tau_\frakM^{p^s} - \id) (\frakM) \subset V \Sk_\tau^+ 
\otimes_\Sk \frakM
\end{equation}
valable pour tout $s \geq t$.

\medskip

Déterminons maintenant des conditions simples sur le nombre $t$ 
permettant d'assurer qu'il satisfait aux deux hypothèses faites au début 
de la démonstration. On a en fait déjà étudié la deuxième hypothèse au 
\S \ref{subsec:borneGK} (voir en particulier la formule 
\eqref{eq:inegs}) : on a montré qu'elle est vérifiée dès que $$t \geq 
s_0(K) \quad \text{et} \quad t > \frac{c_0(K)} e + \textstyle 
\max\big(\frac 1{p-1}, n + \log_p(\frac h e)\big)$$ où $s_0(K)$ et 
$c_0(K)$ étaient des constantes ne dépendant que du corps $K$.

Examinons maintenant la première condition. Elle stipule que 
$\tau^{p^t}(u) \equiv u \pmod{p^n, UV \Sk_\tau^+}$, c'est-à-dire que la 
réduction modulo $p^n$ de $\tau^{p^t}(u) - u$ est dans l'idéal $U V 
\cdot W_n(\m_R)$. Or la différence $\tau^{p^t}(u) - u$ se calcule 
facilement ; elle vaut $u ([\ueps]^{p^t} - 1) \in W(R)$ où on rappelle 
que $\ueps$ est l'élément de $R$ défini par le système compatible 
$(\zeta_{p^s})$ de racines primitives $p^s$-ièmes de l'unité, qui avait 
été fixé au début de l'article. En particulier, étant donné que 
$v_R(\ueps) = \frac p{p-1}$, on voit clairement sur cette écriture que 
$\tau^{p^t}(u) - u \in W(\a_R^{> p^{t+1}/(p-1)})$. Par le lemme 
\ref{lem:incWitt}, on en déduit que la réduction modulo $p^n$ de 
$\tau^{p^t}(u) - u$ est dans l'idéal $U V \cdot W_n(\m_R)$ dès que $p^t 
\geq h p^{n-1}$, c'est-à-dire dès que $t \geq n - 1 + \log_p h$.

Au final, en mettant ensemble les résultats des deux alinéas précédents, 
on trouve qu'il existe une constante $c'(K)$ ne dépendant que du corps 
$K$ telle que tout nombre entier $t \geq c'(K) + n + \log_p h$ vérifie 
les hypothèses requises. L'inclusion \eqref{eq:inclpreuve} est donc 
vraie pour de tels $t$.

\medskip

On revient à présent dans la situation du théorème : on se donne un 
$(\varphi,\tau)$-réseau $\frakM$ de hauteur divisant $U$ à l'intérieur 
d'un $(\varphi,\tau)$-module $M$, on appelle $s_0$ le plus petit entier 
supérieur ou égal à $c'(K) + \log_p h$ et on se donne un entier $s \geq 
s_0$. On pose $n = s - s_0$. L'inclusion \eqref{eq:inclpreuve} 
s'applique alors avec les paramètres $\frakM / p^n \frakM$ (qui est 
manifestement annulé par $p^n$) et $s$ et donne :
$$(\tau_\frakM^{p^s} - \id) (\frakM/p^n\frakM) \subset V \Sk_\tau^+ 
\otimes_\Sk \frakM/p^n\frakM$$
soit, encore :
$$(\tau_\frakM^{p^s} - \id) (\frakM) \subset (V \Sk_\tau^+ + p^n
\Sk_\tau) \otimes_\Sk \frakM$$
c'est-à-dire ce que l'on souhaitait démontrer.
\end{proof}

\begin{rem}
En supposant que $\frakM$ est annulé par $p^n$ (et bien sûr toujours 
qu'il est de hauteur divisant $U$), l'inclusion \eqref{eq:inclpreuve} 
--- qui était la clé de la démonstration ci-dessus --- s'étend aux 
entiers $s < t$ sous la forme suivante : 
$$(\tau_\frakM^{p^s} - \id) (\frakM) \subset (W_n(\O_{F_m}^\perf)
+ V \Sk_\tau^+) \otimes_\Sk \frakM$$
où, pour un entier $m$ donné, on a noté $\O_{F_m}^\perf$ l'anneau des 
entier du perfectisé de $(F_0^\sep)^{H_m}$ avec $H_m = \Gal(\bar K / 
K_\infty(\zeta_{p^m}))$.
Pour démontrer cela, on fixe un entier $s < t$, et on considère un 
élément $\gamma \in G_\infty$ tel que $v_p(\chi(\gamma) - 1) \geq t - 
s$. Le commutateur $\gamma \tau^{p^s} \gamma^{-1} \tau^{-p^s}$ est alors 
égal à $\tau^{p^s (\chi_\tau(\gamma) - 1)}$ dans le quotient $G_K / 
H_\infty$ et appartient donc à $G_t$. Ainsi, d'après ce qui a été vu 
dans la démonstration précédente, cet élément agit trivialement sur 
$\frakM$ modulo $V \Sk_\tau^+$. Il en résulte que pour $x \in 
\frakM$, on a $(\gamma \otimes \id) \circ \tau^{p^s} (x) \equiv 
\tau^{p^s}(x) \pmod {V\Sk_\tau^+}$. Comme ceci est vrai pour tout 
$\gamma$ tel que $v_p(\chi(\gamma) - 1) \geq t - s$, on en déduit 
l'inclusion annoncée.
\end{rem}

\section{Le cas des représentations de $E(u)$-hauteur finie}
\label{sec:Eu}

L'élement $E(u)$ qui apparaît dans la locutation \og $E(u)$-hauteur
finie \fg\ est le polynôme minimal de l'uniformisante $\pi$ sur
$W[1/p]$. Il s'agit d'un polynôme d'Eisenstein de degré $e$. C'est donc
en particulier un élément de $W(R)$ qui n'est pas multiple de $p$.
Ainsi, il fait bien sens de parler de $\varphi$-modules de
$E(u)$-hauteur $\leq r$ pour un certain entier $r$ ou de $E(u)$-hauteur
finie. Une $\Z_p$-représentation de $G_\infty$ ou de $G_K$ est dite de
$E(u)$-hauteur $\leq r$ (resp. de $E(u)$-hauteur finie) si le
$\varphi$-module qui lui est associé l'est, tandis qu'une
$\Q_p$-représentation est ainsi qualifiée si elle admet un réseau stable
par Galois satisfaisant à cette propriété.

Reprenant des idées de Breuil, Kisin a élaboré dans \cite{kisin} une
théorie prometteuse pour étudier les représentations semi-stables (ainsi
que leurs déformations) et a, au passage, montré que celles-ci
entretenaient un lien étroit avec les représentations de $E(u)$-hauteur
finie. Précisément, il a démontré qu'une représentation semi-stable dont
tous les poids de Hodge-Tate sont dans $\{0, 1, \ldots, r\}$ est de
$E(u)$-hauteur $\leq r$. Or, la théorie que nous avons développée dans
cet article donne une description des représentations de $G_K$ de
$E(u)$-hauteur $\leq r$ en termes de $(\varphi,\tau)$-modules. Dans ce
dernier chapitre, nous étudions comment cette nouvelle description
s'insère dans la théorie de Kisin. Comme corollaire, nous obtenons le
théorème \ref{introtheo:Eu} de l'introduction, ainsi qu'une description,
en termes de la théorie de Kisin, des réseaux dans les représentations
semi-stables.

\subsection{Les $(\varphi,\tau)$-modules comme complément de la théorie
de Kisin}

Le polynôme $E(u)$ étant un polynôme d'Eisenstein, son coefficient
constant $E(0)$ s'écrit sous la forme $pc$ où $c$ est un élément
inversible dans $\Z_p$. On pose $U = \frac{E(u)} c$ ; il est alors clair
qu'être de $E(u)$-hauteur $\leq r$ est équivalent à être de $U$-hauteur
$\leq r$. Soit $\t$ un élément de $W(R)$ tel que $\varphi(\t) = U\t$ (un
tel élément existe bien d'après le lemme \ref{lem:existV}).

\subsubsection{Un concentré de théorie de Kisin}
\label{subsec:kisin}

On commence par quelques rappels très sommaires concernant l'anneau
$A_\cris$ de Fontaine. Il est défini comme le complété $p$-adique de
l'enveloppe à puissances divisées de $W(R)$ par rapport à l'idéal
principal engendré par $E(u)$ (et compatibles avec les puissances
divisées canoniques sur l'idéal $(p)$). La série
$$\log([\ueps]) = \sum_{n \geq 1} \frac{(1-[\ueps])^n} n$$
converge dans $A_\cris$ vers un élément $t$ qui joue un rôle très
important en théorie de Hodge $p$-adique. Il vérifie $\varphi(t) = pt$
et $\gamma(t) = \chi(\gamma) t$ pour tout $\gamma \in G_K$.

\medskip

On en vient à présent aux rappels sur la théorie de \cite{kisin}. On
se borne à présenter uniquement ce qui sera utile dans la suite. Soit
$\O$ l'anneau des séries en la variable $u$ convergeant sur le disque
ouvert de centre $0$ et de rayon $1$. Il contient clairement $\Sk [1/p]$
et se plonge dans $A_\cris$ en envoyant comme d'habitude $u$ sur
$[\upi]$. L'élément $\lambda$ défini par le produit convergeant 
suivant :
$$\lambda = \prod_{n=0}^\infty \varphi^n\Big(\frac{E(u)}{E(0)}\Big)
= \prod_{n=0}^\infty \varphi^n\Big(\frac U p\Big) \in \O$$
est solution de l'équation $\frac U p \cdot \varphi(\lambda) = \lambda$. Par
ailleurs, d'après l'exemple 5.3.3 de \cite{liu1}\footnote{Les notations
de \cite{liu1} correspondent à celles de cet article si l'on pose
$c = \varphi(\lambda)$.}, il est
possible de choisir $\t$ de façon à ce que l'égalité $\varphi(\lambda \t)
= -t$ soit satisfaite, ce que nous supposerons dans la suite. L'élément
$\lambda$ permet de munir l'anneau $\O$ d'un opérateur de dérivation
$N_\nabla$ défini par $N_\nabla = -u \lambda \frac d{du}$, puis de
définir la catégorie $\Mod^{\varphi,N_\nabla}_{/\O}$ dont un objet 
consiste en la donnée des points suivants :
\begin{itemize}
\item un $\O$-module $\calM$ libre de rang fini ;
\item un morphisme $\varphi$-semi-linéaire $\varphi : \calM
\to \calM$ qui est tel que le conoyau de $\id \otimes \varphi :   
\O \otimes_{\varphi,\O} \calM \to \calM$ soit annulé par une
puissance de $E(u)$ ;
\item un opérateur $N_\nabla : \calM \to \calM$ vérifiant la loi de
Leibniz (\emph{i.e.} $N_\nabla(ax) = N_\nabla(a) x + a N_\nabla(x)$ pour
tout $a \in \O$ et $x \in \calM$) et la relation $N_\nabla \circ \varphi
= p\cdot \frac {E(u)}{E(0)} \cdot \varphi \circ N_\nabla$.
\end{itemize}
Si l'on note $\Rep^\st_{[0,+\infty[}(G_K)$ la catégorie des
représentations semi-stables de $G_K$ à poids de Hodge Tate positifs ou
nuls, un résultat essentiel de \cite{kisin} est la construction d'un
foncteur pleinement fidèle $\calK : \Rep^\st_{[0,+\infty[} (G_K) \to
\Mod^{\varphi,N_\nabla}_ {/\O}$. On note $\Mod^{\varphi,N_\nabla
,0}_{/\O}$ son image essentielle\footnote{Kisin donne en fait une
caractérisation de cette image essentielle en termes de filtrations par
les pentes à la Kedlaya. Nous ne détaillons par ce point dans ces
rappels car il ne nous sera pas utile dans la suite.}. Kisin démontre
encore deux résultats importants pour ce qui va suivre, à savoir que
pour tout objet $\calM$ de $\Mod^{\varphi,N_\nabla ,0}_{/\O}$,
\begin{itemize}
\item[A.] le $\varphi$-module $\E \otimes_\O \calM$ (muni de
l'opérateur $\varphi$ déduit par extension des scalaires) est étale et
correspond \emph{via} la théorie de Fontaine à la représentation
$\calK^{-1}(\calM)$ restreinte à $G_\infty$, et
\item[B.] pour tout réseau $T \subset \calK^{-1}(\calM)$ stable par
$G_\infty$, il existe dans le $\varphi$-module sur $\E^\ent$ correspondant
(qui vit à l'intérieur de $\E \otimes_\O \calM$) un unique
$\varphi$-réseau libre sur $\Sk$ qui est de $E(u)$-hauteur finie.
\end{itemize}
La propriété B précédente permet de construire un foncteur
$\calR_\varphi : \Mod^{\varphi,N_\nabla,0}_{/\O} \to \Mod^\varphi_{/\Sk}
\otimes \Q_p$,
où $\Mod^\varphi_{/\Sk}$ désigne la catégorie des $\varphi$-réseaux de 
$E(u)$-hauteur finie (dans un $\varphi$-module étale sur $\E^\ent$) 
\emph{libres} sur $\Sk$, et où $\Mod^\varphi_{/\Sk} \otimes \Q_p$ est sa 
catégorie à isogénie près\footnote{Cela signifie que les objets sont les 
mêmes, mais que les ensembles de morphismes sont tensorisés par $\Q$}.
Enfin, si l'on note $\calT_\varphi$ le foncteur qui à un $\varphi$-réseau
$\frakM$ fait correspondre la représentation de $G_\infty$ associée à 
$\E \otimes_\Sk \frakM$, et si l'on définit $\Rep^{E(u)}_{[0,+\infty[}
(G_\infty)$ comme la catégorie des représentations de $G_\infty$ qui
sont de $E(u)$-hauteur finie, la propriété A nous dit que la composée
$\calT_\varphi \circ \calR_\varphi \circ \calK$ n'est autre que la    
restriction de l'action de $G_K$ à $G_\infty$.
Le carré commutatif suivant :
\begin{equation}
\label{eq:diagKisin}
\raisebox{0.5\depth}{\xymatrix @C=50pt {
\Rep^\st_{[0,+\infty[}(G_K) \ar[d]^-{\sim}_-{\calK} \ar[r]^-{\res} &
\Rep^{E(u)}_{[0,+\infty[} (G_\infty) \\
\Mod^{\varphi,N_\nabla,0}_{/\O} \ar[r]^-{\calR_\varphi} &
\Mod^\varphi_{/\Sk} \otimes \Q_p
\ar[u]_-{\calT_\varphi}
}}
\end{equation}
dans lequel $\res$ désigne le foncteur de restriction de l'action de
$G_K$ à $G_\infty$, résume de façon concise --- et, qui plus est,
commode pour notre propos --- les résultats de Kisin qui viennent d'être
rappelés.

\subsubsection{Une démonstration directe de l'unicité dans la propriété
B}

La démonstration que Kisin donne dans \cite{kisin} de la propriété B 
énoncée précédemment utilise de façon essentielle des théorèmes 
difficiles portant sur les $\varphi$-modules sur l'anneau $\O$. Dans ce 
paragraphe, nous donnons une démonstration alternative de la partie \og 
unicité \fg, qui est plus élémentaire (et plus simple), dans le sens où 
elle ne fait intervenir à aucun moment ni l'anneau $\O$, ni la théorie 
de Kisin. Nous démontrons même un résultat très légèrement plus général 
qui s'énonce comme suit.

\begin{prop}
\label{prop:propB}
Soient $\frakM_1$ et $\frakM_2$ deux $\varphi$-réseaux libres sur $\Sk$
de $E(u)$-hauteur finie. Soit également un morphisme $f : \E^\ent
\otimes_\Sk \frakM_1 \to \E^\ent \otimes_\Sk \frakM_2$. Alors $f(\frakM_1)
\subset \frakM_2$.
\end{prop}

On rappelle avant d'entamer la démonstration que si $\frakM$ est un 
$\varphi$-réseau (pas nécessairement libre sur $\Sk$) à l'intérieur d'un 
$\varphi$-module $M$ libre sur $\E^\ent$, il résulte du théorème de 
classification des modules sur $\Sk$ que $\Libre(\frakM) = (\E^\ent 
\otimes_\Sk \frakM) \cap \frakM[1/p]$ est un $\varphi$-réseau libre sur 
$\Sk$. En outre, si $\frakM$ est de $E(u)$-hauteur finie, il en est de 
même de $\Libre(\frakM)$ et les $E(u)$-hauteurs sont, dans ce cas, 
égales.

Quitte à remplacer $\frakM_1$ par $\Libre(f(\frakM_1) + \frakM_2)$ (qui
est encore de $E(u)$-hauteur finie comme on le vérifie sans peine), on
peut supposer que $f$ est l'identité et que $\frakM_1$ contient
$\frakM_2$. Il s'agit alors de montrer que $\frakM_1 = \frakM_2$. En
passant aux déterminants, on peut supposer que $\frakM_1$ et
$\frakM_2$ sont libres de rang $1$ sur $\Sk$. Si le vecteur $e_1$ forme
une base de $\frakM_1$, il existe un élément $a \in \Sk$ tel que $e_2 =
a e_1$ soit une base de $\frakM_2$. D'après une variante du théorème de
préparation de Weierstrass (voir par exemple théorème 2.1, chap. 5 de
\cite{lang}), $a$ s'écrit comme le produit d'un élément inversible de
$\Sk$ et d'un polynôme $A(u)$ à coefficients dans $W$ de la forme $A(u)
= u^d + p(a_{d-1} u^{d-1} + \cdots + a_0)$. On veut montrer que $A(u)$
est inversible dans $\Sk$, c'est-à-dire concrètement que $d = 0$. Soit
$A^\sigma$ le polynôme obtenu en appliquant le Frobenius sur $W$ à tous
les coefficients de $A$. On a $\varphi(A(u)) = A^\sigma (u^p)$. Du fait
que, par hypothèse, $\frakM_2$ est un $\varphi$-module de hauteur
$E(u)$-hauteur finie, on déduit l'existence d'un élément $b \in \Sk$ qui
divise une puissance de $E(u)$ dans $\Sk$ et qui vérifie l'égalité
$\varphi(e_2) = b e_2$.

\begin{lemme}
L'élement $b$ s'écrit sous la forme $E(u)^n b'$ où $n$ est un entier
et $b'$ est inversible dans $\Sk$.
\end{lemme}

\begin{proof}
Le fait que $\Sk$ soit un anneau intègre et que l'idéal principal
engendré par $E(u)$ soit premier (puisque le quotient s'identifie à
l'anneau intègre $\O_K$) implique le résultat dans le cas où $b$ est un
diviseur de $E(u)$ (et non d'une puissance de $E(u)$). Pour le cas
général, on raisonne de la même façon par récurrence sur la plus petite
puissance de $E(u)$ que $b$ divise.
\end{proof}

On utilise maintenant que $\frakM_1$ est, lui aussi, par hypothèse, un
$\varphi$-module de $E(u)$-hauteur finie. La formule $\varphi(e_1) =
A^\sigma (u^p) b e_2$ implique que $A(u)$ divise $A^\sigma(u^p) b$
dans $\Sk$ (puisque $\frakM_1$ est stable par $\varphi$) et que le
quotient $\frac {A^\sigma(u^p) b}{A(u)}$ est divisible par une puissance
de $E(u)$ (puisque $\frakM_1$ est de $E(u)$-hauteur finie). En écrivant
$b$ sous la forme $E(u)^n b'$ comme dans le lemme précédent, on obtient
les deux divisibilités suivantes dans $\Sk$ :
$$A(u) \text{ divise } A^\sigma(u^p) E(u)^n
\quad \text{et} \quad
A^\sigma(u^p) \text{ divise } A(u) E(u)^m$$
pour un certain entier $m$. Quitte à augmenter $n$ et $m$, on peut
supposer que ces deux nombres sont strictement positifs.
Comme $\Sk \subset \O$, on peut évaluer 
les séries appartenant à $\Sk$ en tout élément de l'idéal maximal de 
$\O_{\bar K}$. Les divisibilités qui viennent d'être écrites assurent
ainsi que toute racine du polynôme $A(u)$ (resp. $A^\sigma(u^p)$) qui
appartient à l'idéal maximal de $\O_{\bar K}$ est également racine du
produit $A^\sigma(u^p) E(u)^n$ (resp. $A(u) E(u)^m$).
On suppose à présent, par l'absurde, que $d > 0$. La théorie des
polygônes de Newton affirme que les polynômes $A(u)$ et $A^\sigma(u)$
admettent chacun $d$ racines (comptées avec multiplicité) dans l'idéal
maximal de $\O_{\bar K}$ et que les valuations de ces racines, notées
$v_1, \ldots, v_d$, se correspondent deux à deux. On peut bien entendu
supposer que les $v_i$ sont triées par ordre croissant. Si $x$ désigne
une racine de $A(u)$ de valuation $v_d$ et si $v_d < \infty$, on ne peut
avoir $A^\sigma(x^p) = 0$ puisque la valuation de $x^p$ dépasse tous les
$v_i$. Ainsi, on a nécessairement $E(x) = 0$, ce qui signifie que $x$
est un conjugué de $\pi$. En particulier, $v_d$ ne peut valoir que
$+\infty$ ou $1$. De la même façon, si $y$ désigne une racine
$p$-ième d'une racine de $A^\sigma(u)$ de valuation $v_1$, on a $A(y)
\neq 0$, et donc $E(y) = 0$. Il s'ensuit que $v_1 \in \{p, +\infty\}$.
Clairement la seule façon de concilier les contraintes que l'on vient
d'obtenir est d'avoir $v_i = +\infty$ pour tout $i$, c'est-à-dire
$A(u) = u^d$. Mais, ce cas n'est pas non plus possible
car $A^\sigma(u) = u^{pd}$ ne divise manifestement pas $A(u)E(u)^m = u^d
E(u)^m$ si $d$ est strictement positif. Ainsi se termine donc la
démonstration de la proposition \ref{prop:propB}.

\subsubsection{L'apport des $(\varphi,\tau)$-modules}

La théorie des $(\varphi,\tau)$-modules, qui a été développée dans cet
article, fournit des objets qui s'insèrent de façon très naturelle dans
le diagramme \eqref{eq:diagKisin} et permet de le compléter de façon
satisfaisante. C'est ce que nous nous proposons de présenter dans ce
paragraphe. Plus précisément, nous allons construire le diagramme 
suivant :
\begin{equation}
\label{eq:diagphitau}
\raisebox{0.5\depth-0.3cm}{\xymatrix @C=50pt {
\Rep^\st_{[0,+\infty[}(G_K) \ar[d]^-{\sim}_-{\calK} \ar@{^(->}[r] &   
\Rep^{E(u)}_{[0,+\infty[} (G_K) \ar[r]^-{\res} &
\Rep^{E(u)}_{[0,+\infty[} (G_\infty) \\
\Mod^{\varphi,N_\nabla,0}_{/\O} \ar[r]^-{\calR_{\varphi,\tau}}
\ar@/_0.4cm/[rr]_-{\calR_\varphi} &
\Mod^{\varphi,\tau}_{/\Sk} \otimes \Q_p \ar[u]^-{\sim}
_-{\calT_{\varphi,\tau}} \ar[r]^-{\text{oubli de }\tau} &
\Mod^\varphi_{/\Sk} \otimes \Q_p \ar[u]^-{\sim}_-{\calT_\varphi} }}  
\end{equation}
Les notations $\Rep^{E(u)}_{[0,+\infty[} (G_K)$ et
$\Mod^{\varphi,\tau}_{/\Sk}$ qui apparaissent dans le diagramme
\eqref{eq:diagphitau} font respectivement référence à la catégorie des
représentations de $G_K$ de $E(u)$-hauteur finie et à celle des
$(\varphi,\tau)$-réseaux libres sur $\Sk$ de $E(u)$-hauteur finie. Le
foncteur $\calT_ {\varphi,\tau}$, quant à lui, est celui qui fait
correspondre à un $(\varphi,\tau)$-réseau $\frakM$ de $E(u)$-hauteur
finie la représentation de $G_K$ associée au $(\varphi,\tau)$-module
$\E^\ent \otimes_\Sk \frakM$. Le diagramme \eqref{eq:diagphitau} indique 
que c'est une équivalence de catégories, ce qui est effectivement le
cas :

\begin{prop}
\label{prop:equivTphi}
Les foncteurs $\calT_\varphi$ et $\calT_{\varphi,\tau}$ sont des
équivalences de catégories.
\end{prop}

\begin{proof}
La pleine fidélité des deux foncteurs suit de la proposition 
\ref{prop:propB}. L'essentielle surjectivité de $\calT_\varphi$ est 
vraie par définition des représentations de $E(u)$-hauteur finie. Il ne 
reste donc qu'à démontrer l'essentielle surjectivité de 
$\calT_{\varphi,\tau}$, et il suffit pour cela de justifier que si 
$\frakM$ est un $\varphi$-réseau libre sur $\Sk$, de $E(u)$-hauteur 
finie, dans un $(\varphi,\tau)$-module sur $(\E^\ent, \E_\tau^\ent)$, 
alors $\Sk_\tau \otimes_\Sk \frakM$ est automatiquement stable par 
$\tau$.

Soit $h$ un entier tel que $\frakM$ soit de $E(u)$-hauteur $\leq h$. 
Pour tout entier $n$, on pose $\frakM_n = \frakM/p^n\frakM$. Bien 
entendu, $\frakM_n$ est, lui aussi, de $E(u)$-hauteur $\leq h$. Après 
avoir remarque que le quotient $\frac{\tau(E(u))}{E(u)}$ est un élément 
inversible de $\Sk_{u\tnp,\tau}$, on peut reprendre l'argumentation 
développée au \S \ref{subsec:existres} en remplaçant partour \og 
$u$-hauteur $\leq h$ \fg\ par \og $E(u)$-hauteur $\leq h$ \fg\ et 
démontrer, ce faisant, qu'il existe dans $\frakM_n[1/u]$ un 
$(\varphi,\tau)$-réseau $\frakM'_n$ contenu dans $\frakM_n$. D'autre 
part, par le lemme 4.1.2 de \cite{liu1} (appliqués aux duaux de 
$\frakM_n$ et $\frakM'_n$), il existe une constance $c$ indépendante de 
$n$ telle que $\frakM_n \subset p^c \cdot \frakM'_n$. On en déduit que 
$\Sk_\tau \otimes_{\Sk} \frakM_n$ est stable par l'opération $p^c \tau$ 
puis, en passant à la limite, qu'il en est de même de $\Sk_\tau 
\otimes_\Sk \frakM$. Or, par construction, on sait également que $\tau$ 
stabilise également $\E^\ent_\tau \otimes_\Sk \frakM$. Comme 
$\E^\ent_\tau \cap (p^{-c} \Sk_\tau) = \Sk_\tau$, on en déduit que 
$\tau$ stabilise $\Sk_\tau \otimes_\Sk \frakM$. La proposition est 
démontrée.
\end{proof}

On en vient enfin à la construction du foncteur $\calR_{\varphi, \tau} :
\Mod^{\varphi,N_\nabla,0}_{/\O} \to \Mod^{\varphi,\tau}_{/\Sk} \otimes
\Q_p$. On considère pour cela $\calM$ un objet de
$\Mod^{\varphi,N_\nabla,0} _{/\O}$, et on fixe un réseau $T$ de
$\calK^{-1}(\calM)$ stable par $G_K$. Soit $M$ le
$(\varphi,\tau)$-module sur $(\E^\ent, \E_\tau^\ent)$ qui lui est
associé. D'après la propriété B, il existe dans $M$ un unique
$\varphi$-réseau $\frakM$ qui est libre sur $\Sk$ et de $E(u)$-hauteur
finie, et d'après la démonstration de la proposition
\ref{prop:equivTphi}, le produit tensoriel $\Sk_\tau \otimes_\Sk \frakM$
est stable par $\tau$. L'objet $\frakM$ est donc un $(\varphi,
\tau)$-réseau qui, à cause du choix de $T$, n'est défini qu'à isogénie
près. Autrement dit, c'est un objet de la catégorie
$\Mod^{\varphi,\tau}_{/\Sk} \otimes \Q_p$, et on peut donc poser
$\calR_{\varphi,\tau}(\calM) = \frakM$.

\subsection{Un presque quasi-inverse de $\calR_{\varphi,\tau}$}
\label{subsec:qinvRphitau}

\subsubsection{Énoncé des résultats}

Du diagramme \eqref{eq:diagphitau}, il suit immédiatement que le foncteur 
$\calR_{\varphi,\tau}$ est pleinement fidèle. En contrepartie, il n'est 
généralement pas essentiellement surjectif. Ceci traduit le fait qu'il 
existe des représentations de $E(u)$-hauteur finie qui ne sont pas 
semi-stables. Pour construire un contre-exemple, on choisit $K = 
\Q_p(\sqrt[p] 1)$ et on pose $L = K(p_1)$ où $p_1$ est une racine 
$p$-ième de $p$ fixée. L'extension $L/K$ est alors galoisienne et son 
groupe de Galois $\Gal(L/K)$ est cyclique d'ordre $p$. On peut donc 
faire agir ce dernier non trivialement sur un $\Q_p$-espace vectoriel 
$V_1$ de dimension $p$ (en permutant cycliquement les vecteurs d'une 
base). On obtient de cette façon une représentation $V_1$ de $G_K$ dont 
la restriction à $G_\infty$ est triviale et donc de $E(u)$-hauteur 
finie. La proposition suivante, appliquée avec $V_2 = \Q_p^p$ munie
de l'action triviale, montre qu'elle n'est cependant pas semi-stable.

\begin{prop}
\label{prop:sstextension}
Soient $V_1$ et $V_2$ deux représentations semi-stables de $G_K$. On
suppose qu'il existe une extension \emph{totalement} ramifiée $L/K$
et un isomorphisme $G_L$-équivariant $f : V_1 \to V_2$. Alors $f$ est
$G_K$-équivariant.
\end{prop}

\begin{proof}
Exercice.
\end{proof}

Dans l'exemple que l'on vient de présenter, on observe que, certes, la 
représentation $V_1$ n'est pas semi-stable mais qu'elle coïncide 
néanmoins avec une représentation semi-stable sur le sous-groupe $G_1$. 
Dans le langage des $(\varphi,\tau)$-modules, cela signifie que, certes, 
l'objet $\frakM$ de $\Mod^{\varphi,\tau}_{/\Sk} \otimes \Q_p$ 
correspondant à $V_1$ n'est pas dans l'image essentielle de 
$\calR_{\varphi,\tau}$ mais qu'il existe néanmoins un objet $\calM \in 
\Mod^{\varphi,N_\nabla,0}_{/\O}$ dont l'image par $\calR_{\varphi,\tau}$ 
est isomorphe à $\frakM$ comme $(\varphi,\tau^p)$-module. Il s'avère que 
ceci est un phénomène général, comme le précise le théorème suivant (qui 
sera démontré dans les \S\S \ref{subsec:constcalS} et 
\ref{subsec:propcalS}).

\begin{theo}
\label{theo:calS}
On suppose que $K$ est une extension finie de $\Q_p$\footnote{Il
est probable que le théorème demeure sans cette hypothèse.}. Alors,
il existe un unique foncteur $\calS_{\varphi,N_\nabla} : 
\Mod^{\varphi,\tau}_{/\Sk} \otimes \Q_p \to 
\Mod^{\varphi,N_\nabla,0}_{/\O}$ tel que :
\begin{itemize}
\item pour tout objet $\calM \in \Mod^{\varphi,N_\nabla,0}_{/\O}$, il y 
a un isomorphisme canonique 
$$f_\calM : \calS_{\varphi,N_\nabla} \circ \calR_{\varphi, \tau}(\calM) 
\simeq \calM$$ dans la catégorie $\Mod^{\varphi,N_\nabla,0}_{/\O}$ ;
\item pour tout objet $\frakM \in \Mod^{\varphi,\tau}_{/\Sk} \otimes 
\Q_p$, il y a un isomorphisme canonique
$$f_\frakM : \calR_{\varphi, \tau}\circ \calS_{\varphi,N_\nabla}(\frakM) 
\simeq \frakM$$
qui commute à $\varphi$ et à $\tau^{p^s}$ où $s$ désigne le plus grand
entier tel que $K_s/K$ soit galoisienne.
\end{itemize}
Si, de plus, $\frakM$ est de dimension $d$ et si $t$ est un entier tel 
que $p^t(p-1) > d$, alors l'isomorphisme $f_\frakM$ commute également
à $\tau^{p^t}$.
\end{theo}

Le théorème ci-dessus se transpose bien sûr directement dans le langage 
des représentations. Il affirme qu'il existe un unique foncteur $\calS$ 
allant de la catégorie de représentations de $E(u)$-hauteur finie dans 
celle des représentations semi-stables à poids de Hodge-Tate positifs
qui jouït des deux propriétés suivantes :
\begin{itemize}
\item si $V$ est une représentation semi-stable (et donc, en 
particulier de $E(u)$-hauteur finie d'après le résultat de Kisin),
alors $\calS(V) \simeq V$, et
\item si $V$ est une représentation de $E(u)$-hauteur finie, alors
il existe un isomorphisme $V \to \calS(V)$ qui est $G_s$-équivariant
(où $s$ désigne à nouveau le plus grand entier tel que $K_s/K$ soit
galoisienne) et également $G_t$-équivariant pour tout entier $t$ tel
que $p^t(p-1) > \dim_{\Q_p} V$.
\end{itemize}
Ce résultat implique clairement le théorème \ref{introtheo:Eu} de 
l'introduction. Par ailleurs, en combinant ce qui précède avec la 
proposition \ref{prop:sstextension}, on trouve que $V$ est semi-stable 
si, et seulement si $V$ est isomorphe à $\calS(V)$ (en tant que 
$G_K$-représentation) : le foncteur $\calS$ permet donc, en un sens,
de caractériser les représentations semi-stables parmi celles qui
sont de $E(u)$-hauteur finie.

\subsubsection{Logarithmes tronqués}
\label{subsec:logm}

En guise de préliminaire à la démonstration du théorème \ref{theo:calS}, 
nous étudions dans ce numéro les applications \emph{logarithmes 
tronqués} qui joueront un rôle essentiel dans la suite.

Soit $A$ est une $\Q_p$-algèbre topologique (non nécessairement 
commutative). Si $a$ est un élement de $A$ et si $m$ est un entier 
strictement positif, on définit le \emph{logarithme tronqué d'ordre $m$} 
de $a$ par :
$$\log_m a = \sum_{i=1}^{p^m-1} \frac{(1-a)^i} i.$$
Le but de cette sous-partie est de montrer que la fonction $\log_m$ 
vérifie des propriétés sympathiques sous des hypothèses de convergence 
très faibles. Plus précisément, on se donne $\Lambda \subset A$ un
sous-$\Z_p$-module fermé et on introduit la définition suivante.

\begin{deftn}
Un élément $a \in A$ est dit \emph{$\Lambda$-borné à l'ordre $m$} si 
$\Lambda$ est stable par multiplication à gauche par $a$ et si 
$\frac{(1-a)^i} i \in \Lambda$ pour tout $i \in \{1, \ldots, p^m\}$.

Si $a$ est $\Lambda$-borné à tout ordre, on dira simplement que $a$
est \emph{$\Lambda$-borné}.
\end{deftn}

Pour tout entier $i > 0$, on note $\ell(i)$ la partie entière de $\log_p 
i$ et on convient que $\ell(0) = 0$. Si $A$ est $\Lambda$-borné à l'ordre
$m$, on a :
\begin{equation}
\label{eq:Lambdaborne}
(1-a)^i = p^{\ell(i)} \cdot (1-a)^{i-p^{\ell(i)}} \cdot 
\frac{(1-a)^{p^{\ell(i)}}} {p^{\ell(i)}}
\end{equation}
d'où on déduit que $(1-a)^i \in p^{\ell(i)} \Lambda$ pour tout $i \in 
\{1, \ldots, p^m\}$ et $(1-a)^i \in p^m \Lambda$ pour tout $i > p^m$.

\begin{prop}
\label{prop:logmmult}
Soient $a$ et $b$ deux éléments de $A$ $\Lambda$-bornés à l'ordre $m$.
On suppose que $a$ et $b$ commutent.
Alors $ab$ est $\Lambda$-borné à l'ordre $m$ et on a :
$$\log_m(ab) \equiv \log_m a + \log_m b \pmod{p^{m-1} \Lambda}$$
pour tout entier strictement positif $m$.
\end{prop}

\begin{proof}
En élevant à la puissance $i$ l'égalité $1 - ab = (1-a) + a(1-b)$, 
on obtient :
$$(1-ab)^i = \sum_{j=0}^i \binom i j \cdot a^j (1-a)^{i-j} (1-b)^j.$$
En isolant le terme en $j = 0$ et en utilisant l'égalité 
$\frac 1 i \cdot \binom i j = \frac 1 j \cdot \binom{i-1}{j-1}$, il
vient :
$$\frac{(1-ab)^i} i =
\frac{(1-a)^i} i + \sum_{j=1}^i \binom {i-1}{j-1} \cdot a^j (1-a)^{i-j} 
\cdot \frac{(1-b)^j} j$$
d'où il résulte déjà que $ab$ est $\Lambda$-borné à l'ordre $m$. 
Maintenant, en sommant l'égalité précédente sur tous les $i$ compris 
entre $1$ et $p^m-1$, on trouve :
\begin{eqnarray}
\log_m(ab) & = & \log_m a + \sum_{1 \leq j \leq i \leq p^m-1}
\binom {i-1}{j-1} \cdot a^j (1-a)^{i-j} \cdot \frac{(1-b)^j} j \nonumber \\
& = & \log_m a + \sum_{j=1}^{p^m-1} \Bigg( 
a^j \cdot \sum_{i=j}^{p^m-1} \binom {i-1}{j-1} \cdot (1-a)^{i-j} \Bigg)
\cdot \frac{(1-b)^j} j \label{eq:logab}
\end{eqnarray}
Pour étudier le terme dans le parenthèse ci-dessus, on commence par 
remarquer qu'on a l'identité $\sum_{i=0}^\infty X^i = \frac 1{1-X}$ où 
$X$ est une variable formelle. En dérivant $j-1$ fois cette égalité, on 
obtient la nouvelle formule $\sum_{i=j}^\infty \binom{i-1}{j-1} X^{i-j} 
= \frac 1{(1-X)^j}$, soit encore l'identité
$(1-X)^j \cdot \sum_{i=j}^\infty \binom{i-1}{j-1} X^{i-j} = 1$
valable dans l'anneau $\Z[[X]]$ des séries formelles à coefficients
entiers. Étant donné que $1-X$ est inversible dans cet anneau, on en 
déduit la congruence :
$$(1-X)^j \cdot \sum_{i=j}^{p^m-1} \binom{i-1}{j-1} X^{i-j} 
\equiv 1 \pmod{X^{p^m-j}}.$$
En appliquant ce qui précède avec $X = 1-a$ et en se rappelant, d'une 
part, que $\Lambda$ est stable par multiplication par $a$ (et donc 
aussi par $1-a$) et, d'autre part, que, d'après l'égalité 
\eqref{eq:Lambdaborne}, l'élément $(1-a)^{p^m-j}$ appartient à 
$p^{\ell(p^m-j)} \Lambda$, on obtient :
$$a^j \cdot \sum_{i=j}^{p^m-1} \binom {i-1}{j-1} \cdot (1-a)^{i-j} 
\equiv 1 \pmod {p^{\ell(p^m-j)} \Lambda}.$$
En reportant dans \eqref{eq:logab}, on obtient :
$$\log_m(ab) - (\log_m a + \log_m b) \in \sum_{j=1}^{p^m-1}
\frac{p^{\ell(p^m-j) + \ell(j)}} j \cdot \Lambda.$$
Pour conclure, il suffit donc de démontrer que, pour tout $j$ compris 
entre $1$ et $p^m-1$, on a $\ell(p^m-j) + \ell(j) - v_p(j) \geq m-1$ 
(où $v_p$ désigne la valuation $p$-adique). Or, on peut écrire un tel 
entier $j$ sous la forme $j = p^v \cdot j'$ avec $j'$ premier à $j$. Un 
calcul immédiat donne alors :
$$\ell(p^m-j) + \ell(j) - v_p(j) = v + \ell(p^{m-v} - j') + \ell(j').$$
Si $j' \geq \frac{p^{m-v}} 2$, on a $\ell(j') \geq m-v-1$ et l'inégalité 
voulue est bien vérifiée. Si, au contraire, $j' < \frac{p^{m-v}} 2$, on
a $\ell(p^{m-v} - j') \geq m-v-1$ et l'inégalité voulue est de même
démontrée.
\end{proof}

\begin{prop}
\label{prop:logmcont}
Soient $a$ et $b$ deux éléments de $A$ qui sont $\Lambda$-bornés à 
l'ordre $m$. On suppose que $a$ et $b$ commutent et qu'ils sont congrus 
modulo $p^m \Lambda$. Alors $\log_m(a)\equiv \log_m(b) \pmod {p^{m-1} 
\Lambda}$.
\end{prop}

\begin{proof}
La factorisation
$$(1-a)^i - (1-b)^i = (b-a) \cdot \big( (1-a)^{i-1} + (1-a)^{i-2}(1-b)
+ \cdots + (1-b)^{i-2} \big)$$
montre que la différence $\frac{(1-a)^i}i - \frac{(1-b)^i}i$ appartient 
à $\sum_{j=0}^{i-1} p^{m + \ell(j) + \ell(i-j) - v_p(i)} \Lambda$. Or, si
$j \geq \frac i 2$, on a $\ell(j) \geq \ell(i) - 1$ tandis que, dans 
le cas contraire, on a $\ell(i-j) \geq \ell(i) - 1$. Ainsi on a 
toujours $m + \ell(j) + \ell(i-j) - v_p(i) \geq m + \ell(i) - 1 - v_p(i) 
\geq m - 1$ et la proposition en découle.
\end{proof}

\begin{cor}
\label{cor:logmpuis}
Soit $a \in A$ un élément $\Lambda$-borné à l'ordre $m$. On suppose que
$a^{p^s}$ converge vers $1$ dans $A$ quand $s$ tend vers l'infini.
Alors, pour tout $n \in \Z_p$ :
$$\log_m(a^n) \equiv n \cdot \log_m(a) \pmod {p^{m-1} \Lambda}.$$
\end{cor}

\begin{proof}
On remarque, pour commencer, que l'on peut écrire :
$$a^{p^m} = (1 + (a-1))^{p^m} = \sum_{j=0}^{p^m} \binom{p^m}j (a-1)^j.$$
Sachant que, pour tout $j \in \{1, \ldots, p^m\}$, la valuation 
$p$-adique du coefficient binomial $\binom {p^m} j$ est égale à $m - 
v_p(j)$, on déduit de l'égalité précédente et du fait que $a$ est 
$\Lambda$-borné à l'ordre $m$ que $a^{p^m} \equiv 1 \pmod{p^m \Lambda}$.

On choisit à présent un entier $n'$ congru à $n$ modulo $p^m$ et on 
écrit $n = n' + p^m q$ pour un certain $q \in \Z_p$. On a alors $a^n = 
a^{n'} \cdot (a^{p^m})^q$. Par ailleurs, on sait que $a^{p^m}$ s'écrit 
sous la forme $1 + b$ avec $b \in p^m \Lambda$. Comme $\Lambda$ est 
stable par multiplication par $a$, il l'est aussi par multiplication par 
$b$. Comme il est en outre fermé, on trouve en développant $(1+b)^q$, 
que $(a^{p^m})^q \equiv 1 \pmod{p^m \Lambda}$ puis que $a^n \equiv 
a^{n'} \pmod{p^m \Lambda}$.
Par la proposition \ref{prop:logmmult}, il vient alors $\log_m(a^n) 
\equiv \log_m(a^{n'}) \pmod {p^{m-1} \Lambda}$ et une application 
répétée de la proposition \ref{prop:logmcont} montre enfin que
$$\log_m(a^{n'}) \equiv n' \cdot \log_m(a) \equiv n \cdot \log_m(a) 
\pmod {p^{m-1} \Lambda}$$
ce qui permet de conclure.
\end{proof}

\subsubsection{La construction du foncteur $\calS_{\varphi,N_\nabla}$}
\label{subsec:constcalS}

On en vient à présent à la démonstration du théorème \ref{theo:calS}. On
suppose donc que $K$ est une extension finie de $\Q_p$ ; en particulier,
son corps résiduel $k$ est un corps fini.
Étant donné que les catégories $\Mod^{\varphi,\tau}_{/\Sk} \otimes \Q_p$ 
sont équivalentes entre elles pour les différents $\tau$, on peut 
choisir un $\tau$ particulier pour construire $\calS_{\varphi,N_\nabla}$. 
Pour la suite, on en fixe donc un qui appartient à l'inertie sauvage (de 
sorte que $\tau^{p^n}$ converge vers l'identité dans $G_K$) et qui 
vérifie en outre $\chi(\tau) = 1$ (on rappelle qu'un tel élément existe 
toujours par le lemme 5.1.2 de \cite{liu-compo}).

\begin{figure}
\begin{center}
\includegraphics{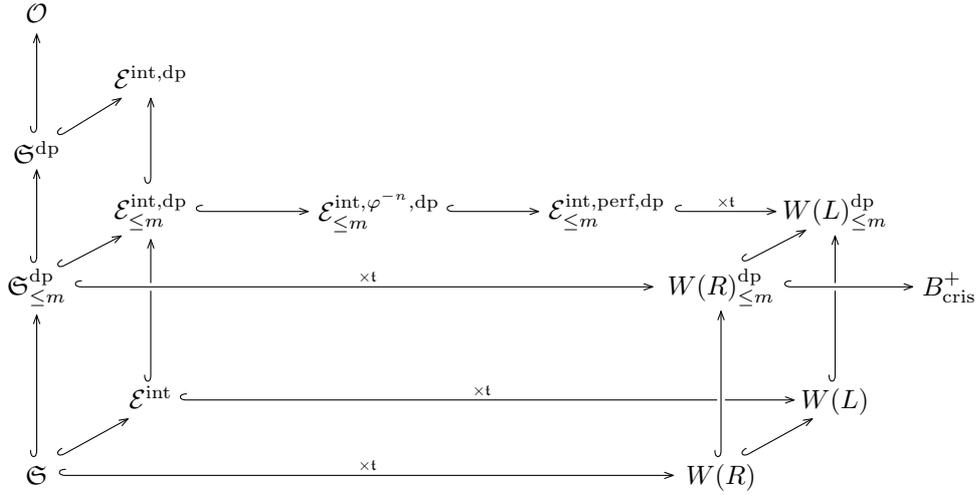}

\caption{Divers $\Sk$-modules intervenant dans la démonstration du
théorème \ref{theo:calS}}
\label{fig:modules}

\bigskip

\begin{minipage}{13cm}
La stratégie de la preuve est la suivante. On construit d'abord des 
versions tronquées $\log_m \tau$ de l'opérateur $\log \tau$ qui 
dépendent d'un entier $m$ et sont définies \emph{a priori} sur l'anneau 
$W(R)^\DP_{\leq m}$. On montre ensuite que, modulo des termes qui 
tendent vers $0$ avec $m$, ces versions tronquées sont définies sur 
$\E^{\ent,\DP}_{\leq m}$ empruntant le chemin suivant : $W(L)^\DP_{\leq 
m}$, $\E^{\ent,\perf,\DP}_{\leq m}$, $\E^{\ent, \varphi^{-n} ,\DP}_{\leq 
m}$ puis $\E^{\ent,\DP}$. On recolle ensuite les $\log_m \tau$ pour 
former un véritable opérateur $\log \tau$ défini sur $\E^{\ent,\DP}$ et 
on démontre pour finir que celui-ci est en réalité défini sur un certain 
localisé explicite de $\Sk^\DP$ qui s'injecte dans $\O$.
\end{minipage} 
\end{center} 
\end{figure}

Soit $\frakM$ un objet de la catégorie $\Mod^{\varphi,\tau}_{/\Sk} 
\otimes \Q_p$, c'est-à-dire un $(\varphi,\tau)$-réseau défini sur $(\Sk, 
\Sk_\tau)$ qui est de $E(u)$-hauteur $\leq r$ pour un certain entier $r$. 
La définition du foncteur $\calR_{\varphi,\tau}$ conduit à poser 
simplement $\calS_{\varphi, N_\nabla}(\frakM) = \O \otimes_\Sk \frakM$ 
pour définir sa structure de $\O$-module, et à munir celui-ci du 
Frobenius $\varphi \otimes \varphi$. Ne reste donc plus qu'à définir 
l'opérateur $N_\nabla$. La discussion menée dans le \S 5.1 de 
\cite{liu-compo} montre que celui-ci doit s'obtenir \emph{via} la 
formule :
$$N_\nabla = \frac 1 {p\t} \cdot \log \tau = \frac 1 {p\t} \cdot 
\sum_{i=1}^\infty \frac{(\id-\tau)^i} i.$$
Toutefois, il n'est pas clair que cette formule ait un sens puisque
aussi bien la division par $\t$ que la convergeance de la somme posent
\emph{a priori} problème. Nous allons démontrer dans la suite que ce
n'est pas le cas, puis que l'opérateur $N_\nabla$ ainsi obtenu est 
défini sur $\O$. Lors de la démonstration, nous allons être amené à
introduire un certain nombre de modules sur $\Sk$ ; pour donner au
lecteur une vision d'ensemble de la situation, nous les avons regroupés
dans le diagramme de la figure \ref{fig:modules}.

\paragraph{L'espace $W(R)^{\DP}$}

On rappelle que $A_\cris$ est défini comme le complété $p$-adique de 
l'enveloppe à puissances divisées de $W(R)$ par rapport à l'idéal 
principal engendré par l'élément $E(u)$ vu dans $W(R)$. Suivant les 
notations habituelles, on pose $B_\cris^+ = A_\cris[1/p]$. Pour tout 
entier $m$, on définit $W(R)^{\DP}_{\leq m}$ comme le sous $W(R)$-module 
de $B_\cris^+$ engendré par $W(R)$ et les éléments de la forme 
$\frac{\t^{p^s}} {p^s} x$ pour $0 \leq s \leq m$ et $x \in W(\m_R)$.

Clairement, les $W(R)^{\DP}_{\leq m}$ forment une suite croissante pour 
l'inclusion ; on note $W(R)^{\DP}$ l'adhérence dans $B_\cris^+$ de 
$\bigcup_{m \geq 0} W(R)^{\DP}_{\leq m}$. Les éléments de cet espace 
sont les éléments de $B_\cris^+$ qui s'écrivent sous la forme $x_0 + 
\sum_{s=1}^\infty \frac{\t^{p^s}} {p^s} x_s$ avec $x_0 \in W(R)$ et $x_s 
\in W(\m_R)$ pour $s \geq 1$ (on notera qu'il n'y a aucune condition de 
convergence à imposer sur les $x_s$ pour que la somme converge).
De l'égalité $\varphi(\t) = U \t$, on déduit que les $W(R)^\DP_{\leq m}$ 
sont tous stables par $\varphi$. Ils sont également stables par l'action 
de $G_K$ étant donné que l'on a, d'une part, $g(\t) = \chi(g) \cdot t$ 
pour tout $g \in G_\infty$ et que, d'autre part, on sait que le quotient 
$\frac{\tau(\t)} {\t}$ est un élément de $W(R)$ (voir le deuxième 
exemple du \S \ref{subsec:exemples}).

\begin{prop}
\label{prop:logborne}
Il existe une constante $c$ (qui dépend de $\frakM$)
telle que, pour
tout entier $i \in \{1, \ldots, p^m\}$, on ait :
$$\frac{(\id - \tau)^i} i (\frakM) \subset p^{-c} \cdot W(R)^{\DP}
_{\leq m} \otimes_\Sk \frakM.$$
Pour cette même constante $c$, on a également :
$$\frac{(\id - \tau)^i} i (\frakM) \subset p^{-c} \cdot W(R)^{\DP}
\otimes_\Sk \frakM$$
pour tout entier $i \geq 1$.
\end{prop}

\begin{proof}
On vérifie que $\frac{(\id - \tau)^i} i = \sum_{j=1}^i (-1)^j \binom 
{i-1}{j-1} \frac{\tau^j - \id} j$.
Il suffit donc de montrer que, pour tout $j \in \{1, \ldots, p^m\}$, 
on a 
$$\left(\frac{\tau^j - \id} j\right) (\frakM) \subset p^{-c} \cdot 
W(R)^{\DP}_{\leq m} \otimes_\Sk \frakM$$
pour une certaine constante $c$.
De plus, étant donné que $W(R)^{\DP}_{\leq m}$ est stable par $\tau$, 
il suffit de vérifier l'inclusion précédente lorsque $j$ est une 
puissance de $p$, \emph{i.e.} $j = p^s$ avec $s \leq m$. Comme $\frakM$ 
est de $E(u)$-hauteur finie, il est de $E(u)$-hauteur $\leq p^h$ pour 
tout $h$ supérieur ou égal à un entier $h_0$. Le théorème 
\ref{theo:congr} affirme alors qu'il existe une constance $c' \geq 0$ ne 
dépendant que de $K$ telle que pour tout $h$ compris entre $h_0$ et 
$s-c'$, on ait l'inclusion :
$$(\tau^{p^s} - \id)(\frakM) \subset (\t^{p^h} W(\m_R) + p^{s-h-c'}
W(R)) \otimes_\Sk \frakM.$$
En prenant l'intersection sur tous les $h$, on obtient :
$$(\tau^{p^s} - \id)(\frakM) \subset p^{s-c'} \cdot \Big( p^{-h_0} W(R) 
+ \sum_{h=h_0}^{s-c'} \frac{\t^{p^h}}{p^h} W(\m_R) \Big) \otimes \frakM.$$
Or, puisque $s \leq m$, ce dernier espace est inclus dans $\frakM 
\subset p^{s-c} \cdot W(R)^{\DP}_{\leq m} \otimes_\Sk \frakM$ avec $c 
= c' + h_0$. La proposition en résulte.
\end{proof}

Le lemme ci-dessus montre que le terme général de la série 
$\sum_{i=1}^\infty \frac{(\id - \tau)^i} i (x)$ (pour $x \in \frakM$) 
définissant $(\log \tau)(x)$ vit dans un espace \og borné \fg\ (dans le 
sens où les dénominateurs en $p$ sont contrôlés). Ceci nous autorise à 
utiliser les résultats du \S \ref{subsec:logm} sur les logarithmes 
tronqués avec $a = \tau$. Plus précisément, soit $A$ l'algèbre des 
applications $\Q_p$-linéaires de $B_\cris^+ \otimes_\Sk \frakM$ dans 
lui-même et, pour tout entier $m$, soit $\Lambda_m$ le sous-ensemble de 
$A$ formé des fonctions qui envoient $\frakM$ sur $p^{-c} \cdot 
W(R)^{\DP}_{\leq m} \otimes_\Sk \frakM$. 
La proposition \ref{prop:logborne} dit alors exactement que $\tau \in A$ 
est $\Lambda_m$-borné à l'ordre $m$. Les propositions \ref{prop:logmmult}, 
\ref{prop:logmcont} et le corollaire \ref{cor:logmpuis} peuvent donc 
être appliqués dans ce contexte. On retiendra en particulier que pour 
tout $n \in \Z_p$, tout entier $m$ et tout $x \in \frakM$, on a 
la congruence :
\begin{equation}
\label{eq:logtaun}
(\log_m \tau^n)(x) \equiv n \cdot (\log_m \tau)(x) \pmod
{p^{m-c-1} W(R)^{\DP}_{\leq m} \otimes_\Sk \frakM}.
\end{equation}
qui nous sera bien utile dans la suite.

\smallskip

La suite de la démonstration consiste à borner, par étapes successives, 
l'image des applications $\log_m \tau$ définie sur $\frakM$ jusqu'à 
arriver à un espace suffisamment petit qui nous permettra en quelque 
sorte de recoller les $\log_m \tau$ pour former une authentique 
application $\log \tau$. On montrera alors que celle-ci prend ses 
valeurs dans $\t \O \otimes_\Sk \frakM$, ce qui nous permettra de 
définir l'opérateur $N_\nabla$ que nous avons évoqué précédemment et,
par suite, de conclure.

\paragraph{Un argument galoisien}

Le premier argument utilisé pour restreindre l'image de $\log_m \tau$ 
met à profit la relation de commutation \eqref{eq:commuttau2} qui 
apparaît dans la définition des $(\varphi,\tau)$-modules. Étant donné
que $\chi(\tau) = 1$, celle-ci s'écrit simplement $\tau^{\chi(g)} = (g
\otimes \id) \circ \tau \circ (g \otimes \id)^{-1}$ (pour $g \in
G_\infty/H_\infty$), l'égalité ayant lieu dans l'espace des
endomorphimes de $\Sk_\tau \otimes_\Sk \frakM$. En prenant le
logarithme et en utilisant \eqref{eq:logtaun}, on obtient, pour tout 
$x \in \frakM$, la congruence
\begin{equation}
\label{eq:congrgalois}
\chi(g) \cdot (\log_m \tau)(x) \equiv (g \otimes \id) \circ (\log_m
\tau)(x) \pmod{p^{m-c-1} W(R)^{\DP}_{\leq m} \otimes_\Sk \frakM}.
\end{equation}
(on rappelle que $G_\infty$ agit trivialement sur $\frakM$). On en 
déduit que $\log_m \tau(x) \in p^{-c} \cdot X_{m,m-1} \otimes_\Sk 
\frakM$ où
$$X_{m,n} = \big\{ \,\, \xi \in W(R)^{\DP}_{\leq m} \quad \big|
\quad 
\forall g \in G_\infty, \,\,\,
g(\xi) \equiv \chi(g) \xi \pmod{p^n W(R)^{\DP}_{\leq m}}
\,\, \big\}$$
et $c$ est la constante introduite dans la proposition \ref{prop:logborne}.
On est ainsi amené à étudier les $X_{m,n}$. Pour cela, on rappelle que 
l'on avait posé $L = \Frac R$ et que l'on note $v_p$ la valuation 
$p$-adique normalisée par $v_p(p) = 1$. On introduit l'espace 
$W(L)^{\DP}_{\leq m} = W(L) + W(R)^{\DP}_{\leq m}$ (où la somme 
est calculée dans $p^{-m} W(L)$) et on note $\E^{\ent,\perf,\DP}_{\leq 
m}$ l'ensemble des séries (formelles) de la forme
$$\sum_{q \in \Z[1/p]} a_q u^q, \quad a_q \in p^{-m} \cdot W$$
vérifiant les conditions de convergence suivantes :
\begin{itemize}
\item[a)] pour tout entier $n \geq 0$, on a $v_p(a_q) \geq -n$ dès que 
$q \leq e \cdot \frac{p^{n+1}-1}{p-1}$ ;
\item[b)] le coefficient $a_q$ tend vers $0$ lorsque $q$ tend vers
$-\infty$ ;
\item[c)] pour tout $\varepsilon > 0$, l'ensemble des $q \in \Z[1/p]$
tels que $|a_q|_p > \varepsilon$ est discret\footnote{Avec la condition
précédente, cela revient à dire qu'il intersecte tous les intervalles
$]-\infty,A]$ selon un ensemble fini.}.
\end{itemize}

\begin{lemme}
\label{lem:t}
L'élement $\t$ appartient à $W(\m_R)$.
\end{lemme}

\begin{proof}
Il s'agit de démontrer que l'image $\bar \t$ de $\t$ dans $W(\bar k)$ est
nulle. Or, en réduisant l'équation $\varphi(\t) = U\t$ dans $W(\bar k)$,
on trouve $\sigma(\bar \t) = p \bar \t$ où $\sigma$ désigne le Frobenius
naturel sur $W(\bar k)$. En prenant les valuations, on obtient $v_K(\bar
\t) = v_K(\bar \t) + 1$, ce qui n'est possible que si $\bar \t$ est nul.
\end{proof}

Le lemme \ref{lem:t} ci-dessus, le fait que $v_R(\t \mod p) = \frac 
e{p-1}$ ainsi que les conditions a), b) et c) qui apparaissent dans la 
définition de $\E^{\ent,\perf,\DP}$ assurent que l'association 
$\sum_q a_q u^q \mapsto \t \sum_q a_q [\upi^q]$ définit un morphisme 
$\E^{\ent,\perf,\DP}_{\leq m} \to W(L)^{\DP}_{\leq m}$ pour tout $m$. 
Celui-ci est injectif car la multiplication par $\t$ est déjà injective 
sur $p^{-m} \cdot W(L)$ (et que tout est plongé dans cet espace) ; dans 
la suite, son image sera notée $\E^{\ent,\perf,\DP}_{\leq m}(1)$ et 
l'image d'un élément $x \in \E^{\ent,\perf,\DP}$ dans $W(L)^{\DP}$ sera 
notée simplement $\t x$.

Pour simplifier les écritures, on pose dans la suite $s_m(n) = e \cdot 
\frac {p^n-1}{p-1}$ lorsque $n$ est élément de $\{1, \ldots, m\}$ et on 
convient que $s_m(0) = -\infty$ pour tout $m$ et $s_m(n) = +\infty$ pour 
$n > m$. Soit également $F_0^\perf$ l'adhérence dans $L = \Frac R$ du 
perfectisé de $F_0$. Pour tout $v \in \R \cup \{ -\infty\}$, on note 
$\a_{F_0^\perf} ^{> v}$ l'ensemble des éléments de $F_0^\perf$ de 
valuation strictement supérieure à $v$ et on convient que 
$\a_{F_0^\perf} ^{> \infty} = 0$.

\begin{lemme}
\label{lem:Eentperfmodp}
Pour tout entier $m$, on a un isomorphisme canonique :
$$\prod_{n = 0}^m \,
\frac{\a_{F_0^\perf}^{> s_m(n)}}{\a_{F_0^\perf}^{> s_m(n+1)}}
\stackrel{\sim}{\longrightarrow}
\frac{\E^{\ent,\perf,\DP}_{\leq m}}{p\:\E^{\ent,\perf,\DP}_{\leq m}}.$$
\end{lemme}

\begin{proof}
La clé consiste à remarquer que $F_0^\perf$ s'identifie à l'ensemble
des séries $\sum_{q \in \Z[1/p]} a_q u^q$ telles que $a_q \in k$ pour
tout $q$ et satisfaisant aux conditions b) et c) mentionnées
précedemment et que, sous cette identification, $\a_{F_0^\perf}^{>v}$
correspondant aux séries pour lesquelle $a_q = 0$ pour tout $q \leq v$.
Le morphisme qui apparaît dans l'énoncé du lemme est alors celui qui 
envoie une famille $(P_n)_{n \geq 0}$ sur la somme $\sum_{0 \leq n 
\geq m} p^{-n} P_n$. Le fait qu'il soit bien défini et qu'il réalise un 
isomorphisme entre les espaces considérées est enfin une vérification 
immédiate.
\end{proof}

De même, on démontre que, pour tout entier $m$, on a un isomorphisme
canonique :
\begin{equation}
\label{eq:WLdpmodp}
\prod_{n = 0}^m \,
\frac{\a_L^{> s(m)}}{\a_L^{> s(m+1)}}
\stackrel{\sim}{\longrightarrow}
\frac{W(L)^{\DP}_{\leq m}}{p \:W(L)^{\DP}_{\leq m}}
\end{equation}
où la flèche provient de la multiplication par $\t$.

\begin{prop}
\label{prop:fixes}
Pour tous entiers $n$ et $m$, on a :
$$\Bigg\{ \,\, x \in \frac{W(L)^\DP_{\leq m}}{p^n \: W(L)^\DP_{\leq m}} 
\quad \Big|
\quad \forall g \in G_\infty, \,\, g(x) = \chi(g) x \,\, \Bigg\} =
\frac{\E^{\ent,\perf,\DP}_{\leq m}(1)}{p^n \: 
\E^{\ent,\perf,\DP}_{\leq m}(1)}.$$
\end{prop}

\begin{proof}
Étant donné que $G_\infty$ agit trivialement sur les éléments $[\upi^q]$ 
et par la formule $g(\t) = \chi(g) \t$ sur l'élément $\t$, l'ensemble de 
droite est inclus dans celui de gauche. Comme ils sont tous deux 
complets pour la topologie $p$-adique (puisque annulés par une puissance 
de $p$), il suffit de démontrer qu'il y a égalité modulo $p$, 
c'est-à-dire lorsque $n = 1$. Dans ce cas, c'est une conséquence 
immédiate du lemme \ref{lem:Eentperfmodp} de la formule 
\eqref{eq:WLdpmodp} et de la proposition \ref{prop:LH}.
\end{proof}

Il résulte de la proposition que l'ensemble $X_{n,m}$ défini plus haut 
est inclus dans $\E^{\ent,\perf,\DP}_{\leq m}(1) + p^n W(L)^\DP_{\leq 
m}$ et donc que, pour tout entier $m$, l'application $\log_m \tau$ 
envoie $\frakM$ sur $p^{-c} \cdot (\E^{\ent,\perf,\DP}_{\leq m}(1) + 
p^{m-1} W(L)^\DP_{\leq m}) \otimes_\Sk \frakM$.

\paragraph{Une relation de Leibniz}

Le but de ce paragraphe est de démontrer que les opérateurs $\log_m
\tau$ vérifient des relations de Leibniz approchées qui nous seront
utiles par la suite. Pour tout entier $m$, on pose :
$$t_m = \log_m [\ueps] = \sum_{i=1}^{p^m-1} \frac {(1-[\ueps])^i} i.$$
Du fait que $v_R(1 - \ueps) = \frac{ep}{p-1} < v_R(\t \mod p) = \frac e 
{p-1}$, on déduit que $t_m \in W(R)^\DP_{\leq m}$ pour tout $m$.

\begin{prop}
\label{prop:Leibniz}
Pour tout entier $m$, tout $a \in \Sk$ et tout $x \in \frakM$, on a
la congruence :
$$(\log_m \tau)(ax) \equiv
t_m \cdot u \: \frac{da}{du} \cdot x + a \cdot (\log_m \tau)(x)
\pmod {p^{m-c-1} W(R)^\DP_{\leq m} \otimes_\Sk \frakM}$$
où la constante $c$ est celle de la proposition \ref{prop:logborne}.
\end{prop}

\begin{proof}
Par continuité et $W$-linéarité, il suffit de démontrer la proposition
lorsque $a$ est une puissance de $u$. Dans ce cas, on part de la
relation de commutation $\tau \circ u^n = u^n [\ueps]^n \circ \tau$.
En revenant à la définition du logarithme tronqué, on en déduit que
$(\log_m \tau) \circ u^n = u^n \circ \log_m ([\ueps]^n \circ \tau)$.
Comme $\tau$ commute à la multiplication par $[\ueps]$ (étant donné
que $\tau$ fixe cet élément), et que la multiplication $[\ueps]$ est 
$\Lambda_m$-borné à l'ordre $m$
(pour l'espace $\Lambda_m$ défini juste après la démonstration de la
proposition \ref{prop:logborne}), la proposition \ref{prop:logmmult}
s'applique et entraîne la congruence :
$$(\log_m \tau) \circ u^n \equiv n u^n t_m + u^n \log_m \tau 
\pmod{p^{m-1} u^n \Lambda_m}$$
Comme $\Lambda_m$ est stable par multiplication par $u^n$, la même 
congruence est \emph{a fortiori} vraie modulo $p^{m-1} \cdot \Lambda_m$ 
et la formule annoncée dans la proposition en résulte.
\end{proof}

\paragraph{L'élément $t_m$}

On rappelle que $W(R)$ est muni d'un morphisme d'anneaux $\theta$ à 
valeurs dans $\O_{\C_p}$ défini comme l'anneau des entiers du complété 
$p$-adique de $\bar K$. Il est bien connu que le noyau de $\theta$ est 
un idéal principal engendré par n'importe lequel de ses éléments dont la 
réduction modulo $p$ a pour valuation $e$. Des exemples de tels éléments 
sont $E(u)$ ou encore $\frac{[\ueps]-1} {[\ueps^{1/p}]-1} = 
\frac{\eta}{\varphi^{-1}(\eta)}$. On introduit l'idéal :
$$F^1 W(R) = \big\{ \, x \in W(R) \,\, \big| \,\, 
\theta(\varphi^i(x)) = 0, \, \forall i \, \big\}.$$

\begin{lemme}
\label{lem:F1WR}
On a $F^1 W(R) = \varphi(\t) \cdot W(R) = (1-[\ueps]) \cdot W(R)$.
\end{lemme}

\begin{proof}
On donne uniquement la démonstration de la première égalité, la
seconde étant absolument analogue.
Pour tout entier $i > 0$, l'élément $\varphi^i(\t)$ est égal à $E(u) 
\varphi(E(u)) \cdots \varphi^{i-1}(E(u)) \t$ et est donc multiple
de $E(u)$. Ainsi $\varphi^i(\t) \in \ker \theta$ pour tout $i > 0$,
et, par suite, $\varphi(\t) \in F^1 W(R)$. On a ainsi démontré 
l'inclusion $F^1 W(R) \supset \varphi(\t) \cdot W(R)$.

Pour démontrer l'inclusion réciproque, on introduit l'idéal $I$ de $R$ 
défini comme le quotient de $F^1 W(R)$ par $p F^1 W(R) = F^1 W(R) \cap p 
W(R)$. Manifestement, $I$ contient un élement de valuation $\frac{ep} 
{p-1}$, qui est la réduction modulo $p$ de $\varphi(\t)$. Par ailleurs, 
on remarque que tout élément $x \in F^1 W(R)$ est nécessairement 
multiple de $\alpha_i = \frac {\varphi^{-i}(\eta)} 
{\varphi^{-(i+1)}(\eta)}$ pour tout entier $i \geq 0$ (étant donné que 
$\varphi(x)$ est dans le noyau de $\theta$). Comme tous les idéaux 
principaux engendrés par les $\alpha_i$ sont premiers (car les quotient 
$W(R)/\alpha_iW(R)$ sont tous isomorphes à $\O_{\C_p}$ qui est 
intègre), on en déduit que $x$ est multiple de $\alpha_0 
\alpha_1 \cdots \alpha_i$ pour tout $i$. Or la réduction modulo $p$ de 
ce produit a pour valuation $e + \frac e p + \cdots + \frac e {p^i} = 
\frac{ep}{p-1} \cdot (1 - p^{-(i+1)})$. Ainsi, $v_R(x \mod p) \geq 
\frac{ep}{p-1} \cdot (1 - p^{-(i+1)})$ et, comme ceci est vrai pour tout 
$i$, on obtient $v_R(x \mod p) \geq \frac{ep}{p-1}$. On en déduit que 
$I$ est l'idéal de $R$ des éléments de valuation $\geq \frac{ep}{p-1}$.
Il résulte de cela que le morphisme d'inclusion $\iota : \varphi(\t) 
\cdot W(R) \to F^1 W(R)$ induit un isomorphisme modulo $p$. Comme les 
espaces de depart et d'arrivée sont complets pour la topologie 
$p$-adique, il suit que $\iota$ est lui-même un isomorphisme, ce qui 
signifie que $F^1 W(R) = \varphi(\t) \cdot W(R)$.
\end{proof}

Il résulte du lemme et du fait que $1-[\ueps]$ peut s'écrire sous
la forme $\t^p \alpha + p \beta$ avec $\alpha$ et $\beta$ dans $W(R)$
(ce qui suit de l'égalité $v_R(1-\ueps) = p v_R(\t \mod p)$), que la
fraction 
$$t'_m = \frac{t_m}{\t} = \sum_{i=1}^\infty E(u) \cdot \frac{1-[\ueps]}
{\varphi(\t)} \cdot \frac{(1-[\ueps])^{i-1}} i$$
appartient à $W(R)^\DP_{\leq m}$. D'autre part, étant 
donné que $[\ueps]$ est $(W(R)^\DP_{\leq m})$-borné à l'ordre $m$ (dans 
la $\Q_p$-algèbre $B_\cris^+$ par exemple), le corollaire 
\ref{cor:logmpuis} implique :
$$\begin{array}{l}
g(t_m) = \log_m ([\ueps]^{\chi(g)}) \equiv \chi(g) \cdot t_m
\pmod {p^{m-1} \cdot W(R)^\DP_{\leq m}} \medskip \\
\varphi(t_m) = \log_m ([\ueps]^p) \equiv p \cdot t_m
\pmod {p^{m-1} \cdot W(R)^\DP_{\leq m}}
\end{array}$$
la première congruence étant vraie pour tout $g \in G_\infty$. Ainsi 
$t'_m$ est fixe par $G_\infty$ modulo $p^{m-2} \cdot W(R)^\DP_{\leq m}$ 
(on perd une puissance de $p$ à cause de la division par $\t$). Soit 
$\Sk^{\perf,\DP}_{\leq m}$ (resp. $\Sk^\DP_{\leq m}$) le sous-ensemble 
de $\E^{\ent,\perf,\DP}_ {\leq m}$ formé des séries $\sum a_q u^q$ pour 
lesquelles $a_q = 0$ dès que $q < 0$ (resp. dès que $q < 0$ ou non 
entier).
En reprenant l'argument de la démonstration de la proposition 
\ref{prop:fixes}, on démontre que $t'_m$ s'écrit comme la somme d'un 
élément $x \in \Sk^{\perf,\DP}_{\leq m}$ et d'un élément de $y \in 
p^{m-2} \cdot W(R)^\DP_{\leq m}$. De plus, la congruence portant sur 
$\varphi(t_m)$ implique que $x \equiv \frac{E(u)}{E(0)} \cdot \varphi(x) 
\pmod {p^{m-3} \cdot \Sk^{\perf,\DP}_{\leq m}}$. En remarquant que tous 
les itérés de l'opérateur $\frac{E(u)}{E(0)} \cdot \varphi$ envoient 
$p^{m-2} \cdot W(R)^\DP_{\leq m}$ sur $p^{m-3} \cdot W(R)^\DP$, on 
démontre en itérant la congruence précédente que l'on a nécessairement 
$x \in \Sk^\DP_{\leq m} + p^{m-3} \cdot W(R)^\DP$. Comme $x$ est aussi 
dans $\Sk^{\perf,\DP}_{\leq m} + p^{m-2} \cdot W(R)^\DP_{\leq m}$, on 
trouve en prenant l'intersection de ces deux espaces que $x \in 
\Sk^\DP_{\leq m} + p^{m-3} \cdot W(R)^\DP_{\leq m}$. Au final, on a donc 
démontré que :
\begin{equation}
\label{eq:tm}
t_m \in \t \cdot \Sk^\DP_{\leq m} + p^{m-3} \cdot W(R)^\DP_{\leq m}.
\end{equation}

\paragraph{Élimination des puissances fractionnaires}

Nous revenons à notre problématique consistant à restreindre l'image de 
$\log_m \tau$. Pour tout entier $m$, on définit 
$\E^{\ent,\varphi^{-n},\DP}_{\leq m}$ le sous-espace de 
$\E^{\ent,\perf,\DP}_{\leq m}$ formé des séries $\sum_q a_q u^q$ pour 
lesquelles $a_q = 0$ dès que $p^n q$ n'est pas un entier. Dans la suite, 
on notera simplement $\E^{\ent,\DP}_{\leq m}$ au lieu de 
$\E^{\ent,\varphi^{-0},\DP}_{\leq m}$. En outre, si $E$ est l'un de
ces espaces, on notera $E(1)$ son image dans $W(R)^\DP_{\leq m}$ par
la multiplication par $\t$. On se propose de démontrer, dans 
ce paragraphe, qu'il existe une constance $c' \geq c$ telle que, pour 
tout entier $m$ suffisamment grand, on ait :
\begin{equation}
\label{eq:sansfrac}
(\log_m \tau)(\frakM) \subset p^{-c'} \cdot 
\Bigg( \E^{\ent,\DP}_{\leq m}(1) \,+\, p^m \cdot W(L)^\DP_{\leq m} 
\,+\, \frac{\t^{p^m}}{p^m} W(R) \Bigg) 
\otimes_\Sk \frakM
\end{equation}
Pour éliminer ainsi les puissances fractionnaires de $u$, l'idée
principale consiste à réduire petit à petit les dénominateurs qui
apparaissent sur les puissances de $u$ en utilisant le Frobenius
$\varphi$ et, plus précisément, les deux faits suivant le concernant : 
\emph{primo}, il commute à $\log_m \tau$ et \emph{secundo} son image
est assez grande (ce qui se traduit concrètement par la finitude de
la $E(u)$-hauteur de $\frakM$). Dans toute la suite, on considère un
entier $h$ tel que $\frakM$ soit de $E(u)$-hauteur $\leq p^h$. On a
alors le lemme suivant qui précise l'idée générale que l'on vient
d'énoncer.

\begin{lemme}
\label{lem:etape}
On suppose qu'il existe un $\Sk$-module $D$ tel que 
$(\log_m \tau)(\frakM) \subset D \otimes_\Sk \frakM$. Alors :
$$E(u)^{p^h} \cdot (\log_m \tau)(\frakM) \subset 
(\Sk \varphi(D) + \Sk t_m + p^{m-c-1} \cdot W(R)^\DP_{\leq m})
\otimes_\Sk \frakM.$$
\end{lemme}

\begin{proof}
On considère un élément $x \in \frakM$. Comme $\frakM$ est supposé
de $E(u)$-hauteur $\leq p^h$, on peut écrire $E(u)^{p^h} x$ sous la
forme :
$$E(u)^{p^h} x = \varphi(x_0) + u  \varphi(x_1) + \cdots + 
u^{p-1} \varphi(x_{p-1})$$
où les $x_i$ sont des éléments de $\frakM$. En appliquant $\log_m \tau$ 
à cette égalité, en utilisant la relation de Leibniz approchée démontrée 
précédemment (voir proposition \ref{prop:Leibniz}) et le fait évident
que $\log_m \tau$ commute à $\varphi$, on obtient la congruence :
$$\begin{array}{l}
\displaystyle
E(u)^{p^h} (\log_m \tau)(x) \equiv - t_m \cdot u \:\frac{d E(u)^{p^h}}{du} 
\cdot x + \sum_{i=0}^{p-1} \big(i t_m u^i \varphi(x_i) + u^i \varphi \circ
(\log_m \tau)(x_i)\big) \\
\hspace{8cm} \pmod{p^{m-c-1} \cdot W(R)^\DP_{\leq m} \otimes_\Sk \frakM}
\end{array}$$
Dans le membre de droite de l'expression ci-dessus, le premier terme
et le premier terme de la somme appartiennent à $t_m \frakM$, tandis 
que le deuxième terme de la somme appartient à $\Sk \varphi(D) 
\otimes_\Sk \frakM$. Le lemme en résulte.
\end{proof}

Il s'agit maintenant d'appliquer le lemme \ref{lem:etape} avec le bon 
$D$. À cette fin, on introduit les ensembles suivants paramétrés par les 
deux entiers $m$ et $n$ :
\begin{eqnarray*}
A_m & = & p^{-h} \cdot W(L) \,+\, W(L)^\DP_{\leq m} \\
B_{m,n} & = & p^{-h} \cdot \E^{\ent,\varphi^{-n},\DP}_{\leq 0}(1) \,+\,
\E^{\ent,\varphi^{-n},\DP}_{\leq m}(1) \\
C_m & = & \frac{\t^{p^m}}{p^m} \cdot W(R).
\end{eqnarray*}
On rassemble ci-dessous quelques propriétés de ces espaces.

\begin{lemme}
\label{lem:Am}
Pour tous entiers $m$, $n$ avec $m \geq h$ et $n \geq 1$, on a :
$$A_m \subset \frac{E(u)^{p^h}} p \cdot A_m
\quad ; \quad
B_{m,n} \subset \frac{E(u)^{p^h}} p \cdot B_{m,n}$$
$$\varphi(A_m) \subset E(u)^{p^h} A_m 
\quad ; \quad
\varphi(B_{m,n}) \subset E(u)^{p^h} B_{m,n-1}
\quad ; \quad
\varphi(C_m) \subset E(u)^{p^h} C_m.$$
\end{lemme}

\begin{proof}
La dernière inclusion est évidente étant donné que $\varphi(\t^{p^m}) = 
\t^{p^m} \cdot U^{p^m}$ (où on rappelle que $U = \frac{E(u)}{E(0)/p}$ et
que $\frac{E(0)} p$ est un élément inversible dans $W(R)$ puisqu'il 
l'est déjà dans $\Z_p$).
La méthode pour démontrer les quatre autres inclusions consiste à 
remarquer que $E(u)^{p^h}$ est inversible dans $W(L)$ et à exprimer son 
inverse de façon appropriée. Précisément, on commence par observer que 
$v_R(E(u)^{p^h} \mod p) = v_R(\t^{(p-1)p^h} \mod p)$, ce qui montre que 
l'on peut écrire $E(u)^{p^h}$ s'écrit comme un produit $\alpha \cdot 
(\t^{(p-1)p^h} - p\beta)$ où $\alpha$ est un élément inversible dans 
$W(R)$ et où $\beta \in W(R)$. On en déduit que :
\begin{equation}
\label{eq:divEuph}
\frac 1 {E(u)^{p^h}} = \alpha^{-1} \cdot 
\sum_{i=0}^\infty \frac{\beta^i}{\t^{(i+1)(p-1)p^h}} \cdot p^i.
\end{equation}
À partir de là, il est facile de démontrer la première inclusion. En 
effet, il suffit de vérifier que si $x$ est dans $p^{-h} \cdot W(L)$ ou 
s'il s'écrit $x = \frac{t^{p^s}}{p^s} \cdot y$ avec $y \in W(\m_R)$, 
alors $\frac p {E(u)^{p^h}} \cdot x \in A_m$, et ceci s'obtient en 
remplaçant la division par $E(u)^{p^h}$ par l'expression donnée par 
la formule \eqref{eq:divEuph} et en constantant que chaque terme de la somme 
obtenue est, lui-même, dans $A_m$. La troisième inclusion (\emph{i.e.}
la première inclusion de la deuxième ligne) se démontre pareilllement.
Pour les inclusions faisant intervenir les $B_{m,n}$, la méthode est
similaire sauf qu'au lieu d'utiliser la formule \eqref{eq:divEuph},
on utilise la formule
$$\frac 1 {E(u)^{p^h}} = 
\sum_{i=0}^\infty \frac{F(u)^i}{u^{(i+1)ep^h}} \cdot p^i$$
qui s'obtient en écrivant $E(u)^{p^h}$ sous la forme $u^{ep^h} + p F(u)$ 
avec $F(u) \in \Sk$.
\end{proof}

On pose à présent 
$D_{m,n} = p^{-c-1} \cdot \big( p^{m-1} A_m + B_{m,n} + C_m 
\big)$.
Sachant que $(\log_m \tau)(\frakM) \subset p^{-c} \cdot (\E^{\ent, 
\perf,\DP}_{\leq m}(1) + p^{m-1} W(L)^\DP_{\leq m}) \otimes_\Sk \frakM$, 
on s'aperçoit, en revenant aux définitions, que pour tout couple 
$(m,s)$, il existe un entier $n$ tel que $(\log_m \tau)(\frakM) \subset 
D_{m,n} \otimes_\Sk \frakM$. Par ailleurs, il résulte du lemme 
\ref{lem:Am} et de l'appartenance $t_m \in B_{m,0} + p^{m-3} A_m$ (qui 
découle de \eqref{eq:tm}) que, pour tout $n \geq 1$, on a :
$$\Sk \varphi(D_{m,n}) + \Sk t_m + p^{m-c-1} \cdot W(R)^\DP_{\leq m}
\subset E(u)^{p^h} \cdot D_{m,n-1}$$
pourvu que l'on ait pris soin de choisir la constante $c$ supérieure
ou égale à $3$ (ce qui est bien sûr toujours possible). Une application
itérée du lemme \ref{lem:etape} montre alors que, si $m \geq h$, on a
$(\log_m \tau)(\frakM) \subset D_{m,s,0} \otimes_\Sk \frakM$, d'où il
découle l'inclusion \eqref{eq:sansfrac}.

\paragraph{Passage à la limite sur $m$}

La clé pour effectuer le passage à la limite est le lemme suivant.

\begin{lemme}
Pour tout entier $m$, l'image $Q_m$ de $W(R)^\DP_{\leq m} \cap 
\E^{\ent,\DP}_{\leq m} (1)$ dans le quotient
$$\frac{W(L)^\DP_{\leq m}}{p^m \cdot W(L)^\DP_{\leq m} + 
\frac{\t^{p^m}}{p^m} \cdot W(R)}$$
est finie.
\end{lemme}

\begin{proof}
Exercice. (On rappelle, à toutes fins utiles, que l'on a supposé que le 
corps résiduel $k$ de $W$ est fini.)
\end{proof}

Il suit du lemme que la limite projective des $Q_m$ (muni de la 
topologie de la limite projective) est un ensemble compact. On en 
déduit, à l'aide de l'inclusion \eqref{eq:sansfrac} que l'on a démontrée 
précédemment, que, pour tout $x \in \frakM$, il existe une suite 
strictement croissante d'entiers $(m_k)_{k \geq 1}$ telle que la suite 
extraite des $(\log_{m_k} \tau)(x)$ converge, pour la topologie de la 
limite projective, vers un élement de $\E^{\ent,\DP}(1) \otimes_\Sk 
\frakM$. Comme $\frac{\t^{p^m}}{p^m}$ converge vers $0$ dans $B_\cris^+$ 
quand $m$ tend vers l'infini, cette même suite converge également pour
la topologie usuelle de $B_\cris^+$ (\emph{i.e.} la topologie 
$p$-adique). 

Soit $e_1, \ldots, e_d$ une base de $\frakM$ sur $\Sk$. Quitte à 
extraire à nouveau, on peut supposer que, pour tout indice $i \in \{1, 
\ldots, d\}$, la suite des $(\log_{m_k} \tau)(e_i)$ est convergente. On 
note $p \t \cdot x_i$ sa limite, de sorte que $x_i$ soit un élément de 
$p^{-c'} \cdot \E^{\ent,\DP} \otimes_\Sk \frakM$. Ceci nous permet de 
définir un opérateur $N_\nabla : \frakM \to p^{-c'} \cdot \E^{\ent,\DP} 
\otimes_\Sk \frakM$ en posant :
$$N_\nabla\Bigg(\sum_{i=1}^d a_i e_i\Bigg) = 
\sum_{i=1}^d N_\nabla(a_i) \otimes e_i + a_i x_i
\qquad (a_i \in \Sk).$$
Par construction, $N_\nabla$ vérifie la relation de Leibniz :
\begin{equation}
\label{eq:leibniz}
N_\nabla(ax) = N_\nabla(a) x + a N_\nabla(x)
\end{equation}
pour tout $a \in \Sk$ et tout $x \in \frakM$. À partir de là, il suit de 
la proposition \ref{prop:Leibniz} que $(\log_{m_k} \tau)(x)$ converge 
vers $p\t \cdot N_\nabla(x)$ pour tout $x \in \frakM$.
Enfin, le fait que $\varphi$ commute à tous les $\log_m \tau$ implique
que l'on a la relation :
\begin{equation}
\label{eq:commphiN}
N_\nabla \circ \varphi(x) = U \cdot (\varphi \otimes \varphi) \circ
N_\nabla(x)
\end{equation}
pour tout $x \in \frakM$.

\paragraph{Élimination des puissances négatives}

Soit $\Sk^{\DP}$ le sous-espace de $\E^{\ent,\DP}$ formé des 
séries $\sum_q a_q u^q$ (avec $q \in \Z$) pour lesquelles $a_q = 0$ dès 
que $q$ est strictement négatif. On note également $N$ le plus grand 
entier (éventuellement nul) pour lequel $E(u) \in \varphi^n(\Sk)$. On
pose
$$D_n(u) = \varphi^{-1}(E(u)) \cdot \varphi^{-2}(E(u)) \cdots 
\varphi^{-n}(E(u))$$
pour tout entier $n \leq N$.

\begin{lemme}
\label{lem:intersection}
Soit $x \in \E^{\ent,\DP}$ vérifiant $\t x \in W(R)^{\DP}$. 
Alors, il existe $y \in \Sk^{\DP}$ tel que $D_N(u)\cdot (x + y) \in 
\Sk$.
\end{lemme}

\begin{proof}
En écrivant $x$ comme la somme d'un élément de $\E^\ent$ et d'un élément 
de $\Sk^{\DP}$, on voit qu'il suffit de démontrer que tout élément $x 
\in \E^\ent$ tel que $\t x \in W(R)$ vérifie :
$$\varphi^{-1}(E(u)) \cdot \varphi^{-2}(E(u)) \cdots \varphi^{-n}(E(u))
\cdot x \in \Sk.$$
On note $\Sk'$ l'ensemble des éléments $x$ de $\E^\ent$ tels que $\t x 
\in W(R)$. L'ensemble des $v_R(x \mod p)$ pour $x$ parcourant $\Sk'$ est 
clairement minoré par $-\frac e {p-1}$ ; on peut donc considérer un 
élément $a \in \Sk'$ tel que $v_R(a \mod p)$ soit minimal. Comme $1 \in 
\Sk'$, on a de surcroît $v_R(a \mod p) \leq 0$.
L'élément $a$ est tel que $a \Sk \subset \Sk'$, et il est facile de 
vérifier que cette inclusion induit un isomorphisme modulo $p$. On en 
déduit que $\Sk' = a\Sk$. Par ailleurs, comme $1 \in \Sk'$, il existe un 
élément $b \in \Sk$ tel que $ab = 1$. D'après une variante du théorème 
de préparation de Weierstrass (voir par exemple théorème 2.1, chap. 5 de 
\cite{lang}), $b$ s'écrit comme le produit d'un élément inversible de 
$\Sk$ et d'un polynôme $B \in W[u]$ de la forme $B(u) = u^d + p(b_{d-1} 
u^{d-1} + \cdots + b_0)$. On a l'égalité
$$\t \cdot \frac{E(u)}{\varphi(B(u))} = c \cdot \varphi \Big( 
\frac \t {B(u)} \Big)$$
qui montre que la fraction $\frac{E(u)}{\varphi(B(u))}$ appartient à
$\Sk'$, et s'écrit donc sous la forme $\frac{C(u)}{B(u)}$ pour une
certaine série $C(u) \in \Sk$. On en déduit que $\varphi(B(u))$ divise
le produit $B(u) E(u)$ dans $\Sk$. Il résulte directement de cette 
divisibilité que $B(u)$ n'est pas multiple de $u$.

Si $B(u)$ n'a aucune racine non nulle dans $\O_{\bar K}$, alors c'est
un polynôme constant qui est inversible dans $\Sk$. Dans ce cas, on
obtient donc $\Sk = \Sk'$ et le lemme est démontré.
Supposons maintenant qu'au contraire $B(u)$ admette une racine non nulle 
dans $\O_{\bar K}$. Le polynôme $B^\sigma$ obtenu en appliquant $\sigma$ 
à $B$ admet alors lui aussi au moins une racine non nulle dans $\O_{\bar 
K}$. Soit $x$ une telle racine de valuation minimale et $y$ une racine 
$p$-ième de $x$. De $\varphi(B(u)) = B^\sigma(u^p)$, on déduit que $y$ 
est une racine du polynôme $\varphi(B(u))$. Par ailleurs, comme $x$ a 
été choisi de valuation minimale, $y$ n'annule pas $B(u)$. Comme il 
annule, par contre, le produit $B(u) E(u)$, on a nécessairement $E(y) = 
0$. Comme $E(u)$ est irréductible dans $\Sk$ (car c'est un polynôme
d'Eisenstein), on a nécessairement :
\begin{equation}
\label{eq:Eudivise}
E(u) \text{ divise } \varphi(B(u)).
\end{equation}
Nous allons à présent montrer que $E(u)$ dans l'image de $\varphi$.
Pour cela, on repart de l'égalité $E(y) = 0$ qui montre que $y$ et $\pi$ 
sont conjugués sous Galois dans $\O_{\bar K}$. Il en est donc de même
de $x = y^p $ et $\pi^p$. Comme $x$ est annulé par $B^\sigma$ qui est
de degré $\leq \frac e{p-1}$, on en déduit que le corps $K' = W[1/p]
(\pi^p)$ est de degré au plus $\frac e{p-1}$ sur $W[1/p]$. Par suite,
l'extension $K/K'$ est de degré au moins $p-1$. Comme elle est, par
ailleurs, de degré au plus $p$ (puisque elle est engendré par $\pi$ qui 
est manifestement annulé par un polynôme de degré $p$ à coefficients 
dans $K'$), on a $[K:K'] \in \{p-1, p\}$. Si on avait $[K:K'] = p-1$,
le polynôme minimal de $\pi$ sur $K'$ serait de degré $p-1$ et serait
également un diviseur de $X^p - \pi^p \in K'[X]$ ; ce dernier polynôme
admettrait donc également un facteur de degré $1$, ce qui n'est pas
possible car il ne peut y avoir dans $K'$ un élément de même valuation
que $\pi$. Ainsi $[K:K'] = p$ et $[K':W[1/p]] = \frac e p$ (et donc,
en particulier, $p$ divise $e$). Le polynôme minimal $Q$ de $\pi^p$ sur 
$W[1/p]$ est donc de degré $\frac e p$ et $\pi$ est annulé par $Q(u^p)$
qui est de degré $e$. On en déduit que $Q(u^p) = E(u)$ et, par suite,
que $E(u)$ est dans l'image de $\varphi$ comme annoncé.

La divisibilité \eqref{eq:Eudivise} implique alors que $\varphi^{-1} 
(E(u))$ divise $B(u)$ dans $\Sk$, \emph{i.e.} $B(u) = \varphi^{-1} 
(E(u)) \cdot B_1(u)$ pour un certain élément $B_1(u) \in \Sk$. De plus, 
le fait que $\varphi(B(u))$ divise $B(u) E(u)$ nous apprend que 
$\varphi(B_1(u))$ divise $B_1(u) \varphi^{-1}(E(u))$. On peut ainsi 
appliquer à nouveau le même raisonnement que précédemment qui conduit
à l'éventualité suivante : soit $B_1(u)$ est inversible, soit $E(u)$
est dans $\varphi^2(\Sk)$ et $B_1(u)$ divise $\varphi^{-2}(E(u))$.
Dans le premier cas, le lemme est démontré tandis que, dans le second
cas, on écrit $B_1(u) = \varphi^{-2}(E(u)) \cdot B_2(u)$ et on applique
à nouveau le même argument jusqu'à ce que $E(u)$ ne soit plus élément
de $\varphi^n(\Sk)$ (ce qui arrive nécessairement).
\end{proof}

Rappelons que nous avions démontré précédemment que $p N_\nabla(\frakM) 
\subset p^{-c'} \cdot \E^{\ent,\DP} \otimes_\Sk \frakM$. Comme, par 
ailleurs, $(\log_m \tau)(\frakM) \subset p^{-c'} \cdot W(R)^{\DP} 
\otimes_\Sk \frakM$ pour tout $m$, l'application $p \t \cdot N_\nabla$ 
prend également ses valeurs dans $p^{-c'} \cdot W(R)^{\DP} \otimes_\Sk 
\frakM$ et le lemme ci-dessus implique que :
\begin{equation}
\label{eq:NnablaO}
N_\nabla(\frakM) \subset p^{-(c'+1)} \cdot \Big( \Sk^{\DP} +
\frac 1 {D_N(u)} \Sk \Big) \otimes_\Sk \frakM \subset
\frac 1 {D_N(u)} \cdot \calM
\end{equation}
où on rappelle que $\calM = \O \otimes_\Sk \frakM$.

Soit, comme précédemment, $h$ un entier pour lequel $\frakM$ est de 
$E(u)$-hauteur $\leq p^h$. Soit $x \in E(u)^{p^h} \frakM$. Il s'écrit 
sous la forme $x = \sum_{i=0} ^{p-1} u^i \varphi(x_i)$ pour des $x_i 
\in \frakM$. La relation \eqref{eq:commphiN} implique alors que
\begin{equation}
\label{eq:NnablaO2}
N_\nabla(x) \in \frac 1 {\varphi(D_N(u)) \cdot E(u)^{p^h}}\cdot \O 
\otimes_\Sk \frakM = 
\frac 1 {D_{N-1}(u) E(u)^{p^h+1}} \cdot \calM
\end{equation}
Or, $D_{N-1}(u)$ n'est pas nul dans le quotient $\O/E(u)\O \simeq K$.
On en déduit que l'intersection des idéaux principaux de $\O$ engendrés
respectivement par $D_{N-1}(u)$ et par $E(u)^{p^h+1}$ est l'idéal principal
engendré par le produit $D_{N-1}(u) E(u)^{p^h+1}$. Ainsi, en mettant ensemble
les inclusions \eqref{eq:NnablaO} et \eqref{eq:NnablaO2}, on trouve
$N_\nabla(\frakM) \subset \frac 1 {D_{N-1}(u)} \cdot \calM$. En 
réappliquant le même argument $N-1$ fois, on obtient $N_\nabla (\frakM) 
\subset \calM$. Ainsi $\calM$ muni des opérateurs $\varphi$ et 
$N_\nabla$ est un objet de $\Mod^{\varphi,N_\nabla,0}_{/\O}$ et on peut 
poser $\calS_{\varphi,N_\nabla}(\frakM) = \calM$ : le foncteur 
$\calS_{\varphi,N_\nabla}$ est construit !

\subsubsection{Démonstration des propriétés annoncées de
$\calS_{\varphi,N_\nabla}$}
\label{subsec:propcalS}

On considère $\frakM$ un objet de la catégorie 
$\Mod^{\varphi,\tau}_{/\Sk} \otimes \Q_p$ et on pose $\calM = 
\calS_{\varphi,N_\nabla}(\frakM)$. Soient encore $s$ le plus grand 
entier tel que l'extension $K_{s}/K$ soit galoisienne et $t$ le plus 
petit entier tel que $p^t(p-1)$ dépasse strictement le rang de $\frakM$ 
comme $\Sk$-module. Notre premier objectif est de démontrer que les 
représentations galoisiennes $V$ et $V_\st$ (\og $\st$ \fg\ pour 
semi-stable) associées respectivement à $\frakM$ et $\calM$ coïncident 
en restriction au sous-groupe $G_n$ pour $n = s$ et $n=t$.
Pour cela, on pose $\frakM_\st = \calR_{\varphi,\tau}(\calM)$ : c'est un 
$(\varphi,\tau)$-réseau de $E(u)$-hauteur finie qui vit à l'intérieur du 
$(\varphi,\tau)$-module correspondant à $V_\st$. De plus, par construction, 
les $\varphi$-modules $\frakM[1/p]$ et $\frakM_\st[1/p]$ sont isomorphes. 
Ainsi, plutôt que de voir $\frakM[1/p]$ et $\frakM_\st[1/p]$ comme deux 
objets différents, on pourra et on préfèrera considérer $\frakM_\st 
[1/p]$ comme le $\varphi$-module $\frakM[1/p]$ muni d'un nouvel 
opérateur $\tau_\st$. De la même façon, on verra dans la suite $V_\st$ 
comme le $\Q_p$-espace vectoriel $V$ muni d'une action différence de 
$\tau \in G_K$. Avec ces notations, on cherche à démontrer que 
$\tau^{p^n} = \tau_\st^{p^n}$ (soit sur $\frakM[1/p]$, soit sur $V$, 
c'est équivalent) pour $n = s$ et $n = t$. Or, on possède deux 
informations essentielles sur les opérateurs $\tau$ et $\tau_\st$, à 
savoir :

\begin{enumerate}
\item par construction, il existe une suite strictement croissante 
d'entiers $(m_k)_ {k \geq 0}$ telle que, pour tout $x \in \frakM[1/p]$, 
la suite des $\log_{m_k}(\tau)(x)$ converge vers $p \t \cdot 
N_\nabla(x)$ (quand $k$ tend vers l'infini) ;
\item d'après les résultats de Liu, pour tout $x \in \frakM[1/p]$,
la suite des $\log_m (\tau_\st)(x)$ converge vers $p \t \cdot
N_\nabla(x)$ (quand $m$ tend vers l'infini).
\end{enumerate}
Pour avancer, il semble donc qu'il faille relier d'une façon ou d'une 
autre les opérateurs $\log_m \tau$ et $\log_m \tau_\st$ agissant sur 
$\frakM[1/p]$ à l'action de l'élément $\tau \in G_K$ sur les 
représentations galoisiennes $V$ et $V_\st$. C'est l'objet du la
proposition suivante.

\begin{prop}
\label{prop:logtauN}
Soit $T$ une $\Z_p$-représentation galoisienne libre de $E(u)$-hauteur 
finie et soit $\frakN$ l'unique $(\varphi,\tau)$-réseau de $E(u)$-hauteur 
finie vivant à l'intérieur de son $(\varphi,\tau)$-module. On note 
$\tau_\frakN$ (resp. $\tau_{\Sk^\ur}$, resp. $\tau_T$) l'opérateur 
$\tau$ agissant sur $\frakN$ (resp. $\Sk^\ur$, resp. $T$). Alors, il 
existe une constance $C$ telle que, pour tout entier $m$, on ait :
$$(\log_m \tau_\frakN)(f) \equiv (\log_m \tau_{\Sk^\ur}) \circ f - 
f \circ (\log_m \tau_T) \pmod {p^{m-C} \cdot W(R)^{\DP}}$$
pour tout $f \in \frakN = \Hom_{\Z_p[G_\infty]}(T, \Sk^\ur)$ (où la
congruence signifie que la différence des deux fonctions prend ses
valeurs dans $p^{m-C} \cdot W(R)^{\DP}$).
\end{prop}

\begin{proof}
On rappelle que par définition $\tau_\frakN(f) = \tau_{\Sk^\ur}
\circ f \circ \tau_T^{-1}$. Ainsi :
$$\frac{(\id - \tau_\frakN)^i} i(f) = \frac 1 i \cdot 
\sum_{\alpha = 0}^i (-1)^\alpha \binom i \alpha \cdot \tau_{\Sk^\ur}^\alpha 
\circ f \circ \tau_T^{-\alpha}.$$
En écrivant $\tau_T^{-1} = \id-g$ et en injectant cela dans l'écriture
précédente, on obtient :
$$\frac{(\id - \tau_\frakN)^i} i(f) = \frac 1 i \cdot 
\sum_{0 \leq \beta \leq \alpha \leq i} (-1)^{\alpha+\beta} \binom i \alpha
\cdot \binom \alpha \beta \cdot \tau_{\Sk^\ur}^\alpha \circ f \circ 
g^\beta.$$
En isolant les termes pour lesquels $\beta = 0$ et en remarquant que,
si $\beta \neq 0$, on a $\frac 1 i \binom i \alpha \binom \alpha \beta
= \frac 1 \beta \binom {i-1}{\beta-1} \binom{i-\beta}{\alpha - \beta}$,
on aboutit à la nouvelle expression :
\begin{eqnarray}
\frac{(\id - \tau_\frakN)^i} i(f) 
& = & \frac{(\id - \tau_{\Sk^\ur})^i} i \circ f + 
\sum_{0 \leq \beta \leq \alpha \leq i} (-1)^{\alpha+\beta} \binom {i-1}
{\beta-1} \cdot \binom {i-\beta}{\alpha-\beta} \cdot \tau_{\Sk^\ur}^\alpha 
\circ f \circ \frac{g^\beta}{\beta} \nonumber \\
& = & \frac{(\id - \tau_{\Sk^\ur})^i} i \circ f + 
\sum_{\beta=0}^i \binom {i-1} {\beta-1} \cdot \tau_{\Sk^\ur}^\beta \circ 
\Bigg( \sum_{\alpha=\beta}^i
\binom {i-\beta}{\alpha-\beta} (-\tau_{\Sk^\ur})^{\alpha-\beta} \Bigg)
\circ f \circ \frac{g^\beta}{\beta} \nonumber \\
& = & \frac{(\id - \tau_{\Sk^\ur})^i} i \circ f + 
\sum_{\beta=0}^i \binom {i-1} {\beta-1} \cdot \tau_{\Sk^\ur}^\beta \circ 
(\id - \tau_{\Sk^\ur})^{i-\beta} \circ f \circ \frac{g^\beta}{\beta}
\label{eq:tauN}
\end{eqnarray}
En sommant maintenant sur $i$ variant entre $1$ et $p^m-1$, on obtient :
\begin{equation}
\label{eq:logtauN}
(\log_m \tau_\frakN)(f) = (\log_m \tau_{\Sk^\ur}) \circ f + 
\sum_{\beta=1}^{p^m-1} \tau_{\Sk^\ur}^\beta \circ
\Bigg(\sum_{i=\beta}^{p^m-1} \binom {i-1} {\beta-1} 
(\id - \tau_{\Sk^\ur})^{i-\beta} \Bigg) \circ
f \circ \frac{g^\beta}{\beta}.
\end{equation}
D'autre part, d'après la proposition \ref{prop:logborne}, on sait que 
$\frac{(\id - \tau_\frakN)^i} i(f) \in \Lambda \otimes_\Sk \frakM$ avec 
$\Lambda = p^{-c}\cdot W(R)^{\DP}$ pour une certaine constante $c$. 
Étant donné que tous les éléments de $\frakM$ sont des fonctions à valeurs 
dans $\Sk^\ur \subset W(R)$, on déduit de la formule \eqref{eq:tauN} 
ci-dessus, que la composée $\frac{(\id - \tau_{\Sk^\ur})^i} i \circ f$ 
prend ses valeurs dans $\Lambda$ pour tout entier $i$. De même que l'on 
a obtenu la formule \eqref{eq:Lambdaborne}, il s'ensuit que $\frac{(\id 
- \tau_{\Sk^\ur})^i} i \circ f$ prend ses valeurs dans $p^{\ell(i)} 
\Lambda$ pour tout $i$. On déduit de ce fait que :
\begin{equation}
\label{eq:taubeta}
\sum_{i=\beta}^{p^m-1} \binom {i-1} {\beta-1}
(\id - \tau_{\Sk^\ur})^{i-\beta} \circ f \equiv \tau^{-\beta} \circ f
\pmod {p^{\ell(p^m-j)} \Lambda}
\end{equation}
De plus, il existe une constante $c'$ telle que $g^\beta$ prenne ses 
valeurs dans $p^{[c'\beta]} T$ pour tout $\beta$. Par ailleurs, pour $j 
\in \{1, \ldots, p^m-1\}$, la différence $\ell(p^m-j) - [c'j] - v_p(j)$ 
est minorée par $\ell(p^m-j) - [c'j] - \log_p(j)$, quantité elle-même 
minorée par $m-C$ pour une certaine constante $C$ (qui dépend de $c'$). 
À partir de là, on déduit la congruence désirée en injectant 
\eqref{eq:taubeta} dans \eqref{eq:logtauN}.
\end{proof}

Soit $T$ un réseau de $V$ stable par les deux actions de $G_K$, à savoir 
celle provenant de $V$ et celle provenant de $V_\st$. Si l'on note 
respectivement $\tau$ et $\tau_\st$ les endomorphismes de $T$ donnés par 
l'action de $\tau \in G_K$ sur $V$ et $V_\st$, la proposition 
\ref{prop:logtauN} appliquée avec $T$ (vu successivement comme un réseau 
de $V$ et de $V_\st$) permet de déduire des informations 1 et 2 degagées 
plus haut que
\begin{equation}
\label{eq:verszero}
f \circ \big(\log_{m_k} \tau - \log_{m_k} \tau_\st\big)
\text{ tend vers } 0.
\end{equation}
Or, comme on a choisi $\tau$ dans le sous groupe d'inertie sauvage, les 
suites $\tau^{p^n}$ et $\tau_\st^{p^n}$ tendent vers $1$, et les 
opérateurs $\log \tau$ et $\log \tau_\st$ sont bien définis. En passant à 
la limite dans \eqref{eq:verszero}, on obtient $f \circ \log \tau = f 
\circ \log \tau_\st$. Comme ceci est vrai pour tout $f \in 
\Hom_{\Z_p[G_\infty]}(T, \Sk^\ur)$, on obtient finalement $\log \tau = 
\log \tau_\st$ comme endomorphismes $\Z_p$-linéaires de $T$.
Montrons à présent que cela implique que $\tau^{p^t} = \tau_\st^{p^t}$ 
pour tout entier $t$ tel que $p^t(p-1) > \dim_{\Q_p} V$. Le lemme
clé est le suivant.

\begin{lemme}
\label{lem:rdc}
Soit $M$ un $\Z_p$-module libre de rang $d$ et $f : M \to M$ un
endomorphisme $\Z_p$-linéaire. On suppose qu'il existe un entier
$n$ tel que $f^{p^n} \equiv \id \pmod p$. Alors pour tout entier 
$t \geq 0$ tel que $d \geq p^{t-1}(p-1)$, on a :
$$v_p\big((\id-f^{p^t})^i\big) \geq \frac{p^t} d \cdot i - 1$$
où on convient que $v_p(g) \geq v$ signifie que $g$ est divisible
par $p^v$.
\end{lemme}

\begin{proof}
De $f^{p^n} \equiv \id \pmod p$, on déduit que $f-\id$ induit un 
endomorphisme nilpotent sur $M/pM$ et donc que $(f-\id)^d \equiv 0 \pmod 
p$. On pose $g = f-\id$. On a alors $v_p(g^d) \geq 1$ et, plus 
généralement pour tout entier $j \geq 0$, $v_p(g^j) \geq \frac d j 
- 1$. D'autre part, un calcul immédiat montre que $(\id - f^{p^t})^i$ 
s'exprime en fonction de $g$ comme suit :
$$(\id - f^{p^t})^i = - \sum_{1 \leq j_1, \ldots, j_i \leq p^t}
\binom {p^t}{j_1} \cdot \binom {p^t}{j_2} \cdots \binom {p^t}{j_i}
\cdot g^{j_1 + \cdots + j_i}$$
où la notation signifie que la somme est prise sur tous les $i$-uplets 
$(j_1, \ldots, j_i)$ d'éléments de $\{1, \ldots, p^t\}$. On en déduit 
que :
$$v_p((\id-f^{p^t})^i) \geq \min_{1 \leq j_1, \ldots, j_i \leq p^t}
\frac{j_1 + \cdots + j_i} d + (t - v_p(j_1)) + \cdots + (t -
v_p(j_i)) - 1.$$
On observe facilement que le minimum est atteint lorsque chaque $j_i$ 
est une puissance de $p$ : $j_i = p^{a_i}$ pour un certain $a_i \in
\{0, \ldots, t\}$. Pour $j_i$ de cette forme, la quantité à minimiser 
vaut $-1 + \sum_{\alpha=1}^i \frac{p^{a_i}} d + (t - a_i)$ et 
l'hypothèse $d \geq p^{t-1}(p-1)$ assure que chaque terme de la somme 
précédente est minimal lorsque $a_i = t$. Le lemme en découle.
\end{proof}

On est maintenant en mesure de conclure. On note $d$ la dimension de 
$V$. On rappelle en outre que $t$ désigne le plus petit entier tel que
$p^t(p-1) > d$. On a donc $p^{t-1}(p-1) \leq d$ ; autrement dit,
$t$ satisfait l'inégalité du lemme \ref{lem:rdc}. Appliqué avec
$M = T$ et $f = \tau$ puis avec $f = \tau_\st$, il donne :
$$v_p\big((\id-\tau^{p^t})^i\big) \geq \frac{p^t} d \cdot i - 1
\quad \text{et} \quad
v_p\big((\id-\tau_\st^{p^t})^i\big) \geq \frac{p^t} d \cdot i - 1.$$
Étant donné que $\frac{p^t} d > \frac 1{p-1}$ et que le rayon de 
convergence de l'exponentielle est $|p|^{1/(p-1)}$, on déduit des
estimations ci-dessus que les séries $\exp(\log \tau^{p^t})$ et $\exp(
\log \tau_\st^{p^t})$ convergent respectivement vers $\tau^{p^t}$ et 
$\tau_\st^{p^t}$. Comme les logarithmes sont égaux (puisque $\log 
\tau^{p^t} = p^t \log \tau = p^t \log \tau_\st = \log \tau_\st^{p^t}$), 
il s'ensuit $\tau^{p^t} = \tau_\st^{p^t}$ comme voulu.

\medskip

Soit maintenant $s$ le plus grand entier pour lequel l'extension $K_s/K$ 
est galoisienne. Le sous-groupe $G_s$ est alors distingué dans $G_K$ et 
c'est le plus petit contenant $G_\infty$ ayant cette propriété. En 
d'autres termes, les conjugués de $G_\infty$ par les éléments de $G_K$ 
engendrent $G_s$. On souhaite montrer que $V$ est $V_\st$ sont 
isomorphes en tant que représentations de $G_s$. Soit $g \in G_K$. On 
peut appliquer ce que l'on vient de démontrer en remplaçant le système 
comptatible de racines $p^{s'}$-ièmes de l'uniformisante $(\pi_{s'})$ 
par $(g \pi_{s'})$. Ce faisant, on obtient l'existence d'une représentation 
semi-stable $V_{\st,g}$ qui coïncide avec $V$ en restriction au 
sous-groupe $g G_t g^{-1}$. Si maintenant $g_1$ et $g_2$ sont deux 
éléments de $G_K$, les représentations $V_{\st,g_1}$ et $V_{\st,g_2}$ 
coïncident sur l'intersection $(g_1 G_t g_1^{-1}) \cap (g_2 G_t 
g_2^{-1})$ qui est un sous-groupe d'indice fini de $G_K$ qui découpe 
une extension totalement ramifiée. D'après la proposition 
\ref{prop:sstextension}, les représentations $V_{\st,g_1}$ et 
$V_{\st,g_2}$ coïncident en fait sur tout $G_K$. Il en résulte que
$V$ coïncide avec $V_\st$ sur tous les conjugués de $G_t$ et donc,
par suite, sur le sous-groupe engendré par ces conjugués. Or celui-ci
n'est autre que $G_s$ ; on a donc bien démontré ce que l'on voulait.

\medskip

Pour établir complètement le théorème \ref{theo:calS}, il ne reste plus 
qu'à démontrer que si $\calM$ est un objet de 
$\Mod^{\varphi,N_\nabla,0}_{/\O}$ alors $\calS_{\varphi,N_\nabla} \circ 
\calR_{\varphi,\tau}(\calM)$ est isomorphe à $\calM$. Mais, si l'on note 
$V$ (resp. $W$) la représentation galoisienne associée à $\calM$ (resp. 
à $\calS_{\varphi,N_\nabla} \circ \calR_{\varphi,\tau}(\calM)$), on 
sait, par ce qui vient d'être fait, que $V$ et $W$ coïncident sur $G_s$. 
Or, $V$ et $W$ sont deux représentations semi-stables. Une nouvelle 
application de la proposition \ref{prop:sstextension} permet donc de 
conclure.

\subsection{Réseaux à l'intérieur des représentations semi-stables}
\label{subsec:classres}

Un problème classique et important consiste à décrire les réseaux à
l'intérieur des représentations semi-stables. Dans ce dernier
paragraphe, nous aimerions expliquer sommairement comment cela peut être
réalisé à l'aide de structures entières à l'intérieur des objets de
$\Mod^{\varphi, N_\nabla,0}_{/\O}$.

\subsubsection{Une borne plus précise pour l'action de $N_\nabla$}

Le point de départ de notre analyse consiste à remarquer que la 
construction que nous avons faite du foncteur $\calS_{\varphi, 
N_\nabla}$ (voir \S \ref{subsec:constcalS}) donne en réalité un peu plus 
que ce qui a été énoncé jusqu'à présent. Précisément, elle assure que, 
si $\frakM$ est un $(\varphi,\tau)$-réseau de $E(u)$-hauteur finie, on 
n'a pas seulement $N_\nabla(\frakM) \subset \O \otimes_\Sk \frakM$, mais 
aussi
\begin{equation}
\label{eq:inclR1}
N_\nabla(\frakM) \subset \calR^\ent_1 \otimes_\Sk \frakM
\end{equation}
où $\calR^\ent_1$ désigne le sous-$\Sk$-module de $\O$ formé des séries $f
= \sum_{n \geq 0} a_n u^n$ pour lesquelles la quantité $v_p(a_n) +
\log_p(n)$ est minorée pour $n \geq 1$. Ce n'est pas la première fois
qu'un tel espace est introduit ; il a été, par exemple, déjà considéré dans
\cite{colmez}, \S II.1, référence dans laquelle l'auteur montre en
particulier qu'il s'agit d'un espace de Banach pour la valuation
$v_{\calR^\ent_1}$ définie, en reprenant les notations précédentes, par
$v_{\calR^\ent_1}(f) = \inf_{n \geq 1} v_p(a_n) + \ell(n)$ où
$\ell(0) = 0$ et $\ell(n) = \log_p(n)$ si $n \geq 1$.
La proposition II.3.1 de \emph{loc. cit.} montre en outre que
$\calR^\ent_1$ s'identifie à l'espace des distributions continues sur
$\Z_p$ d'ordre $1$, c'est-à-dire définissant une forme linéaire continue
sur l'espace des fonctions de classe $\mathcal C^1$.
Dans ce numéro, nous allons démontrer que l'on peut améliorer encore 
l'inclusion \eqref{eq:inclR1} en remplaçant $\calR^\ent_1$ par une 
partie bornée complètement explicite de cet espace.

On rappelle, pour commencer, que l'on dispose d'un morphisme surjectif 
$\theta : W(R) \to \O_{\C_p}$ dont le noyau est l'idéal principal 
engendré par $E(u)$. On rappelle également que l'anneau de Fontaine 
$A_\cris$ est alors défini comme le complété $p$-adique de l'enveloppe à 
puissances divisées de $W(R)$ relativement à $\ker \theta$ et que l'on a 
posé $B_\cris^+ = A_\cris[1/p]$. La série $\sum_{i=1}^\infty 
\frac{(1-[\ueps])^i} i$ définit un élément $t \in A_\cris$.

\begin{lemme}
\label{lem:tiAcris}
Pour tout entier $i$, on a $(t^i B_\cris^+) \cap W(R) = \varphi(\t)^i 
\cdot W(R)$.
\end{lemme}

\begin{proof}
Étant donné que $\varphi(\lambda \t) = -t$ et que $\varphi(\lambda)$ est 
inversible dans $B_\cris^+$, l'inclusion $\varphi(\t)^i \cdot W(R) \subset 
(t^i B_\cris^+) \cap W(R)$ est immédiate. Pour l'inclusion réciproque, on 
raisonne par récurrence sur $i$. Pour $i=0$, il n'y a rien à dire. 
Supposons que l'on ait $(t^i B_\cris^+) \cap W(R) = \varphi(\t)^i \cdot
W(R)$ et 
considérons un $x \in (t^{i+1} B_\cris^+) \cap W(R)$. Bien sûr, $x$ est 
alors aussi élément de $(t^i B_\cris^+) \cap W(R)$ et donc, d'après 
l'hypothèse de récurrence, de $\varphi(\t)^i \cdot W(R)$. On peut donc 
écrire $x = \varphi(\t)^i y$ avec $y \in W(R)$. D'autre part, $x$ 
s'écrit également par hypothèse $x = t^{i+1} z$ avec $z \in B_\cris^+$. 
On obtient $t^{i+1} z = \varphi(\t)^i y$, ce qui donne $y = (-1)^{i+1} 
\varphi(\lambda^{i+1} \t) z$. En particulier, $y$ est multiple de $E(u)$ 
(puisque $\varphi(\t)$ l'est) et même, plus généralement, $E(u)$ divise 
$\varphi^j(y)$ pour tout entier $j$. Autrement dit $\theta(\varphi^j(y)) 
= 0$ pour tout $j$ et donc, par le lemme \ref{lem:F1WR}, $y \in 
\varphi(\t) \cdot W(R)$. Finalement, on obtient $x \in \varphi(\t)^{i+1} 
\cdot W(R)$. 
\end{proof}

On considère, à présent, $V$ une représentation semi-stable de $G_K$ et 
$T \subset V$ un $\Z_p$-réseau stable par l'action de Galois. Soit 
$\frakM$ l'unique $(\varphi,\tau)$-réseau de $E(u)$-hauteur finie à 
l'intérieur du $(\varphi,\tau)$-module associé à $T$. D'après le 
résultat principal de \cite{liu2}, l'opérateur $\tau$ induit un 
endomorphisme de $\hat \calR \otimes_ {\varphi,\Sk} \frakM$ avec
$$\hat \calR = W(R) \cap \Bigg\{ \,
\sum_{n=0}^\infty f_n \cdot \frac{t^n}{p^{q(n)} q(n)!}
\, \text{ avec } f_n \in S[1/p], \, f_i \to 0 \,\Bigg\}
\subset A_\cris$$
où $q(n)$ désigne le quotient de la division euclidienne de $n$ par 
$p-1$. On en déduit que pour tout entier $i$, l'opérateur $(\id-\tau)^i$ 
envoie $\Sk \otimes_{\varphi,\Sk} \frakM$ sur $(t^i B_\cris^+) \otimes_ 
{\varphi,\Sk} \frakM$. Comme, en outre, $(\id-\tau)^i$ stabilise $W(R) 
\otimes_{\varphi,\Sk} \frakM$, on déduit du lemme \ref{lem:tiAcris} que 
$(\id-\tau)^i$ envoie $\Sk \otimes_{\varphi,\Sk} \frakM$ sur 
$\varphi(\t)^i \cdot W(R) \otimes_{\varphi,\Sk} \frakM$ et finalement 
que :
$$(\id-\tau)^i(\frakM) \subset \t^i \cdot W(R) \otimes_\Sk \frakM.$$
On en déduit la proposition suivante :

\begin{prop}
\label{prop:Nnabla}
Soient $V$ une représentation semi-stable de $G_K$, et $T$ un
$\Z_p$-réseau de $V$ stable par $G_K$. On note $\frakM$ un\footnote{Il
est en fait unique d'après la proposition \ref{prop:propB}.}
$(\varphi,\tau)$-réseau de $E(u)$-hauteur finie à l'intérieur du
$(\varphi,\tau)$-module associé à $T$. Soit $\calM$ l'objet de
$\Mod^{\varphi,N_\nabla,0}_{/\O}$ associé à $V$. Alors
$$N_\nabla(\frakM) \subset \Sk_\nabla \otimes_\Sk \frakM$$
où $\Sk_\nabla$ est l'ensemble des séries de la forme 
$\sum_{n \geq 0} \frac{P_n(u)}{p^{n+1}} u^{e (p^n-1)/(p-1)}$ 
où les $P_n(u)$ sont des polynômes à coefficients dans $W$.
\end{prop}

\begin{rem}
Comme le module $\Sk_\nabla$ est inclus dans $\calR^\ent_1$, la 
proposition apparaît comme un raffinement de l'inclusion 
\eqref{eq:inclR1}. Mieux encore, $\Sk_\nabla$ apparaît comme un 
sous-ensemble borné dans $\calR^\ent_1$ ; de façon explicite, il est 
inclus dans la boule de centre $0$ et de rayon $\frac{p^2} e$.
\end{rem}

\subsubsection{Reconstruction de $\tau$ à partir de $N_\nabla$}
\label{subsec:reconsttau}

Les bornes que l'on vient d'obtenir permettent de donner un sens à 
l'expression $\exp(p\t N_\nabla)$ et ainsi de reconstruire $\tau$ à 
partir de la connaissance de $N_\nabla$. Il faut toutefois être prudent
car les itérés successifs de $p \t N_\nabla$ ne sont pas clairement 
définis puisque $N_\nabla$ n'est défini \emph{a priori} que sur $\calM = 
\O \otimes_\Sk \frakM$ qui n'est manifestement pas stable par 
multiplication par $\t$.

Pour contourner ce problème, on procède comme suit. Pour tout entier 
positif ou nul $i$, on définit $\calR^\ent_i$ l'ensemble des séries $f = 
\sum_{n \geq 0} a_n u^n$ pour lesquelles la quantité $v_p(a_n) + i 
\log_p(n)$ est minorée pour $n \geq 1$. Si $i$ et $j$ sont deux entiers, 
le produit d'une fonction de $\calR^\ent_i$ par une fonction de 
$\calR^\ent_j$ appartient à $\calR^\ent_{i+j}$. La relation de Leibniz 
permet de prolonger l'opérateur $N_\nabla$ à $\O \otimes_\Sk \frakM$, 
tandis que la remarque que l'on vient de faire montre que ce 
prolongement, que l'on appelle encore $N_\nabla$, envoie $\calR^\ent_i 
\otimes_\Sk \frakM$ sur $\calR^\ent_{i+1} \otimes_\Sk \frakM$. Les 
formules $N_\nabla^{(0)} = \id$ et
$$N_\nabla^{(i+1)} = i u \cdot \frac{d \lambda}{du} \cdot N_\nabla^{(i)}
+ N_\nabla \circ N_\nabla^{(i)}$$
définissent alors des applications $\Sk$-linéaires $N_\nabla^{(i)} :
\frakM \to \calR^\ent_i \otimes_\Sk \frakM$ qui, d'un point de vue
intuitif, doivent être pensées comme les quotients $\frac{(p\t
N_\nabla)^i}{(p\t)^i}$. Sachant que $N_\nabla(\frakM) \subset \Sk_\nabla 
\otimes_\Sk \frakM$, il n'est pas difficile de démontrer que la somme
infinie $\sum_{i \geq 0} \frac{(p\t)^i}{i!} \cdot N_\nabla^{(i)}$
converge vers un opérateur $W$-linéaire $\tau : \frakM \to
B_\cris^+ \otimes_\Sk \frakM$. Une récurrence sur $i$ montre que,
pour tout $x \in \frakM$ et tout entier $i$, on a la formule :
$$N_\nabla^{(i)}(ux) = u \cdot \sum_{j=0}^i \binom i j (-\lambda)^{i-j}
N_\nabla^{(j)}(x).$$
Il en résulte par un nouveau calcul que $\tau(ux) = [\ueps] \cdot
\tau(x)$ pour tout $x \in \frakM$, relation qui permet de prolonger
$\tau$ par semi-linéairité à tout $B_\cris^+ \otimes_\Sk
\frakM$. Par ailleurs, on montre, à nouveau par récurrence sur $i$, que 
$N_\nabla^{(i)}$ et $\varphi$ satisfont à la relation de commutation 
$N_\nabla^{(i)} \circ \varphi = U^i \cdot \varphi \circ N_\nabla^{(i)}$, 
à partir de quoi il suit que $\tau$ et $\varphi$ commutent. Finalement, 
on laisse en exercice au lecteur le soin de vérifier que la série 
définissant le logarithme de $\tau$ définit un endormorphisme $\frakM 
\to B_\cris^+ \otimes_\Sk \frakM$ qui coïncide avec $p\t \cdot N_\nabla$.

\begin{rem}
La reconstruction de $\tau$ que l'on vient de présenter implique de 
nouvelles contraintes sur l'espace d'arrivée de l'opérateur $\tau$. 
Soit, en effet, pour tout entier $i$, $\Sk_\nabla^{(i)}$ le 
sous-$\Sk$-module de $\O$ engendré par les produits de $i$ éléments de 
$\Sk_\nabla$. Par convention, on pose également $\Sk_\nabla^{(0)} = 
\Sk$. Pour tout $i$, $\Sk_\nabla^{(i)}$ est alors un sous-ensemble borné 
de $\calR^\ent_i$ et les $p^i \Sk_\nabla^{(i)}$ se plongent également 
dans $W(R)[\t/p]^\wedge$ en envoyant, comme d'habitude, $u$ sur 
$[\upi]$. Il est alors immédiat de vérifier que les morphismes 
$N_\nabla^{(i)}$ définis précédemment prennent leurs valeurs dans 
$\Sk_\nabla^{(i)} \otimes_\Sk \frakM$, et donc que :
\begin{equation}
\label{eq:tausst}
\tau(\frakM) \subset \Bigg( \sum_{i \geq 0} p^i \Sk_\nabla^{(i)}
\cdot \t^i \Bigg)^{\!\wedge} \otimes_\Sk \frakM
\end{equation}
où le produit $p^i \Sk_\nabla^{(i)}$ est vu dans $W(R)[\t/p]^\wedge$
et l'exposant ${}^\wedge$ signifie que l'on prend l'adhérence dans
$W(R)[\t/p]^\wedge$.
\end{rem}

\subsubsection{Structures entières à l'intérieur des objets de
$\Mod^{\varphi, N_\nabla,0}_{/\O}$}

Le résultat de la proposition \ref{prop:Nnabla} conduit naturellement à
la définition (sans doute provisoire) suivante. 

\begin{deftn}
Un \emph{$(\varphi,N_\nabla)$-réseau} est la donnée d'un objet $\frakM$ de
$\Mod^\varphi_{/\Sk}$ muni d'un morphisme $N_\nabla : \frakM \to
\Sk_\nabla \otimes_\Sk \frakM$ vérifiant la relation de Leibniz
\eqref{eq:leibniz} et la relation de commutation \eqref{eq:commphiN}.
\end{deftn}

\noindent
Un morphisme entre deux $(\varphi,N_\nabla)$-réseaux $\frakM$ et $\frakM'$
est un morphisme $f$ entre les objets de $\Mod^\varphi_{/\Sk}$ faisant
commuter le diagramme suivant :
$$\xymatrix @C=50pt {
\frakM \ar[r]^-{N_\nabla} \ar[d]_-{f} & 
\Sk_\nabla \otimes_\Sk \frakM \ar[d]^-{\id \otimes f} \\
\frakM' \ar[r]^-{N_\nabla} & \Sk_\nabla \otimes_\Sk \frakM' 
}$$
On note $\Mod^{\varphi, N_\nabla}_{/\Sk}$ la catégorie des $(\varphi,
N_\nabla)$-réseaux. Si, par ailleurs, $\Rep^{\ent,\st}_{[0,+\infty[}
(G_K)$ désigne la catégorie des réseaux stables par $G_K$ à l'intérieur
des représentations semi-stables à poids de Hodge-Tate positifs ou nuls,
on a construit précédemment un foncteur $\calK^{\ent} :
\Rep^{\ent,\st}_{[0,+\infty[} (G_K) \to \Mod^{\varphi, N_\nabla}_{/\Sk}$
qui s'insère dans le carré commutatif suivant :
$$\xymatrix @C=50pt {
\Rep^{\ent,\st}_{[0,+\infty[} (G_K) \ar[r] \ar[d]_-{\calK^{\ent}} &
\Rep^\st_{[0,+\infty[} (G_K) \ar[d]_-{\sim}^-{\calK} \\
\Mod^{\varphi, N_\nabla}_{/\Sk} \ar[r] & 
\Mod^{\varphi, N_\nabla,0}_{/\O}
}$$

\begin{prop}
Le foncteur $\calK^\ent$ est pleinement fidèle. De plus, un
$(\varphi,N_\nabla)$-réseau $\frakM$ est dans son image essentielle si,
et seulement si l'endomorphisme $\tau$ de $B_\cris^+ \otimes_\Sk
\frakM$ défini au \S \ref{subsec:reconsttau} stabilise $W(R) \otimes_\Sk
\frakM$.
\end{prop}

\begin{proof}
C'est immédiat.
\end{proof}

\subsection{Bornes pour la ramification sauvage des représentations
semi-stables}
\label{subsec:bornesst}

Nous concluons cette section en expliquant sommairement comment la 
combinaison de plusieurs idées qui ont été développées dans les pages 
précédentes permettent de compléter les méthodes de \cite{carliu} et
de démontrer la conjecture 1.2.(1) dont voici l'énoncé.

\begin{theo}
\label{theo:bornesst}
Soit $r$ un entier positif. Soit $T$ le quotient de deux réseaux stables 
par $G_K$ dans une représentation semi-stable à poids de Hodge-Tate dans 
$\{0, \ldots, r\}$. On se donne un entier $n$ tel que $p^n T = 0$. Si 
$\frac r{p-1}$ s'écrit sous la forme $p^{\alpha} \beta$ avec $\alpha' 
\in \mathbb N$ et $\frac 1 p < \beta' \leq 1$, alors pour tout 
$$\mu > 1 + e (n + \alpha) + \max(e\beta - \frac 1 {p^{n+\alpha}},\, 
\frac e{p-1})$$ 
le sous-groupe de ramfication $G_K^{(\mu)}$ agit trivialement sur $T$.
\end{theo}

\begin{proof}
Elle est analogue à celle présentée dans \cite{carliu} sauf que, comme
dans le \S \ref{subsec:bornes} de cet article, on remplace les quotients
$W_n(R)/[\mathfrak a_R^{>a}]$ et $W_n(R)/[\mathfrak a_R^{>b}]$ (avec les
notations de \emph{loc. cit.}) par $W_n(R)/(E(u)^r \t^r) W_n(\m_R)$ et
$W_n(R)/\t^r W_n(\m_R)$ respectivement.
La clé réside en fait dans une généralisation appropriée du lemme 3.2.1
de \cite{carliu} qui stipule que, si $\frakM$ est l'unique
$(\varphi,\tau)$-réseau de $E(u)$-hauteur $\leq r$ dans le
$(\varphi,\tau)$-module associé à un réseau dans une représentation
semi-stable à poids de Hodge-Tate dans $\{0, \ldots, r\}$, alors pour
tout $s > n-1 + \log_p r$ et pour tout $\sigma \in G_s$, on a :
$$(\sigma - \id) (\frakM) \subset (\t^r \Sk_\tau^+ + p^n \Sk_\tau) 
\otimes_\Sk \frakM.$$
Il suffit bien sûr d'établir cette inclusion lorsque $g = \tau^{p^s}$.
Dans ce cas, elle découle d'une expression de $\tau^{p^s}$ en fonction
de l'opération $N_\nabla$ analogue à celle qui a été établie pour
$\tau$ dans la démonstration de la proposition \ref{prop:Nnabla}. Si
l'on préfère, elle peut également s'obtenir, de même que dans 
\cite{carliu}, comme une conséquence de la théorie des $(\varphi,
\hat G)$-modules de Liu. Le reste de la démonstration est absolument 
similaire à \cite{carliu} ; on ne le répète donc pas ici.
\end{proof}

\section{Quelques perspectives}

\subsubsection*{Le cas $p = 2$}

Tout au long de cet article, nous avons supposé que le nombre premier
$p$ était impair. Bien que cette hypothèse ait été utilisée à plusieurs
reprises au fil des démonstrations, l'auteur est d'avis que l'ensemble
des résultats obtenus devrait s'étendre (sans doute avec quelques
modifications mineures) au cas $p=2$. L'exercice reste cependant à
faire.

\subsubsection*{Lien avec les $(\varphi,\Gamma)$-modules}

Comme cela a été démontré dans cet article, la catégorie des
$(\varphi,\tau)$-modules est équivalente à celle des représentations
galoisiennes de $G_K$. Ainsi elle est aussi équivalente à la catégorie
des $(\varphi,\Gamma)$-modules de Fontaine. Expliciter cette équivalence
sans passer par les représentations de $G_K$ nous semble une question
naturelle et intéressante.
L'obtention d'un tel résultat pourrait permettre de déduire du théorème
\ref{theo:calS} un nouveau critère pour reconnaître les représentations
semi-stables en termes de leurs $(\varphi,\Gamma)$-modules. Se poserait
alors la question de comparer celui-ci avec celui qui découle de la
théorie de Berger (voir \cite{berger}). Le rapprochement des deux points
de vue pourrait peut-être conduire à une meilleure compréhension des
représentations semi-stables.

Dans le même veine, on peut chercher à rendre explicite les liens entre
les différentes catégories de $(\varphi,\tau)$-modules que l'on obtient
en faisant varier l'uniformisante $\pi$, la famille des $\pi_n$, ou
encore l'élément $\tau$. Ceci devrait permettre une meilleure
compréhension des $(\varphi, \tau)$-modules. Un moyen, qui semble
raisonnable, pour aborder cette question consiste à adapter les
constructions du \S 3.1 de \cite{caruso-crelle} en gardant à l'esprit
que $\tau$ joue le rôle de l'opérateur $\exp(tN)$ où $t = \log [\ueps]
\in A_\cris$.

\subsubsection*{Surconvergence des $(\varphi,\tau)$-modules}

Un résultat important de la théorie des $(\varphi,\Gamma)$-modules est
le théorème de Cherbonnier-Colmez (voir \cite{chercolmez}) qui affirme
que tout $(\varphi,\Gamma)$-module étale sur $\E$ admet un \og
$(\varphi,\Gamma)$-réseau \fg\ défini sur un anneau des séries
surconvergentes, c'est-à-dire convergentes sur une coronne
infinitésimale sur le bord du disque unité.
Dans cet article, nous nous sommes contentés de considérer des réseaux
définis sur l'anneau $\Sk$ dont les éléments convergent dans tout le 
disque unité. Comme nous l'avons vu, cela suffit pour l'application aux 
représentations semi-stables à poids de Hodge-Tate positifs ou nuls car, 
d'après les résultats de Kisin, les $(\varphi,\tau)$-modules qui leur 
sont associés admettent toujours de tels réseaux.
Par contre, comme nous l'avons montré dans le \S \ref{subsec:critere},
ce n'est pas le cas de tous les $(\varphi,\tau)$-modules et, notamment
de celui correspondant à la représentation $\Z_p(-1)$ que l'on a \emph{a 
priori} pas envie d'écarter. Pour pouvoir continuer à travailler avec 
cet exemple basique (et bien d'autres), il paraît donc important de 
comprendre si --- et, le cas échéant, comment --- le théorème de 
Cherbonnier-Colmez s'étend aux $(\varphi,\tau)$-modules.

\subsubsection*{Utilisation du logarithme dans le cas de torsion}

Dans le \S \ref{sec:Eu}, nous avons vu que, dans le cas des $(\varphi,
\tau)$-modules libres sur $\E^\ent$, la considération du logarithme de
$\tau$ s'est relévé être un outil puissant pour décrire cet opérateur.
C'est par exemple elle qui nous a permis d'obtenir l'inclusion
\eqref{eq:tausst}.
On peut donc se demander dans quelle mesure des arguments similaires
s'appliquent dans le cas des $(\varphi, \tau)$-modules de $p$-torsion.
Bien entendu, la considération du logarithme devient alors bien plus
délicate à cause des divisions par $p$ qu'il faudra désormais contrôler
avec plus d'attention.

\subsubsection*{Le problème du rélèvement des représentations}

Dans \cite{carliu}, les auteurs ont posé dans le \S 5 un certain nombre
de questions sur la possibilité de relever en caractéristique nulle des
représentations galoisiennes de torsion. Typiquement, on se donne une
$\F_p$-représentation $\bar T$ de $G_K$, et on se demande s'il existe un
réseau $T$ dans une représentation cristalline (resp. semi-stable) à
poids de Hodge-Tate contraints, et un morphisme surjectif $T \to \bar T$.
La théorie des $(\varphi,\tau)$-modules nous semble être une approche
intéressante pour aborder ce type de problèmes (au moins dans le cas des
représentations semi-stables) puisqu'elle permet à la fois de décrire
les représentations libres et de torsion, et de reconnaître de façon
particulièrement simple les réseaux dans les représentations
semi-stables.

\end{document}